\documentclass[10pt,a4paper]{article}
\usepackage{amsopn}

\usepackage{graphics}
\usepackage{graphicx}
\usepackage{epsfig}

\usepackage{amssymb}
\usepackage{amsthm}
\usepackage{amsmath}
\usepackage{amsfonts}
\usepackage{lipsum}
\usepackage[left=2.5cm,right=2.5cm,bottom=3.5cm,top=3.5cm]{geometry}
\usepackage{empheq}
\usepackage{epstopdf}
\usepackage{changepage}
\usepackage{rotating}
\usepackage{adjustbox}
\usepackage{graphbox}
\usepackage{color}
\usepackage{xcolor}
\usepackage{dsfont}
\usepackage{dutchcal}
\usepackage{hyperref} 
\usepackage{float}
\usepackage{resizegather}
\usepackage{multicol}
\usepackage{booktabs} 
\usepackage{colortbl} 
\usepackage{multirow} 
\usepackage{longtable,tabu} 
\usepackage{hhline}   
\usepackage{anyfontsize} 
\usepackage{everysel}
\usepackage{psfrag}
\usepackage{titlesec}
\usepackage{bbm}
\usepackage{bm}
\usepackage{pgf,tikz}

\DeclareMathOperator{\dive}{\nabla \cdot}

\DeclareMathOperator{\grae}{\nabla}

\newcommand{\x}{\mathbf{x}}

\newcommand{\q}{\mathbf{q}}
\newcommand{\n}{\mathbf{n}}

\newcommand{\Flux}{\mathbf{\mathcal{F}}}
\newcommand{\vels}{v}
\newcommand{\press}{\eta}
\newcommand{\vel}{\mathbf{v}}

\definecolor{nverde}{RGB}{0,61,0} 
\definecolor{cr1}{RGB}{200,0,0}
\definecolor{cr2}{RGB}{0,0,200}
\definecolor{cr12}{RGB}{100,0,100}

\newtheorem{theorem}{Proposition}
\newtheorem{weakproblem}{Weak problem}
\theoremstyle{remark}

\hypersetup{
	pdftitle={A staggered semi-implicit hybrid finite volume / finite element scheme for the shallow water equations at all Froude numbers},    
	pdfauthor={S. Busto, M. Dumbser},     
	pdfcreator={S. Busto, M. Dumbser},   
	colorlinks=true,       
}
\title{A staggered semi-implicit hybrid finite volume / finite element scheme \\ for the shallow water equations at all Froude numbers} 

\allowdisplaybreaks

\begin{document}

\begin{center}
	\textbf{ \Large{A staggered semi-implicit hybrid finite volume / finite element  \\ scheme for the shallow water equations at all Froude numbers} }
	
	\vspace{0.5cm}
	{S. Busto\footnote{saray.busto@unitn.it}, M. Dumbser\footnote{michael.dumbser@unitn.it}}
	
	\vspace{0.2cm}
	{\small		
		\textit{$^{(1,2)}$ Laboratory of Applied Mathematics, DICAM, University of Trento, via Mesiano 77, 38123 Trento, Italy}
		
		\textit{$^{(1)}$ Departamento de Matem\'atica Aplicada a la Ingenier\'ia Industrial, Universidad Polit\'ecnica de Madrid, Jos\'e Gutierrez Abascal~2, 28006, Madrid, Spain}
	}
\end{center}

\hrule

\thispagestyle{plain}
\begin{center}
	\textbf{Abstract}
\end{center} 

	We present a novel staggered semi-implicit hybrid finite volume / finite element method for the numerical solution of the shallow water equations at all Froude numbers on  unstructured meshes.
	A semi-discretization in time of the conservative Saint-Venant equations with bottom friction terms leads to its decomposition into a first order hyperbolic subsystem containing the nonlinear convective term and a second order wave equation for the pressure. For the spatial discretization of the free surface elevation an unstructured mesh composed of triangular simplex elements is considered, whereas a dual grid of the edge-type is employed for the computation of the depth-averaged momentum vector. The first stage of the proposed algorithm consists in the solution of the nonlinear convective subsystem using an explicit Godunov-type finite volume method on the staggered dual grid. Next, a classical continuous finite element scheme provides the free surface elevation at the vertices of the primal simplex mesh. The semi-implicit strategy followed in this paper circumvents the contribution of the surface wave celerity to the CFL-type time step restriction, hence making the proposed algorithm well-suited for the solution of low Froude number flows. {It can be shown that in the low Froude number limit the proposed algorithm reduces to a semi-implicit hybrid FV/FE projection method for the incompressible Navier-Stokes equations.} At the same time, the conservative formulation of the governing equations also allows the discretization of high Froude number flows with shock waves. As such, the new hybrid FV/FE scheme can be considered an all-Froude number solver, able to deal simultaneously with both, subcritical as well as supercritical flows. Besides, the algorithm is well balanced by construction. The accuracy of the overall methodology is studied numerically and the C-property is proven theoretically and is also validated via numerical experiments. The numerical solution of several Riemann problems, including flat bottom and a bottom with jumps, attests the robustness of the new method to deal also with flows containing bores and discontinuities. Finally, a 3D dam break problem  over a dry bottom is studied and our numerical results are successfully compared with numerical reference solutions obtained for the 3D free surface Navier-Stokes equations and with available experimental reference data.
	
\vspace{0.2cm}

\vspace{0.4cm}

\hrule

\section{Introduction}\label{sec:intro}
Free surface shallow water flows are ubiquitous in geophysics as well as in environmental and coastal engineering. The most widespread model for the mathematical description of such flows are the shallow water or Saint-Venant equations. The basic assumptions of this model are a hydrostatic pressure distribution and a constant velocity profile in the vertical direction. Furthermore, the free surface is assumed to be a single-valued function. There are two main classes of well-established numerical methods for the discretization of the shallow water equations: the family of explicit Godunov-type finite volume schemes, as well as the family of semi-implicit finite volume / finite difference methods. For explicit Godunov-type finite volume schemes, the so-called C-property first found by Berm{\'u}dez and V{\'a}zquez-Cend{\'o}n in \cite{BVC94} is of fundamental importance for any numerical discretization of the shallow water equations in order to guarantee that lake-at-rest equilibrium solutions of the system are exactly preserved at the discrete level for arbitrary bottom topography. Achieving this well-balanced property in the context of Godunov-type schemes is non-trivial and a lot of research has been carried out in this direction in the last decades. Recent developments concerning well-balanced Godunov-type finite volume schemes for the Saint-Venant equations can be found, for example, in  {\cite{GL96wb,VC99,GNVC00,LeVeque98,Gos00,DT11,Ric11,Noelle1,Noelle2,GNVC00,GarciaNavarro2,Bouchut,Bollermann,CL21,SGNNP21,DBC16,AAHKW19,KLZ20,ECS21,CKLLW21}} 
{while modern work on multilayer models is presented in \cite{DL13,GCM20,Liu21}}. A review of {classical} shock-capturing finite volume schemes for shallow water flows is available in the textbook \cite{toro-book-swe}. A modern and very elegant concept for the construction of well-balanced Godunov-type finite volume schemes is based on the framework of path-conservative schemes introduced by Par\'es and Castro in  {\cite{Pares2004,Par06,CGP06,Castro2d,Castro2007,Castro2008,CastroPardoPares,CLP13}}, based on the definition of weak sol0utions for non-conservative hyperbolic PDE introduced in \cite{DLMtheory}. High order explicit discontinuous Galerkin finite element schemes for shallow water flows have been presented in {\cite{Rhebergen2008,ADERNC,ShuDGSWE,GassnerSWE,WWGW18}}. 

While very efficient and accurate for critical and supercritical flows including shock waves, explicit Godunov-type finite volume schemes have the substantial drawback that they become computationally inefficient and inaccurate for low Froude number flows, since the time step is limited by a CFL-type stability condition based on the surface wave celerity. For low Froude number flows, the use of semi-implicit schemes is preferable from an efficiency and accuracy point of view. An overview of the state of the art about semi-implicit schemes for free surface flows can be found, for example, in {\cite{Casulli1990,CasulliCheng1992,CasulliCattani,Casulli1999,CasulliWalters2000,Casulli2009,CasulliStelling2011,BrugnanoCasulli,BrugnanoCasulli2,BrugnanoSestini,CasulliZanolli2012,KramerStelling,GFN21,FD21,BTP21}}. Higher order extensions of staggered semi-implicit finite volume schemes are presented, e.g., in \cite{BDR13,BPR19}, while high order semi-implicit discontinuous Galerkin finite element methods for the solution of the shallow water equations were proposed in {\cite{GR10,TBR13,DC13,TD14sw,IoriattiSWE,SFT20}}. {An implicit-explicit approach has also been recently used in the context of multilayer shallow water equations, see e.g. \cite{ABPS10,BFGN18,GB21}.
}

{
Another important aspect that needs our attention when developing numerical schemes for low Froude number flows is the asymptotic preserving property (AP). That is, we should check that the proposed method 
provides a consistent approximation of the limiting continuous equations when Fr$\rightarrow 0$.
Let us note that the first AP schemes have been originally proposed in the context of kinetic energy equations in the late ninety's, see \cite{Jin95,JPT98,Jin99} for the pioneering works and \cite{CKLS19,EKKT21} for more recent advances. Since then, AP schemes have been developed for a wide variety of PDE systems including hyperbolic conservation laws, \cite{BT11,BCT13,BLR14,BCT21}. Moreover, an increasing number of authors has focused on the development of AP schemes for Euler and Navier-Stokes equations, see \cite{PM05,NBALM12,CordierDegond,DC16,BLY17,CDN17,BRD18,BQRX19,AbateIolloPuppo,AS21,TPK20} and references therein. 
This research has set the basis for novel developments of AP schemes for the shallow water equations. Several works already presented in this context are, e.g., \cite{BALN14,CDV17,KLL19,Liu20,BBBB21,LCK19}. 
}

The above lists of references do not pretend to be complete, since the research topic is extremely vast and many relevant contributions have been made in the last decades. 

The mathematical modelling and the numerical simulation of non-hydrostatic and turbulent free surface flows is more complex than the solution of the classical Saint-Venant equations. In this context the reader is referred to the work on Bousinesq-type dispersive models  \cite{madsen:1991,madsen:1992,JSMarie,escalanteNH} and their recent hyperbolic reformulations in \cite{favrie:2017,Escalante2018,ricchiutoHyp,DispersiveSWE,Busto2021}. Fully three-dimensional models that allow even a multi-valued free surface profile have been discussed in \cite{CasulliVOF,DIM3D,SPH3D,levelset3Da,levelset3Db,Kleefsman,Rieber,Loehner,Colagrossi,Monaghan1994} and references therein, while recent hyperbolic models for turbulent shallow water flows have been forwarded and studied in \cite{Gavrilyuk2018,Ivanova2019,Bhole2019,SWTurbulence}. 

Most of the numerical methods mentioned so far are based on finite volume or a finite difference discretization of the governing equation. Classical finite element schemes for shallow water flows have been, for example, discussed in \cite{WaltersCasulli1998} and \cite{Telemac1991,Telemac2011}. 

Usually, there is a clear separation between finite volume-type schemes and continuous finite element methods for computational fluid dynamics. In a recent series of papers, a novel semi-implicit hybrid finite volume / finite element scheme on staggered unstructured meshes has been developed for the solution of the incompressible \cite{BFSV14,BFTVC17} and the compressible \cite{BBDFSVC20,BRVD21} Navier-Stokes equations. The idea in this family of schemes is to apply first a semi-discretization to the governing equations in time only. The nonlinear convective terms are then discretized with an explicit finite volume scheme on an edge-based / face-based staggered dual mesh, while the pressure terms are discretized implicitly at the aid of a classical continuous finite element scheme on the primal simplex mesh. In the present paper, we extend this methodology to the shallow water equations. Being a semi-implicit scheme written in terms of the conservative variables, the proposed hybrid FV/FE solver is able to deal efficiently with flows at \textit{all Froude numbers}, ranging from potential type flows to supercritical flows with shock waves. In the numerical examples section of this paper we actually show flows ranging from Froude numbers of the order {Fr$=10^{-7}$} up to Fr$=3$. Furthermore, since the pressure is discretized at the aid of a classical continuous finite element method, there is no orthogonality requirement for the mesh, in contrast to the method forwarded in  \cite{CasulliWalters2000}.

The rest of the paper is organized as follows. In Section \ref{sec:goveq} we recall the shallow water system and introduce the main related notations. Section \ref{sec:method} presents the semi-implicit hybrid finite volume / finite element algorithm. First, we introduce the time discretization of the equations. Next, the staggered mesh arrangement is detailed and the spatial discretization is outlined. The method is then validated in Section \ref{sec:numericalresults}, including a comparison with exact analytical solutions and experimental results. Some concluding remarks and notes on future developments are drafted in Section \ref{sec:conclusions}. 

\section{Governing partial differential equations} \label{sec:goveq}
\begin{figure}[ht!]
	\centering
	\includegraphics[width=0.5\linewidth]{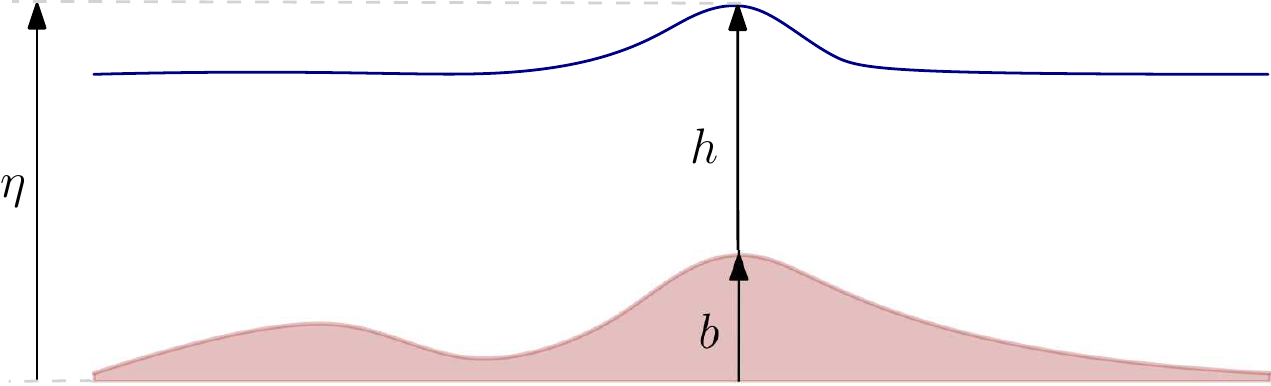}
	\caption{Sketch of the computational domain and of the notation used in this paper.}
	\label{fig:swenotation}
\end{figure}

The shallow water equations including bottom friction effects can be written as  
\begin{eqnarray}
	\label{eqn.eta} 
	\frac{\partial \eta}{\partial t }       + \nabla \cdot \q &=& 0, \\ 
	\label{eqn.q} 
	\frac{\partial \q}{\partial t } + \nabla \cdot \left( \vel \otimes \q \right) + g h \, \nabla \eta &=& -\gamma \q, \\ 
	\label{eqn.bottom}
	{\frac{\partial b}{\partial t}  } & {=} & {0}, 
\end{eqnarray}
with {$\mathbf{x}=(x,y) \in \Omega \subset \mathbb{R}^2$ the spatial coordinate vector, $t\in \mathbb{R}_0^+$ the time}, $h=h(\x,t)$ the water depth, $b=b(\x)$ the bottom elevation (constant in time), $\eta = h + b$ the free surface elevation, $\vel(\x,t)=(u,v)$   the velocity field, $\q = h \vel$ the corresponding conservative variable and $g$ the gravity constant. The short hand notation $\gamma$ in the friction term reads
\begin{equation}
	\gamma = \frac{g}{\kappa_{s}^{2} h^{\frac{4}{3}}}\left|\vel\right|, 
\end{equation}
where $k_{s}$ is the Strickler coefficient that can be rewritten in terms of the Manning coefficient $n$ as $\kappa_{s}=\frac{1}{n}$. For a sketch of the problem and the notation used in this paper, see Figure \ref{fig:swenotation}. 

{
The eigenvalues of the full shallow water system \eqref{eqn.eta}-\eqref{eqn.bottom} in $x$ direction are 
\begin{equation}
\lambda_1 = u - \sqrt{g h}, \quad 
\lambda_{2} = 0, \quad 
\lambda_{3} = u, \quad 
\lambda_4 = u + \sqrt{g h},  
\label{eqn.lambda.full} 
\end{equation}
where the zero eigenvalue $\lambda_2=0$ corresponds to the bottom jump, the eigenvalue $\lambda_3=u$ to the transport of the transverse velocity (shear wave) and the outer eigenvalues $\lambda_{1,3}=u\mp \sqrt{gh}$ correspond to the propagation of free surface gravity waves. }

{ 
Throughout this paper we will make use of the following \textit{splitting} of the shallow water system into a \textit{convective subsystem}, which will be discretized \textit{explicitly} in time, and a \textit{pressure subsystem}, which will be discretized \textit{implicitly} in time. Such a splitting is similar to the flux-vector splitting of the compressible Euler equations proposed by Toro and V\'azquez-Cend\'on in \cite{TV12} and corresponds to the usual splitting used in the context of semi-implicit methods for the shallow water equations, see \cite{Casulli1990,CasulliCheng1992,CasulliCattani}: 
\paragraph{Convective subsystem}
The convective subsystem, which will be discretized \textit{explicitly} in this paper and which will thus be responsible for the CFL time step restriction of the final semi-implicit scheme, reads  
\begin{eqnarray}
\label{eqn.eta.c} 
\frac{\partial \eta}{\partial t }        &=& 0, \\ 
\label{eqn.q.c} 
\frac{\partial \q}{\partial t } + \nabla \cdot \left( \vel \otimes \q \right)  &=& 0, \\ 
\label{eqn.bottom.c}
{\frac{\partial b}{\partial t}  } & {=} & {0}. 
\end{eqnarray}
The corresponding eigenvalues of the convective subsystem \eqref{eqn.eta.c}-\eqref{eqn.bottom.c} in $x$ direction are 
\begin{equation}
\lambda^c_1 = 0, \quad 
\lambda^c_{2} = 0, \quad 
\lambda^c_{3} = u, \quad 
\lambda^c_4 = 2u.  
\label{eqn.lambda.c} 
\end{equation}
One immediately notes the absence of the surface wavespeed $c=\sqrt{gh}$ in equation \eqref{eqn.lambda.c} compared to the eigenvalues of the full system shown in \eqref{eqn.lambda.full}. As such, the method proposed in this paper will be particularly suitable for \textit{low Froude} number flows. 
}

{
\paragraph{Pressure subsystem}
The pressure subsystem, which will be discretized \textit{implicitly} in this paper and which therefore does \textit{not} contribute to the CFL time step restriction is given by 
\begin{eqnarray}
\label{eqn.eta.p} 
\frac{\partial \eta}{\partial t }       + \nabla \cdot \q &=& 0, \\ 
\label{eqn.q.p} 
\frac{\partial \q}{\partial t } + g h \, \nabla \eta &=& 0, \\ 
\label{eqn.bottom.p}
{\frac{\partial b}{\partial t}  } & {=} & {0}. 
\end{eqnarray}
The eigenvalues of the pressure subsystem \eqref{eqn.eta.p}-\eqref{eqn.bottom.p} in $x$ direction  are 
\begin{equation}
\lambda^p_1 = -\sqrt{g h}, \quad 
\lambda^p_{2} = 0, \quad 
\lambda^p_{3} = 0, \quad 
\lambda^p_4 = +\sqrt{g h},   
\label{eqn.lambda.p} 
\end{equation}
similar to those found in \cite{TV12} for the flux vector splitting of the compressible Euler equations.
}

\section{The hybrid FV/FE method}\label{sec:method}
As already mentioned in the introduction, the hybrid FV/FE methodology presented in this paper is based on the idea of a judicious combination of finite volume (FV) and continuous finite element (FE) methods to efficiently discretize nonlinear hyperbolic systems of PDEs. To this end, the first step to be performed is a semi-discretization in time that leads to the split of the equations into a set of nonlinear first order hyperbolic transport equations and a second order wave-equation for the pressure, {following the ideas of the flux splitting discussed in the previous section}. Then, the spatial derivatives are discretized employing  edge-based staggered unstructured meshes. Explicit Godunov-type finite volume methods are  employed for the spatial discretization of the nonlinear convective part of the equations, obtaining an intermediate approximation of the conservative momentum in the staggered dual grid. Next, a second order wave equation for the pressure is solved using classical continuous finite elements. Its solution provides the necessary correction for the intermediate solution of the momentum obtained in the previous step. 
In particular, discretization only in time of system \eqref{eqn.eta}-\eqref{eqn.q}, while keeping space still continuous, yields the following semi-discrete scheme: 
\begin{eqnarray}
	\label{eqn.eta.aux} 
	\eta^{n+1} &=& \eta^n  - \Delta t \,  \nabla \cdot \q^{n+1}, \\ 
	\label{eqn.q.aux} 
	\q^{n+1} &=& \q^n - \Delta t \, \nabla \cdot \left( \vel^n \otimes \q^n \right) - \Delta t \, g \, h^n \, \nabla \eta^{n+1} -  \Delta t \, \gamma^{n} \q^{n+1},   
\end{eqnarray}
where $q^{n}$, $\vel^{n}$, $\eta^n$ are the known discrete solutions at time $t^{n}$, and $\q^{n+1}$, $\eta^{n+1}$ are the approximations of $\q(\x,t^{n+1})$, $\eta(\x,t^{n+1})$.
In the following, we will use the abbreviation 
\begin{equation}
	\label{eqn.qstar.sd} 
	\q^* = \q^n - \Delta t \, \nabla \cdot \left( \vel^n \otimes \q^n \right)
\end{equation}

for the intermediate approximation of $\q$. Hence the above system can be rewritten more compactly as    
\begin{eqnarray}
	\label{eqn.eta.sd} 
	\eta^{n+1} &=& \eta^n  - \Delta t \,  \nabla \cdot \q^{n+1}, \\ 
	\label{eqn.q.sd} 
	\q^{n+1} &=& \q^* - \Delta t \, g \, h^n \, \nabla \eta^{n+1} - \Delta t \, \gamma^{n} \q^{n+1}.   
\end{eqnarray}
Inserting the semi-discrete momentum equation \eqref{eqn.q.sd} into the semi-discrete mass conservation equation (\ref{eqn.eta.sd}) yields the following semi-discrete wave-equation for the unknown free surface elevation $\eta^{n+1}$: 
\begin{equation} 
	\label{eqn.wave.eta} 
	\eta^{n+1} - \Delta t^2 \,  g \, \nabla \cdot \left(\frac{h^n}{1+\Delta t\, \gamma^{n}} \nabla \eta^{n+1} \right) = \eta^n  - \Delta t \,  \nabla \cdot \left( \frac{1}{1+\Delta t\, \gamma^{n}}\, \q^{*}\right) .  
\end{equation} 	
Since $\frac{h^n}{1+\Delta t \gamma^{n}} \geq0$ the operator on the left hand side is symmetric and positive definite. 
Equation \eqref{eqn.wave.eta} can now be discretized with a classical continuous finite element method for $\eta^{n+1}$ on the unstructured primal triangular mesh, while the nonlinear advection terms contained in $\q^{*}$ are discretized conveniently on the dual grid via an explicit Godunov-type finite volume approach, see \cite{toro-book-swe} for details. Once the new free surface location $\eta^{n+1}$ is known, the quantity $\q^{n+1}$ can be readily computed via (\ref{eqn.q.sd}) in a post-projection stage, i.e. 
\begin{equation}
	\label{eqn.q.update} 
	\mathbf{q}^{n+1} = \frac{1}{1+\Delta t\, \gamma^{n}} \left( \mathbf{q}^* - \Delta t \, g \, h^n \, \nabla \eta^{n+1}\right) .   
\end{equation}
The algorithm fits perfectly well into the existing general framework of hybrid FV/FE methods established in \cite{BFSV14,BFTVC17,BBDFSVC20,BRVD21} for the incompressible and compressible Navier-Stokes equations. 
Besides, an important point to be checked when dealing with shallow water equations is the well balanced property in the sense of \cite{BVC94,GL96wb,LeVeque98,BDDV98,GNVC00,Ric11}, i.e. the capability of the scheme to preserve steady state solutions of the lake-at-rest type, i.e. solutions with flat free surface $\eta=const$, vanishing velocity $\mathbf{v} = 0$ and arbitrary bottom topography $b(\x)$. For the proposed scheme this property is fulfilled by construction:  

\begin{theorem}[C-property]
	The scheme \eqref{eqn.qstar.sd}, \eqref{eqn.wave.eta}, \eqref{eqn.q.update}  is well balanced.
\end{theorem}
\begin{proof}
	Let us assume an initial condition of the form $\eta^n=\eta_0= \textnormal{const}.$ and $\vel^n=0$, then we have $\q^n=0$ and thus also $\q^*=0$, even with a Rusanov-type flux in the discretization of the nonlinear convective terms, since any numerical dissipation in the numerical fluxes associated with the advective term $\nabla \cdot (\vel^n \otimes \q^n)$ will vanish for $\q^n=0$. With $\q^*=0$ it is then obvious that $\eta^{n+1}=\eta^{n}=\eta_0$ is solution of \eqref{eqn.wave.eta} and since the discrete gradient of a constant is zero, we also have $\q^{n+1}=0$ due to \eqref{eqn.q.update}. Therefore, $\eta=\eta_0 = \textnormal{const}.$ and $\q=\vel^n=0$, which corresponds to the well-known C-property for lake-at-rest solutions of the original system \eqref{eqn.eta}-\eqref{eqn.q}.
\end{proof}

\subsection{Staggered unstructured grids} \label{sec:numdisc}

\begin{figure}
	\centering
	\includegraphics[width=0.32\linewidth]{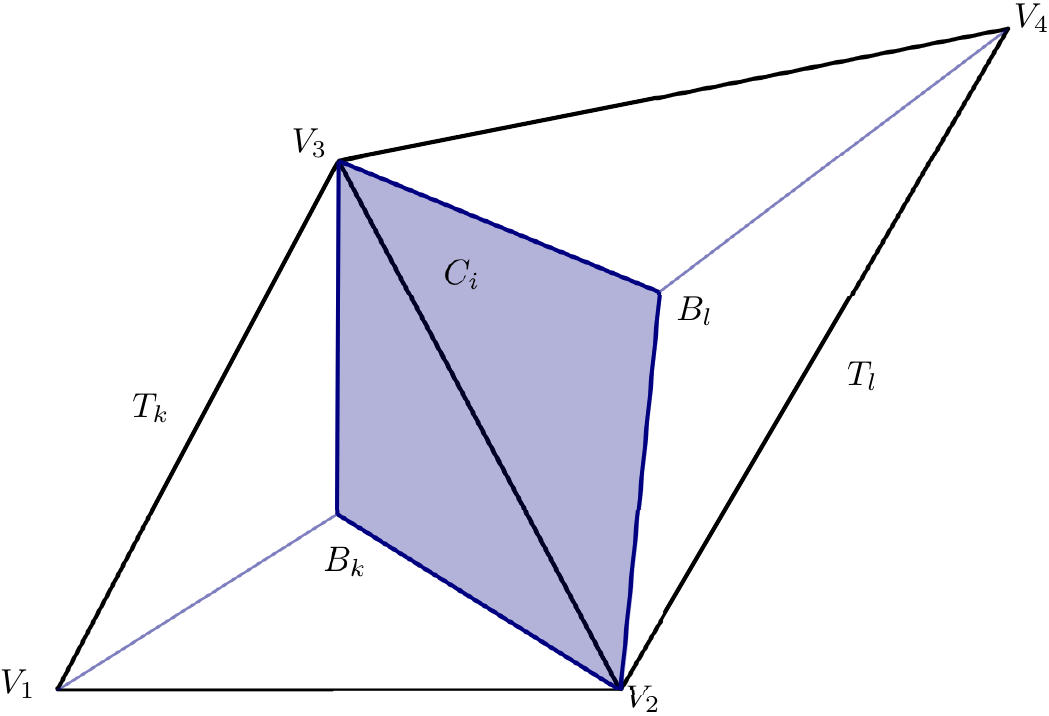}\hfill
	\includegraphics[width=0.32\linewidth]{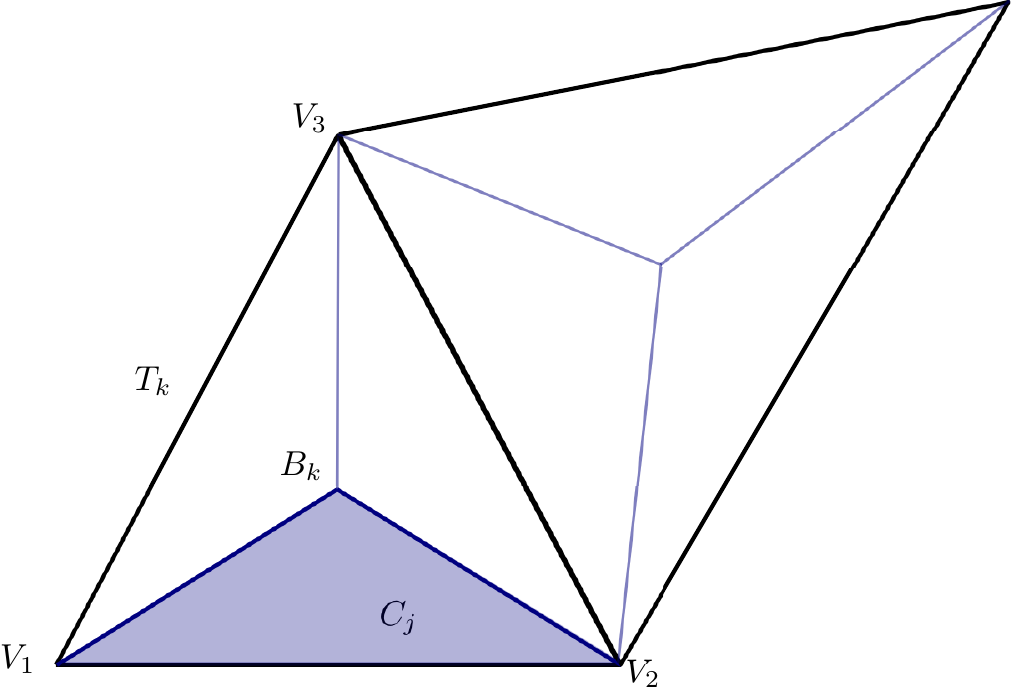}\hfill
	\includegraphics[width=0.32\linewidth]{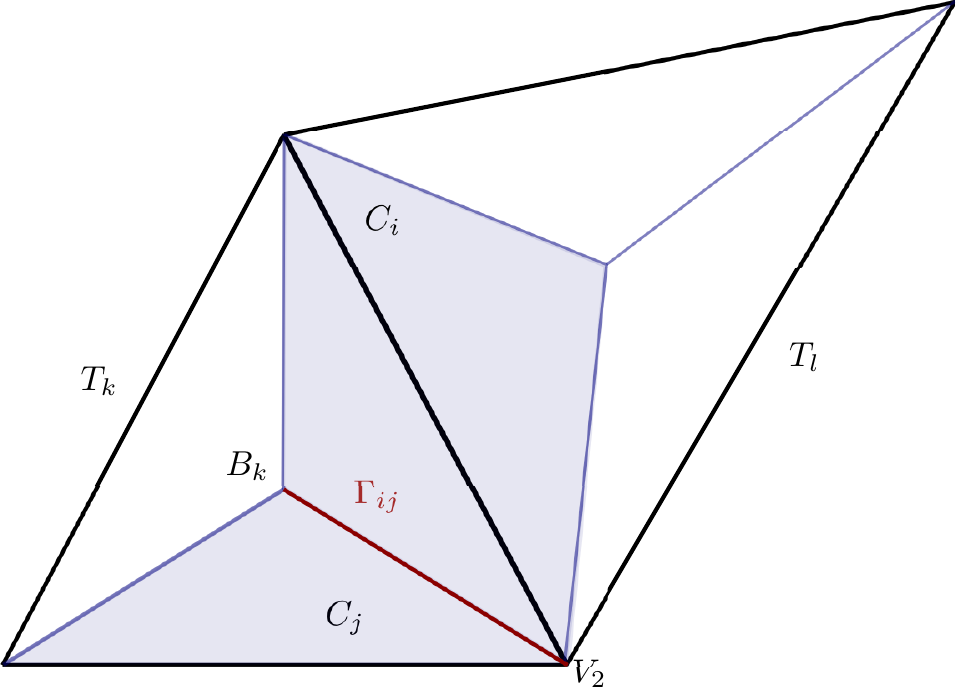}
	\caption{Construction of the dual staggered mesh from the primal triangular grid. Left: interior dual cell $C_{i}$ (shadowed in blue). Centre: boundary dual cell $C_{i}$ (shadowed in blue). Right: boundary edge of the dual elements $\Gamma_{ij}$ (red).}
	\label{fig.dualinterior}
\end{figure}

Before detailing the proposed spatial discretization, we will introduce some mesh related notations. 
Let us consider a spatial domain $\Omega\in\mathbb{R}^{2}$ covered up with a set of non overlapping triangles $T_{k}$, $k=1,\dots, nel$, $nel$ denoting the total number of elements. For each interior edge, $\Gamma_{kl}$, shared by two primal elements $T_{k},\,T_{l}$, we can construct a dual element of the edge type, $C_{i}$, by gathering the two triangles determined by the vertex of the edge and the barycentres of  $T_{k}$ and $T_{l}$, see Figure \ref{fig.dualinterior}. Accordingly, a dual cell at the domain boundary is the triangle defined by the two vertexes of the boundary edge,  $\Gamma_{k_{b}}\subset \partial \Omega$, and the barycentre of the primal element containing it. We will denote by $N_{i}$ the node of the dual elements $C_{i}$, $\left|C_{i}\right|$ the area of the cell and $\mathcal{K}_{i}$ the set of indexes identifying the neighbouring cells that have a common edge with $C_{i}$, being $\Gamma_{ij}$ the edge with $\tilde{\n}_{ij}$ its outward unit normal vector and $\Gamma=\cup_{j\in \mathcal{K}_{i}} \Gamma_{ij}$ the boundary of the dual element. Moreover, we define $\n_{ij}:=\tilde{\n}_{ij}\left\|\n_{ij}\right\|$ the weighted normal vector where $\left\|\n_{ij}\right\|$ is the length of edge $\Gamma_{ij}$.

\subsection{Transport stage}
The first stage of the algorithm involves the spatial discretization of equation \eqref{eqn.qstar.sd}  following a finite volume approach. Integration on a spatial control volume $C_{i}$ and use of the Gauss theorem gives
\begin{eqnarray}
	\label{eqn.qstar.fd} 
	\q^*_{i} = \q^n_{i} - \frac{\Delta t}{\left|C_{i}\right|}\int_{\Gamma_{i}} \Flux \left( h^n, \q^n \right) \tilde{\n}_{i}\, dS, \qquad \Flux \left( h^n, \q^n \right) := \frac{1}{h^{n}} \, \q^{n}\otimes \q^{n} = \vel^{n}\otimes \q^{n}.
\end{eqnarray}
Denoting by 
\begin{equation}
	{\mathcal{Z} }\left(h^{n}, \q^{n} ,  \tilde{\n}_{i}\right)  := {\Flux} (h^{n}, \q^{n}) \, {\tilde{\n}}_{i}
\end{equation}
the global normal flux through the boundary of the cell and splitting the integral on the cell boundary into the sum of the integrals on each edge, we get
\begin{equation}
	\int_{ \Gamma_{i} } \Flux(
	h^{n}, \q^{n})  \tilde{\n}_{i} \, dS
	= \displaystyle \sum_{j \in {\cal K}_{i}} \displaystyle
	\int_{\Gamma_{ij}} {
		\mathcal{Z} }\left( h^{n}, \q^{n}, \tilde{\n}_{ij} \right)   \,\mathrm{dS},
\end{equation}
that can be approximated using the modified Rusanov flux function, \cite{Rus62,BBDFSVC20,BRVD21},
\begin{gather}	
	\boldsymbol{\phi} \left( h_{i}^{n}\, ,h_{j}^{n}\, , \q_{i}^{n}\, ,\q_{j}^{n}\, , {\n}_{ij} \right) = \displaystyle\frac{1}{2} \left( {\cal Z}
	(h_{i}^n\, ,\q_{i}^n\, ,\n_{ij})+{\cal Z}(h_{j}^n\, , \q_{j}^n \, , \n_{ij})\right) -\frac{1}{2}
	\alpha^{n}_{RS,\, ij} \left( \q_{j}^n-\q_{i}^n\right) , \notag\\
	\alpha^{n}_{RS,\, ij} = \alpha_{RS}(h_{i}^n,h_{j}^n,\q_{i}^n,\q_{j}^n,\n_{ij}):=\max \left\lbrace 2\left|\vel_i^{n}\cdot\n_{ij}\right|,2\left|\vel_j^{n}\cdot\n_{ij}\right|\right\rbrace + c_{\alpha} \left\|\boldsymbol{\eta}_{ij} \right\|, 
	\label{eq.flux}
\end{gather}
where an artificial viscosity coefficient, $c_{\alpha}\in\mathbb{R}^{+}_{0}$, has been added to the classical computation of the maximum signal speed on the edge. This coefficient will be activated only on particular tests to ameliorate the stability properties of the scheme in the presence of large water depth jumps and small velocities.

To increase the accuracy of the proposed scheme a local ADER methodology is employed, \cite{BFTVC17}. The first step consists in the computation of the ENO-based reconstruction polynomials of a variable $w$ at the neighbouring of each edge,
	\begin{equation}
		P^{i}_{ij}(N)=w_i+(N-N_{i})\left( \grae\, w\right)^{i}_{ij},\quad P^{j}_{ij}(N)=w_j+(N-N_{j})\left( \grae\, w\right)^{j}_{ij},
	\end{equation}	
	which is achieved by profiting from the dual mesh arrangement to obtain a reduced stencil. First, the gradients are computed at each primal element using the Crouzeix-Raviart finite elements basis. Next, at each edge we select the gradient as
	\begin{gather*}
		\left( \grae\, w\right)^{i}_{ij}= \left\lbrace
		\begin{array}{lc}
			\left(\grae w \right)_{T_{ijL}}, & \textrm{if }\left| \left(\grae w \right)_{ T_{ijL}}\cdot \left(N_{ij}-N_{i}\right)\right| \leq \left|\left(\grae w \right)_{ T_{ij}}\cdot \left(N_{ij}-N_{i}\right) \right|,\\[8pt]
			\left(\grae w \right)_{ T_{ij}}, & \mathrm{otherwise};
		\end{array}
		\right.
	\end{gather*}
	\begin{gather*}
		\left( \grae w \right)^{j}_{ij}= \left\lbrace
		\begin{array}{lc}
			\left(\grae w\right)_{ T_{ijR}}, & \textrm{if }\left|\left( \grae w \right)_{ T_{ijR}}\cdot \left(N_{ij}-N_{j}\right)\right| \leq \left|\left(\grae w \right)_{ T_{ij}}\cdot \left(N_{ij}-N_{j}\right)\right|,\\[8pt]
			\left(\grae w\right)_{ T_{ij}},  & \mathrm{otherwise},
		\end{array}
		\right.
	\end{gather*}
	This procedure introduces the non linearity needed to circumvent Godunov's theorem obtaining an essentially non-oscillatory second order scheme. An alternative approach would be the use of classical slope limiters, such as the minmod limiter of Roe, \cite{Roe85,Toro}, or the more restrictive Barth and Jespersen limiter, see \cite{BJ89}. Further strategies to be considered in future work include  \textit{a posteriori} limiting techniques like the MOOD methodology, \cite{CDL11}.
	Evaluation of the reconstruction polynomials at the barycentre of each edge, $\x_{N_{ij}}$ provides the boundary extrapolated data as
	\begin{gather} w_{i\, N_{ij}} =p^{i}_{ij}(N_{ij})
		= w_{i}+(N_{ij}-N_{i})   \left( \grae w\right)^{i}_{ij},\qquad
		w_{j\, N_{ij}}=p^{j}_{ij}(N_{ij})=
		w_{j}+(N_{ij}-N_{j})   \left( \grae w\right)^{j}_{ij}.  
	\end{gather}
	Next, we approximate the half in time approximation of conservative variables by using the mid point rule combined with a Cauchy-Kovalevskaya procedure applied to equation \eqref{eqn.qstar.sd},
	\begin{equation}
		\overline{\q}_{i\, N_{ij}}=\q_{i\, N_{ij}} + \q_{ N_{ij}}^{\ast},\qquad \overline{\q}_{j\, N_{ij}}= \q_{j\, N_{ij}} + \q_{N_{ij}}^{\ast},\label{eq:Wevolij}
	\end{equation}
	where	
	\begin{equation}
		\q_{N_{ij}}^{\ast}:=-\frac{\Delta t}{2\mathcal{L}_{ij}}\left[ \mathcal{Z}(h_{i\, N_{ij}}^{n}, \q_{i\, N_{ij}}^{n},\n_{ij}  )+\mathcal{Z}(h_{j\, N_{ij}}^{n}, \q_{j\, N_{ij}}^{n},\n_{ij}  )\right].
		\label{eq:wevol}
	\end{equation}
	Finally, the obtained approximations are substituted in \eqref{eq.flux} yielding 
	\begin{equation}
		\boldsymbol{\phi} \left(h_{i\, N_{ij}}^{n}\, ,h_{j\, N_{ij}}^{n}\, ,  \overline{\q}_{i\, N_{ij}}\, ,\overline{\q}_{j\, N_{ij}}\, , {\n}_{ij} \right) = \displaystyle\frac{1}{2} \left( {\cal Z}
		(h_{i\, N_{ij}}^n\, ,\overline{\q}_{i\, N_{ij}}\, ,\n_{ij})+{\cal Z}(h_{j\, N_{ij}}^n\, , \overline{\q}_{j\, N_{ij}} \, , \n_{ij})\right) -\frac{1}{2}
		\alpha^{n}_{RS,\, ij} \left(\overline{\q}_{j\, N_{ij}}-\overline{\q}_{i\, N_{ij}}\right).
	\end{equation}
For further details on ADER methodology including its theoretical analysis, modified approaches and applications to different PDE systems, we may refer the reader to \cite{TMN01,Toro,DET08,BD14,BTVC16,BCDGP20}.

\subsection{Pre-projection stage}
The intermediate solution of the conservative variables $\q^{\ast}$ is initially obtained within the transport stage as an average value at each dual element. However, the value of $\q^{\ast}$ takes part on the right hand side of equation \eqref{eqn.wave.eta}, which is solved using finite elements in the primal grid. Consequently, we need to interpolate the solution between both meshes. In particular, we will approximate $\q^{\ast}$ at each dual element $T_{k}$ by adopting a point interpretation of the computed data and considering the average of its values at the dual volumes associated to each face:
\begin{equation}
	\q^{\ast}_{k} = \frac{1}{3}\sum_{i\in\mathcal{K}_{k}} \q^{\ast}_{i},
\end{equation}
where $\mathcal{K}_{k}$ is the set of indexes pointing to the faces of the primal element $T_{k}$. Let us note that the nodes of the dual elements are assumed to be located at the barycenter of the primal faces following a finite differences approach, which allows for an exact interpolation of data between the staggered grids.

\subsection{Projection stage}\label{sec.projection}
Within the projection stage the solution of \eqref{eqn.wave.eta} is obtained using classical continuous finite elements. Let us consider the space $V_{0}:=\left\lbrace z\in H^{1}\left(\Omega\right): \int_{\Omega} z dV = 0\right\rbrace$ and $z\in V_{0}$ a test function. Multiplication of \eqref{eqn.wave.eta} by the test function, integration in $\Omega$, and use of Green's formula yields
\begin{eqnarray}
	\int_{\Omega} \eta^{n+1} z dV 
	+ \Delta t^2 g \int_{\Omega} \frac{h^{n}}{1+\Delta t\, \gamma^{n}}\grae \eta^{n+1}\cdot \grae z dV = 
	\Delta t^2 g \int_{\Gamma} \frac{h^{n}}{1+\Delta t\, \gamma^{n}}\grae \eta^{n+1} \n \, z dA \notag\\
	+ \int_{\Omega} \eta^{n} z dV 
	+ \Delta t \int_{\Omega}  \frac{1}{1+\Delta t\, \gamma^{n}}\, \q^{\ast}\cdot \grae z dV 
	- \Delta t \int_{\Gamma} \frac{1}{1+\Delta t\, \gamma^{n}}\, \q^{\ast}\cdot \n\, z dA.
	\label{eqn.wave.eta.disc}
\end{eqnarray}
Besides, the dot product of \eqref{eqn.q.sd} by the test function multiplied by the normal vector integrated on the boundary of the computational domain, $\Gamma=\partial\Omega$, gives
\begin{equation}\label{eqn.q.sd.discef}
	\int_{\Gamma} \frac{1}{1+\Delta t\, \gamma^{n}}\, \q^{\ast} \cdot \n\, z dA = \int_{\Gamma} \q^{n+1} \cdot \n\, z dA + \Delta t g \int_{\Gamma} \frac{h^{n}}{1+\Delta t\, \gamma^{n}} \grae \eta^{n+1} \cdot \n\, z dA.
\end{equation}
Therefore, the weak problem for the free surface computation, obtained by substituting  \eqref{eqn.q.sd.discef} in \eqref{eqn.wave.eta.disc}, reads
\begin{weakproblem}\label{wp.orig}
	Find the free surface $\eta^{n+1}\in V_{0}$ that satisfies
	\begin{equation}\label{eqn.weakproblem}
		\int_{\Omega} \eta^{n+1} z dV 
		+ \Delta t^2 g \int_{\Omega}  \frac{h^{n}}{1+\Delta t\, \gamma^{n}} \grae \eta^{n+1}\cdot \grae z dV = 
		\int_{\Omega} \eta^{n} z dV 
		+ \Delta t \int_{\Omega} \frac{1}{1+\Delta t\, \gamma^{n}}\, \q^{\ast}\cdot \grae z dV 
		- \Delta t \int_{\Gamma}\q^{n+1}\cdot \n\, z dV
	\end{equation}
	for all $z \in V_{0}$.
\end{weakproblem}
We now discretize the above problem using  $\mathbb{P}_{1}$ finite elements. To this end, we need to approximate 
the water depth appearing in the stiffness matrix. Since the free surface is already computed at the previous time step at the vertex, $V_{k_{m}}, \, m=1,\dots,3$, of each primal cell, $T_{k}$, we can obtain $h^{n}_{k}$ as an average at each element,
\begin{equation}\label{eqn.weak}
	 h^{n}_{k} = \frac{1}{3} \sum_{m=1}^{3}\left( \eta^{n}_{k_m} - b_{k_{m}}\right). 
\end{equation}
Let us note that, to guarantee consistency in the discretization of the different terms in \eqref{eqn.weakproblem}, it is important to approximate the first term in the right hand side of the related system by applying the mass matrix to the solution at the previous time step. On the other hand, the third term in the right hand side of  \eqref{eqn.weakproblem} depends on the boundary conditions defined for the velocity field. If a Dirichlet boundary condition is imposed the exact value of the velocity at the barycenter of each boundary face is set. However, for outflow or symmetry boundary conditions, where the normal velocity is left free, the value at the boundary is not known a priori and needs to be extrapolated from the interior. The resulting system can then be efficiently solved using a conjugate gradient method. 

The stability of the final algorithm can be improved by adding a Rusanov-type dissipation in the left hand side of the free surface elevation system  \eqref{eqn.weakproblem}. Accordingly, the following extra term can be considered: 
\begin{equation}
	 \frac{\Delta t \, |T_{k}|}{|\partial T_{k}|} \int_{T_{k}} \left| \hat{\lambda}_{k}^{n} \right| \grae \phi_l \cdot \grae\phi_r dV\, \hat{\eta}^{n+1}_{k,\, l}, 
\end{equation}
where $\left| \hat{\lambda}_{k}^{n} \right|$ denotes the maximum absolute eigenvalue of the original shallow water system, i.e.
\begin{equation}
	\left| \hat{\lambda}_{k}^{n} \right| = \max \left\lbrace \left| 2 \vel_{k}^{n}\cdot \n -\sqrt{g h_{k}^{n}}\right|, \left|2 \vel_{k}^{n}\cdot \n\right|, 
	\left| 2 \vel_{k}^{n}\cdot \n + \sqrt{g h_{k}^{n}} \right| \right\rbrace ,
\end{equation}
$\left\lbrace\phi_{l}\right\rbrace$ are the basis functions and $\hat{\boldsymbol{\eta}}_{k}$ is the vector of local degrees of freedom in  $T_{k}$.

\subsection{Post-projection stage}

Once the position of the free surface has been obtained at each primal node, we can use \eqref{eqn.q.update} to correct the intermediate conservative variable related to the velocity and get $\q^{n+1}$ at each dual cell
\begin{equation}
	\label{eqn.q.update.disc} 
	\mathbf{q}^{n+1}_{i} =  \frac{1}{1+\Delta t\, \gamma^{n}_{i}} \left(   \mathbf{q}^*_{i} - \Delta t \, g \, h^n_{i} \, \left( \nabla \eta^{n+1}\right)_{i} \right) .   
\end{equation}
The gradients involved in the former equation are approximated by a weighted interpolation of their values from the primal grid to the dual one as follows:
\begin{equation}
	\left( \nabla \eta^{n+1}\right)_{i} = \frac{1}{\left|C_{i} \right|} \sum_{k\in\mathcal{J}_{i}} \left|T_{ki} \right| \left( \nabla \eta^{n+1}\right)_{T_{k}}, 
\end{equation}
where $\mathcal{J}_{i}$ represents the set of primal elements related to $C_{i}$, $\left|T_{ki}\right|$ denotes the area of the subtriangle in $T_{k}$ used for the construction of $C_{i}$ and $\left( \nabla \eta^{n+1}\right)_{T_{k}}$
is the value of the free surface gradient approximated in $T_{k}$ using a Galerkin approach. 

\subsection{Theta method}
The {implicit part of} the former algorithm is only of first order in time. To increase the accuracy of the overall scheme we introduce a semi-implicit discretization thanks to the $\theta$-method. Let us consider $\theta\in[0.5,1]$ and define
\begin{eqnarray}
	\q^{n+\theta} = (1-\theta)\q^{n} + \theta \q^{n+1}, \label{eqn.q.theta}\\
	\eta^{n+\theta} = (1-\theta)\eta^{n} + \theta \eta^{n+1} \label{eqn.eta.theta}
\end{eqnarray}
the approximations of the solution at $t^{n+\theta}$. Replacing \eqref{eqn.q.theta} in \eqref{eqn.eta.sd}  and \eqref{eqn.eta.theta} in \eqref{eqn.q.sd}, we get modified equations for the computation of the free surface and the momentum at the new time step:
\begin{eqnarray}
	\eta^{n+1} &=& \eta^{n}-\Delta t \left(1-\theta\right)\dive \q^{n} - \Delta t\, \theta \dive \q^{n+1}, \label{eqn.eta.sd_theta}\\
	 \q^{n+1} &=& \frac{1}{ 1+\Delta t\, \gamma^{n}}\left[ \q^{\ast} -\Delta t\, g\, h^{n}\left(1-\theta\right) \grae \eta^{n} -\Delta t\,  g \, h^{n} \theta \grae \eta^{n+1}\right] .
	\label{eqn.q.sd_theta}
\end{eqnarray}
Thus, formal substitution of \eqref{eqn.q.sd_theta} in \eqref{eqn.eta.sd_theta} gives 
\begin{eqnarray}
	\eta^{n+1}  
	-\Delta t^{2}\, \theta^{2}\, g \dive\left( \frac{h^{n}}{ 1+\Delta t\, \gamma^{n}}\grae \eta^{n+1}\right)  
	&=& \eta^{n} - \Delta t \left(1-\theta\right) \dive \q^{n} 
	-\Delta t\, \theta  \dive \left( \frac{1}{ 1+\Delta t\, \gamma^{n}} \, \q^{\ast} \right) \notag\\
	&&
	+ \Delta t^{2}\, \theta \left(1-\theta\right) \dive \left(  \frac{g h^{n}}{ 1+\Delta t\, \gamma^{n}} \grae \eta^{n} \right).
\end{eqnarray}
Defining the new auxiliary variable
\begin{equation}
	\q^{\ast\ast} =  \left(1-\theta\right)  \q^{n} +  \frac{\theta}{ 1+\Delta t\, \gamma^{n}}  \, \q^{\ast}
\end{equation}
leads to 
\begin{eqnarray}
	\eta^{n+1}  
	-\Delta t^{2} \,\theta^{2}\, g \dive\left(  \frac{h^{n}}{ 1+\Delta t\, \gamma^{n}}\grae \eta^{n+1}\right) 
	= \eta^{n} 
	- \Delta t \dive   \q^{\ast\ast}  	+ \Delta t^{2}\, \theta \left(1-\theta\right) g \dive \left( \frac{h^{n}}{ 1+\Delta t\, \gamma^{n}} \grae \eta^{n}\right) .
\end{eqnarray}
Following the steps already introduced in the projection stage, Section \ref{sec.projection}, we can rewrite the above equation in variational formulation obtaining
\begin{weakproblem}\label{wp.theta}
	Find $\eta^{n+1}\in V_{0}$ such that
	\begin{eqnarray}
		\int_{\Omega}\eta^{n+1} z dV  
		+\Delta t^{2}\, \theta^{2}\, g \int_{\Omega} \left( \frac{h^{n}}{ 1+\Delta t\, \gamma^{n}}\grae 	\eta^{n+1}\right)\cdot \grae z dV = 
		\int_{\Omega}\eta^{n} z dV 
		+ \Delta t \int_{\Omega}  \q^{\ast\ast} \grae z dV \notag\\
		- \Delta t \int_{\Omega}  \q^{n+\theta}  \cdot \n\,  z dV 		- \Delta t^{2}\, \theta \left(1-\theta\right) g \int_{\Omega} \frac{h^{n}}{ 1+\Delta t\, \gamma^{n}} \grae 	\eta^{n} \cdot \grae z dV
		\label{eqn.wp2} 
	\end{eqnarray}
	for all $z\in V_{0}$.
\end{weakproblem}
From the numerical point of view, the main difference between weak problems \ref{wp.theta} and \ref{wp.orig} is the inclusion of the last term in the right hand side of the system corresponding to a stiffness matrix for the solution at the previous time step $\eta^{n}$. Moreover, the boundary integral of the conservative velocity is now computed at the intermediate time $t^{n+\theta}=t^{n}+\theta$.
Finally, applying the $\theta$-method leads to a modification of the post projection stage, where the correction of $\q^{n+1}$ must be performed using $\eta^{n+\theta}$ recovered from \eqref{eqn.eta.theta} and not from $\eta^{n+1}$, the solution of the Weak Problem \ref{wp.theta}. Defining $\theta=0.5$ yields the second order Crank-Nicolson method, whereas for $\theta=1$ we revert to the backward Euler scheme as used in \eqref{eqn.weakproblem}. {We emphasize, however, that despite the use of the $\theta$-method the overall scheme always remains formally first order accurate in time, due to the first order operator splitting between the explicit convective subsystem and the implicit pressure subsystem. Nevertheless, the $\theta$-method with $\theta \approx 0.5$ yields better results for \textit{smooth} wave propagation problems than the implicit Euler scheme with $\theta=1$. Since the Crank-Nicolson method is well-known to produce spurious oscillations in the vicinity of discontinuities, it is recommended to use $\theta=1$ for problems involving shock waves. } 

{
\subsection{Pressure correction formulation}
Similar to classical pressure correction algorithms for the incompressible Navier-Stokes equations \cite{PS72,Pat80,BFTVC17}, the present semi-implicit method can also be rewritten in pressure-correction formulation as follows:  
Introducing the pressure-correction 
\begin{equation}
\delta \eta^{n+1} = \eta^{n+1} - \eta^n
\end{equation} 
and the intermediate momentum which takes into account the pressure gradients at time $t^n$ 
\begin{equation}
\label{eqn.q2star.pc} 
\q^{**} = \q^n - \Delta t \, \nabla \cdot \left( \vel^n \otimes \q^n \right) - \Delta t \, g h^n \,  \nabla  \eta^n
        = \q^* - \Delta t \, g h^n \,  \nabla \eta^n,
\end{equation}
one can rewrite the semi-discrete scheme \eqref{eqn.eta.sd}-\eqref{eqn.q.sd} as 
\begin{eqnarray}
\label{eqn.eta.pc} 
\eta^{n+1} &=& \eta^n  - \Delta t \,  \nabla \cdot \q^{n+1}, \\ 
\label{eqn.q.pc} 
\q^{n+1} &=& \q^{**} - \Delta t \, g \, h^n \, \nabla \delta \eta^{n+1} - \Delta t \, \gamma^{n} \q^{n+1}.   
\end{eqnarray}
Inserting the semi-discrete momentum equation into the semi-discrete continuity equation and subtracting $\eta^n$ from both sides of the equation yields the semi-discrete wave equation for the pressure correction  
\begin{equation} 
\label{eqn.wave.eta.pc} 
\delta \eta^{n+1} - \Delta t^2 \,  g \, \nabla \cdot \left(\frac{h^n}{1+\Delta t\, \gamma^{n}} \nabla \delta \eta^{n+1} \right) =  - \Delta t \,  \nabla \cdot \left( \frac{1}{1+\Delta t\, \gamma^{n}}\, \q^{**}\right) .  
\end{equation} 	
The weak problem associated with \eqref{eqn.wave.eta.pc} is 
\begin{weakproblem}\label{wp.pc}
	Find the pressure correction  $\delta \eta^{n+1}\in V_{0}$ that satisfies
	\begin{equation}\label{eqn.weakproblem.pc}
	\int_{\Omega} \delta \eta^{n+1} z dV 
	+ \Delta t^2 g \int_{\Omega}  \frac{h^{n}}{1+\Delta t\, \gamma^{n}} \grae \! \delta \eta^{n+1}\cdot \grae \! z dV =  
	 \Delta t \int_{\Omega} \frac{1}{1+\Delta t\, \gamma^{n}}\, \q^{\ast \ast}\cdot \grae z dV 
	- \Delta t \int_{\Gamma}\q^{n+1}\cdot \n\, z dV
	\end{equation}
	for all $z \in V_{0}$.
\end{weakproblem}
}

{
\subsection{Asymptotic behaviour of the semi-implicit hybrid FV/FE method in the low Froude number limit}\label{sec:AP}
It is a well known result from asymptotic analysis of the compressible Euler equations, see \cite{KlaMaj,KlaMaj82,Klein95,MeisterMach,Klein2001,MRKG03}, that for low Mach numbers and in the absence of compression from the boundary the pressure tends to a global constant, plus some small fluctuations around this constant of order $M^2$. For free surface shallow water flows, the role of the Mach number is taken over by the Froude number
 Fr$ = \| \mathbf{v} \| / \sqrt{g h}$. Assuming zero friction ($\gamma=0$), flat bottom ($b=0$ and thus $\eta = h$), $h = h_0 + Fr^2 h_2$, and using $\mathbf{q} = h \mathbf{v}$ and $p=g h$ we can therefore rewrite the weak problem \eqref{eqn.weakproblem.pc} as  
 	\begin{equation*}\label{eqn.weakproblem.pc.asy1}
\int_{\Omega} \delta \eta^{n+1} z dV 
+ \Delta t^2 g \int_{\Omega}  (h_0^{n} + Fr^2 h_2^n) \grae \! \delta \eta^{n+1}\cdot \grae \! z dV =  
 \Delta t \int_{\Omega} \, \q^{\ast \ast}\cdot \grae z dV 
- \Delta t \int_{\Gamma}\q^{n+1}\cdot \n\, z dV. 
\end{equation*}
 Neglecting the quadratic terms in the Froude number and dividing by the global constant $h_0$ yields 
 	\begin{equation*}\label{eqn.weakproblem.pc.asy2}
\frac{1}{h_0} \int_{\Omega} \delta \eta^{n+1} z dV 
+ \Delta t^2  \int_{\Omega}  \grae \! \delta p^{n+1}\cdot \grae \! z dV =  
 \Delta t \int_{\Omega} \, \mathbf{v}^{\ast \ast}\cdot \grae z dV 
- \Delta t \int_{\Gamma}\mathbf{v}^{n+1}\cdot \n\, z dV.
\end{equation*}
In the limit Fr$ \to 0$, which corresponds to $c \to \infty$ and therefore to $h_0 \to \infty$ we therefore obtain 
\begin{equation}\label{eqn.weakproblem.pc.asy}
   \int_{\Omega}  \grae \! \delta p^{n+1}\cdot \grae \! z dV =  
\frac{1}{\Delta t}  \int_{\Omega} \, \mathbf{v}^{\ast \ast}\cdot \grae z dV 
- \frac{1}{\Delta t} \int_{\Gamma}\mathbf{v}^{n+1}\cdot \n\, z dV,
\end{equation}
which is the classical Poisson equation for the pressure correction of the projection stage of the hybrid FV/FE scheme for the incompressible Navier-Stokes equations presented in \cite{BFTVC17}. We can therefore conclude that our scheme is compatible with a hybrid FV/FE discretization of the incompressible Navier-Stokes equations when Fr $ \to 0$. As such, our algorithm is asymptotic preserving (AP), which was expected, since we are using a pressure-based semi-implicit discretization on staggered meshes that is typical for the simulation of incompressible flows. In the next section we will confirm the AP property for Fr $\to 0$ also numerically. 
Via completely analogous computations one can also show that the discrete pressure equation of the scheme without pressure correction, i.e. \eqref{eqn.wp2}, tends to the discrete pressure Poisson equation 
\begin{equation}\label{eqn.weakproblem2.pc.asy}
\int_{\Omega}  \grae \! p^{n+1}\cdot \grae \! z dV =  
\frac{1}{\Delta t}  \int_{\Omega} \, \mathbf{v}^{\ast}\cdot \grae z dV 
- \frac{1}{\Delta t} \int_{\Gamma}\mathbf{v}^{n+1}\cdot \n\, z dV
\end{equation} 
in the low Froude number limit Fr $ \to 0$. 
}

\section{Numerical results}\label{sec:numericalresults}
Several numerical test cases are addressed aiming at validating the numerical method presented in the previous section. 
In all test cases the following CFL-type condition is employed to determine the time step size:  
\begin{equation*}
	\Delta t = \min_{C_{i}}\left\lbrace \Delta t_{i}\right\rbrace, \qquad \Delta t_{i} =  \frac{\textrm{CFL}\; \mathcal{L}_{i}}{ | \lambda_{i} | }
\end{equation*}
with $|\lambda_{i}|$ the maximum absolute eigenvalue in cell $C_{i}$ and $\mathcal{L}_{i}$ its incircle diameter. The friction coefficient is assumed to be zero unless a specific value is indicated. Similarly, the default option for the theta method is $\theta = 1$.

\subsection{Numerical convergence analysis}
As first test case, we consider a steady shallow water vortex problem, see \cite{DC13}, for which the exact solution is known:
\begin{equation*}
\vels_{1}   = -v_\varphi \sin(\phi), \qquad
\vels_{2}   =  v_\varphi\cos(\phi),  \qquad
h			= {h_0} - \frac{1}{2 g} e^{-(r^2-1)}, \qquad
b			= 0.
\end{equation*}
Here, $r = \sqrt{x^2+y^2}$, $\tan\left( \phi \right) = \frac{y}{x}$ and $v_\varphi = r\, e^{-\frac{1}{2}\left(r^{2}-1\right)}$.
The {first set of simulations of this test} is run using the hybrid FV/FE scheme on the computational domain $\Omega=[-5,5]\times[-5,5]$ discretized with a set of meshes of  $N_{x}=N_{y}\in \left\lbrace32,64,128,256,512\right\rbrace$ divisions along each boundary, setting {$h_0=1$, hence the corresponding Froude number is Fr$=0.32$}. For the discretization of the nonlinear convective terms the second order LADER methodology is employed. The final simulation time is set to $T = 0.1$ and {the time step is fixed for each mesh with $\Delta t\in\left\lbrace 10^{-2}, 5.0 \cdot 10^{-3}, 2.5 \cdot 10^{-3}, 1.25 \cdot 10^{-3}, 6.25 \cdot 10^{-4} \right\rbrace$}. Periodic boundary conditions are considered everywhere. Table \ref{TGV_errors} reports the errors obtained for the main flow variables: the {velocity components, $\vels_{1}$, $\vels_{2}$}, and the free surface elevation, $\eta$. We observe that the sought second order of accuracy is reached.

\begin{table}[ht]
	\renewcommand{\arraystretch}{1.2}
	\begin{center}
		\caption{Shallow water vortex test problem. Spatial $L^{2}$ error norms {in space} obtained at time $T=0.1$, and convergence rates for the hybrid FV/FE method using the explicit local ADER scheme for the nonlinear convective terms.} \label{TGV_errors}
		
		\vspace{4pt}		
		{\begin{tabular}{ccccccc}
			\hline 
			Mesh &$L^{2}_{\Omega}\left(\vels_{1}\right)$ & $\mathcal{O}\left(\vels_{1} \right)$                  
			&$L^{2}_{\Omega}\left(\vels_{2} \right)$ & $\mathcal{O}\left(\vels_{2} \right)$ &  $L^{2}_{\Omega}\left(\press\right)$ & $\mathcal{O}\left(\press\right)$ \\ 
			\hline 
			32  & $1.3311 \cdot 10^{-2}$ & $   $ & $1.3499 \cdot 10^{-2} $ & $   $ & $1.4099 \cdot 10^{-3} $ & $    $ \\
			64  & $3.4160 \cdot 10^{-3}$ & $2.0$ & $3.3813 \cdot 10^{-3} $ & $2.0$ & $3.9320 \cdot 10^{-4} $ & $1.8$ \\
			128 & $8.6462 \cdot 10^{-4}$ & $2.0$ & $8.4725 \cdot 10^{-4} $ & $2.0$ & $1.0217 \cdot 10^{-4} $ & $1.9$ \\
			256 & $2.1740 \cdot 10^{-4}$ & $2.0$ & $2.1211 \cdot 10^{-4} $ & $2.0$ & $2.5971 \cdot 10^{-5} $ & $2.0$ \\
			512 & $5.4650 \cdot 10^{-5}$ & $2.0$ & $5.3219 \cdot 10^{-5} $ & $2.0$ & $6.5433 \cdot 10^{-6} $ & $2.0$ \\
			\hline 
		\end{tabular}}
	\end{center}
\end{table}

{The second set of simulations aims at showing, numerically, the asymptotic preserving property (AP) of the scheme in the low Froude number limit (Fr $\to 0$), already discussed in Section \ref{sec:AP} from the theoretical point of view. To this end, we run the pressure-correction version of the proposed methodology for Froude numbers from $10^{-1}$ down to $10^{-7}$. Hence, the corresponding values for $h_{0}$ fall between $10^{1}$ and $10^{13}$. Even if the three coarsest grids are used, for  $h_{0}\geq 10^{5}$ the magnitude of the water depth compared with the errors obtained for this variable is greater than $10^{9}$, i.e. more than nine orders of magnitude, making necessary the use of quadruple precision arithmetics for the computations. As it can be observed in Table \ref{SWV_Fr}, the numerical obtained results clearly demonstrate the AP property of the scheme, since the second order of accuracy is maintained, independently of the Froude number. Furthermore, even the absolute errors of all variables become independent of the Froude number in the limit Fr $\to 0$. 
}

\begin{table}[ht!]
	\renewcommand{\arraystretch}{1.2}
	\begin{center}
		{\caption{Shallow water vortex test problem. Numerical analysis of the asymptotic preserving property (AP) of the pressure-correction algorithm with \mbox{ Fr $\in\left\lbrace 10^{-1}, 10^{-2}, 10^{-3}, 10^{-4}, 10^{-5}, 10^{-6},  10^{-7} \right\rbrace$ at time $T=0.1$.}}\label{SWV_Fr}} 

		\vspace{4pt}		
		{\begin{tabular}{ccccccc}
				\hline 
				Mesh &$L^{2}_{\Omega}\left(\vels_{1}\right)$ & $\mathcal{O}\left(\vels_{1} \right)$                  
				&$L^{2}_{\Omega}\left(\vels_{2} \right)$ & $\mathcal{O}\left(\vels_{2} \right)$ &  $L^{2}_{\Omega}\left(\press\right)$ & $\mathcal{O}\left(\press\right)$ \\ 
				\hline 
				\multicolumn{7}{c}{ Fr $= 10^{-1}$, $h_0 = 10^1$ \hspace{10pt}{\scriptsize (double precision)}}
				\\ \hline 
				32  & $1.0777 \cdot 10^{-2}$ & $   $ & $1.0352 \cdot 10^{-2} $ & $   $ & $8.3811 \cdot 10^{-3}$ & $   $ \\
				64  & $2.7161 \cdot 10^{-3}$ & $2.0$ & $2.6160 \cdot 10^{-3} $ & $2.0$ & $2.3257 \cdot 10^{-3}$ & $1.8$ \\
				128 & $6.8350 \cdot 10^{-4}$ & $2.0$ & $6.6641 \cdot 10^{-4} $ & $2.0$ & $6.0772 \cdot 10^{-4}$ & $1.9$ \\
				\hline 				
				\multicolumn{7}{c}{Fr $= 10^{-2}$, $h_0 = 10^3$ \hspace{10pt}{\scriptsize (double precision)}}
				\\ \hline 
				32  & $1.0544 \cdot 10^{-2}$ & $   $ & $9.7097 \cdot 10^{-3} $ & $   $ & $6.2106\cdot 10^{-3}$ & $   $ \\
				64  & $2.6177 \cdot 10^{-3}$ & $2.0$ & $2.4821 \cdot 10^{-3} $ & $2.0$ & $1.6155\cdot 10^{-3}$ & $1.9$ \\
				128 & $6.5243 \cdot 10^{-4}$ & $2.0$ & $6.3628 \cdot 10^{-4} $ & $2.0$ & $4.0767\cdot 10^{-4}$ & $2.0$ \\
				\hline 
				\multicolumn{7}{c}{Fr $= 10^{-3}$, $h_0 = 10^5$ \hspace{10pt}{\scriptsize (quadruple precision)}}
				\\ \hline 
				32  & $1.0544 \cdot 10^{-2}$ & $   $ & $9.7098 \cdot 10^{-3} $ & $   $ & $6.2218 \cdot 10^{-3}$ & $   $ \\
				64  & $2.6177 \cdot 10^{-3}$ & $2.0$ & $2.4823 \cdot 10^{-3} $ & $2.0$ & $1.6226 \cdot 10^{-3}$ & $1.9$ \\
				128 & $6.5245 \cdot 10^{-4}$ & $2.0$ & $6.3633 \cdot 10^{-4} $ & $2.0$ & $4.1224 \cdot 10^{-4}$ & $2.0$ \\
				\hline
				\multicolumn{7}{c}{Fr $= 10^{-4}$, $h_0 = 10^7$ \hspace{10pt}{\scriptsize (quadruple precision)}}
				\\ \hline 
				32  & $1.0544 \cdot 10^{-2}$ & $   $ & $9.7098 \cdot 10^{-3} $ & $   $ & $6.2218 \cdot 10^{-3}$ & $   $ \\
				64  & $2.6177 \cdot 10^{-3}$ & $2.0$ & $2.4823 \cdot 10^{-3} $ & $2.0$ & $1.6226 \cdot 10^{-3}$ & $1.9$ \\
				128 & $6.5245 \cdot 10^{-4}$ & $2.0$ & $6.3633 \cdot 10^{-4} $ & $2.0$ & $4.1224 \cdot 10^{-4}$ & $2.0$ \\
				\hline 
				\multicolumn{7}{c}{Fr $= 10^{-5}$, $h_0 = 10^9$ \hspace{10pt}{\scriptsize (quadruple precision)}}
				\\ \hline 
				32  & $1.0544 \cdot 10^{-2}$ & $   $ & $9.7098 \cdot 10^{-3} $ & $   $ & $6.2218 \cdot 10^{-3}$ & $   $ \\
				64  & $2.6177 \cdot 10^{-3}$ & $2.0$ & $2.4823 \cdot 10^{-3} $ & $2.0$ & $1.6226 \cdot 10^{-3}$ & $1.9$ \\
				128 & $6.5245 \cdot 10^{-4}$ & $2.0$ & $6.3633 \cdot 10^{-4} $ & $2.0$ & $4.1224 \cdot 10^{-4}$ & $2.0$ \\
				\hline 
				\multicolumn{7}{c}{Fr $= 10^{-6}$, $h_0 = 10^{11}$ \hspace{10pt}{\scriptsize (quadruple precision)}}
				\\ \hline 
				32  & $1.0544 \cdot 10^{-2}$ & $   $ & $9.7098 \cdot 10^{-3} $ & $   $ & $6.2218 \cdot 10^{-3}$ & $   $ \\
				64  & $2.6177 \cdot 10^{-3}$ & $2.0$ & $2.4823 \cdot 10^{-3} $ & $2.0$ & $1.6226 \cdot 10^{-3}$ & $1.9$ \\
				128 & $6.5245 \cdot 10^{-4}$ & $2.0$ & $6.3633 \cdot 10^{-4} $ & $2.0$ & $4.1224 \cdot 10^{-4}$ & $2.0$ \\
				\hline 
				\multicolumn{7}{c}{Fr $= 10^{-7}$, $h_0 = 10^{13}$ \hspace{10pt}{\scriptsize (quadruple precision)}}
				\\ \hline 
				32  & $1.0544 \cdot 10^{-2}$ & $   $ & $9.7098 \cdot 10^{-3} $ & $   $ & $6.2218 \cdot 10^{-3}$ & $   $ \\
				64  & $2.6177 \cdot 10^{-3}$ & $2.0$ & $2.4823 \cdot 10^{-3} $ & $2.0$ & $1.6226 \cdot 10^{-3}$ & $1.9$ \\
				128 & $6.5245 \cdot 10^{-4}$ & $2.0$ & $6.3633 \cdot 10^{-4} $ & $2.0$ & $4.1224 \cdot 10^{-4}$ & $2.0$ \\
				\hline 
		\end{tabular}}
	\end{center}
\end{table}

\subsection{Numerical verification of the well-balanced property (C-property)}\label{sec:wellbalance}
The C-property in the sense of \cite{BVC94,LeVeque98},  
also known as well-balance property, characterises the capability of numerical methods to exactly preserve stationary equilibrium solutions of the shallow water equations of the form $\eta(\x,t) = const$ and  $\vel(\x,t) = \mathbf{0}$ for arbitrary bottom topography $b(\x)$. One can easily observe that the numerical method proposed in this paper is well balanced by construction. To check it also numerically, we consider the classical benchmark proposed by LeVeque in \cite{LeVeque98} and characterised by the following variable bottom bathymetry
\begin{equation*}
	b(\x,t) = 0.8 e^{-5 \left(x+0.1\right)^2 - 50 y^2},
\end{equation*}
an initial free surface of the form
\begin{equation*}		
	\eta\left(\mathbf{x},0\right) = \left\lbrace \begin{array}{lc}
		1 + \epsilon & \mathrm{ if } \; -0.95 \leq x \leq  -0.85,\\
		1  & \mathrm{ elsewhere, } 
	\end{array}\right.
\end{equation*}
and fluid at rest conditions for the velocity field, $\vel(\x,0) = \mathbf{0}$. The computational domain, $\Omega=[-2,1]\times[-0.5,0.5]$, is discretised with $54787$ triangular primal elements and we have set $\theta=0.6$. As boundary conditions we define the exact solution in $x$-direction and periodic boundary conditions along the $y$-direction.

Initially, we set $\epsilon=0$ and run the test using different machine precisions. The results reported in Table~\label{key}\ref{tab.WB} confirm the ability of the proposed method to preserve the water at rest solution up to machine precision. 
\begin{table}[ht!]
	\renewcommand{\arraystretch}{1.1}
	\caption{Well balance test ($\epsilon=0$). $L^{2}$ errors obtained for different machine precisions at time $T=0.1$.}
	\label{tab.WB}
	\begin{center}
		\begin{tabular}{ccc}
			\hline 
			Precision & $L^{2}(\eta)$ & $L^{2}( h \vel)$ \\ \hline
			Single    & $1.88211 \cdot 10^{-6} $ & $5.94205 \cdot 10^{-6} $ \\
			Double    & $1.49243 \cdot 10^{-13}$ & $6.78885 \cdot 10^{-13}$ \\
			Quadruple & $8.22214 \cdot 10^{-26}$ & $3.70937 \cdot 10^{-25}$ \\ 
			\hline 
		\end{tabular} 
	\end{center}
\end{table}

\begin{figure}
	\begin{center}
		\includegraphics[trim=10 10 20 10,clip,width=0.49\linewidth]{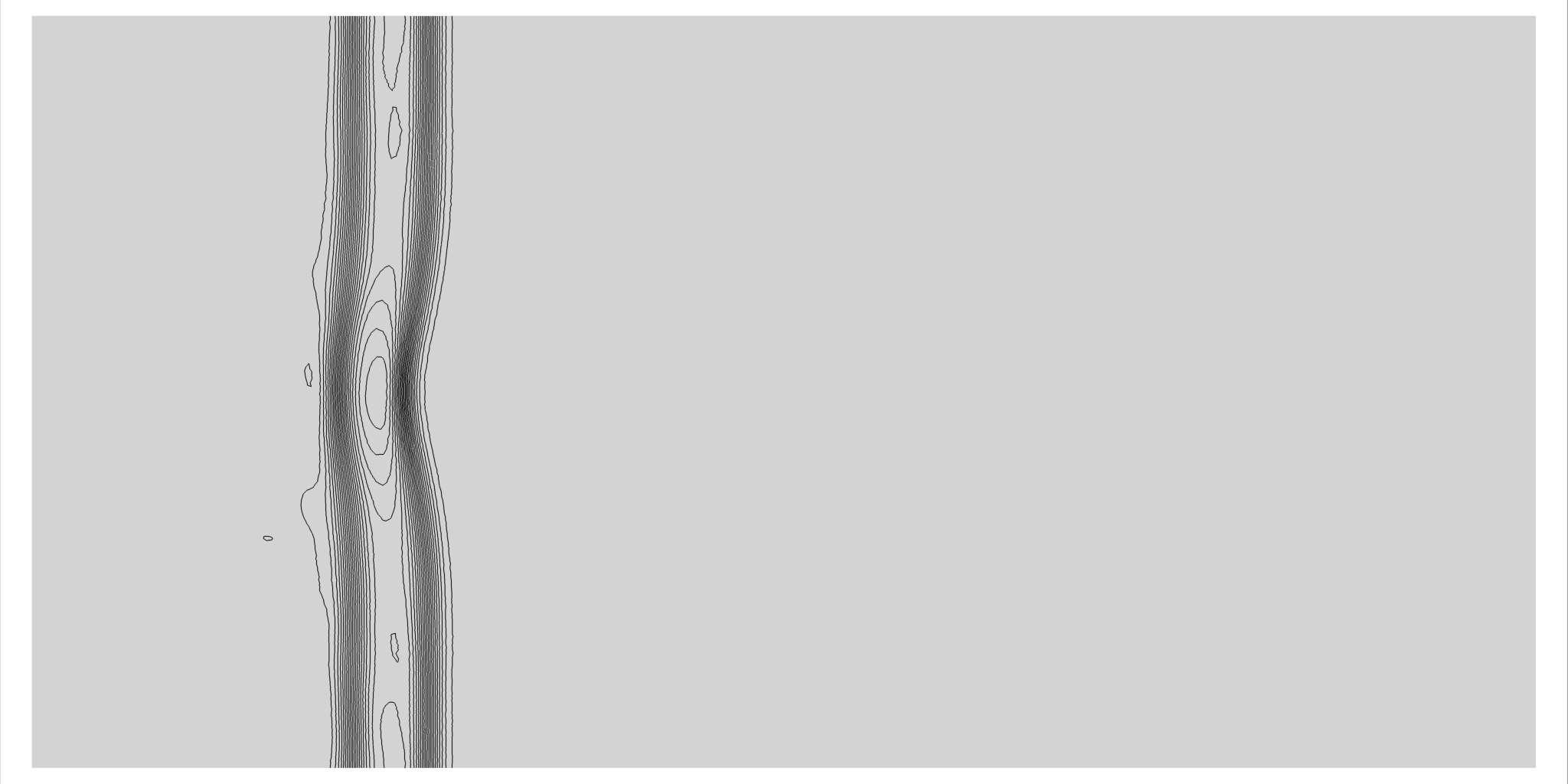}
		\includegraphics[trim=10 10 20 10,clip,width=0.49\linewidth]{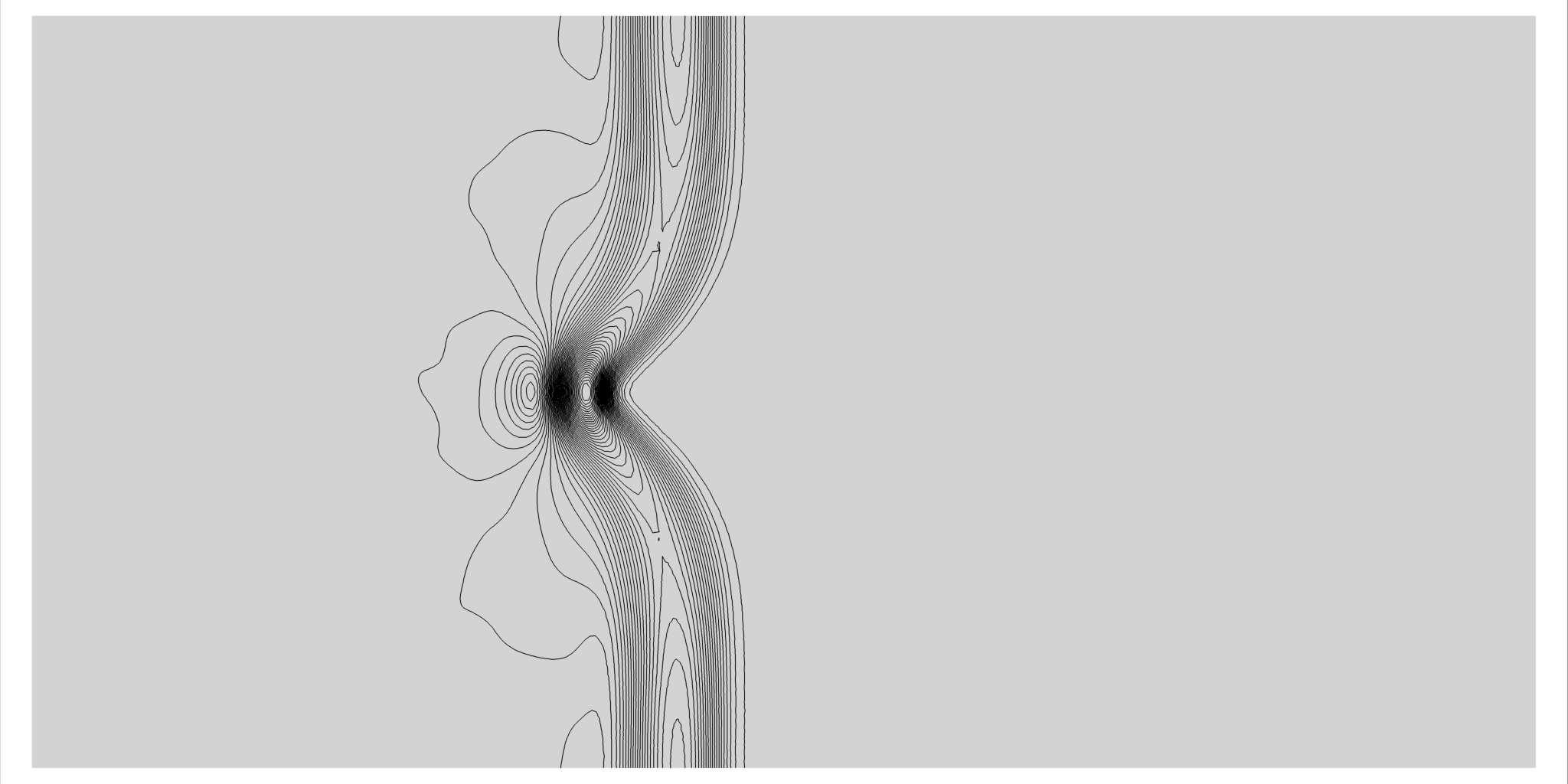} \\[5pt]
		\includegraphics[trim=10 10 20 10,clip,width=0.49\linewidth]{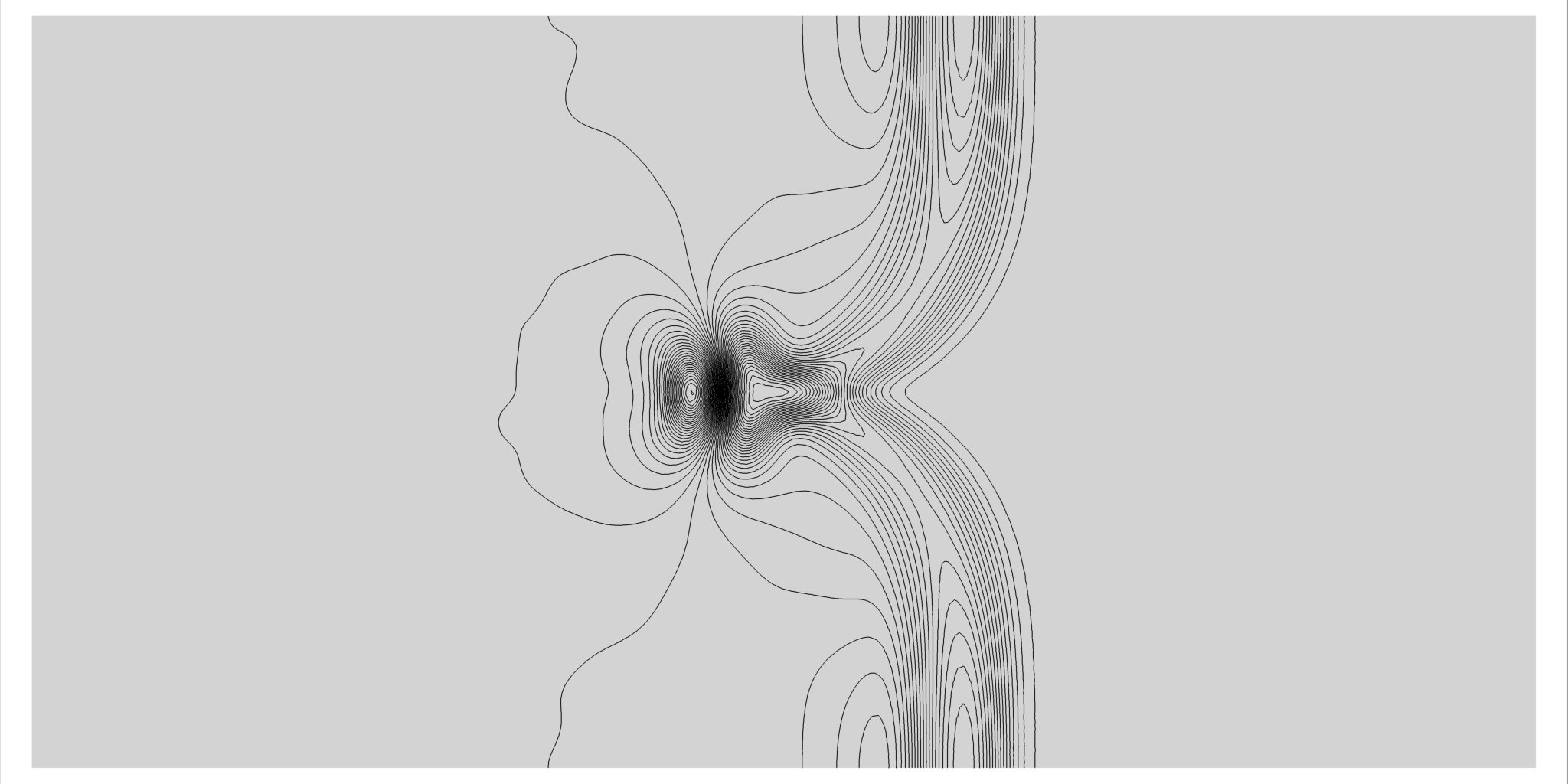}
		\includegraphics[trim=10 10 20 10,clip,width=0.49\linewidth]{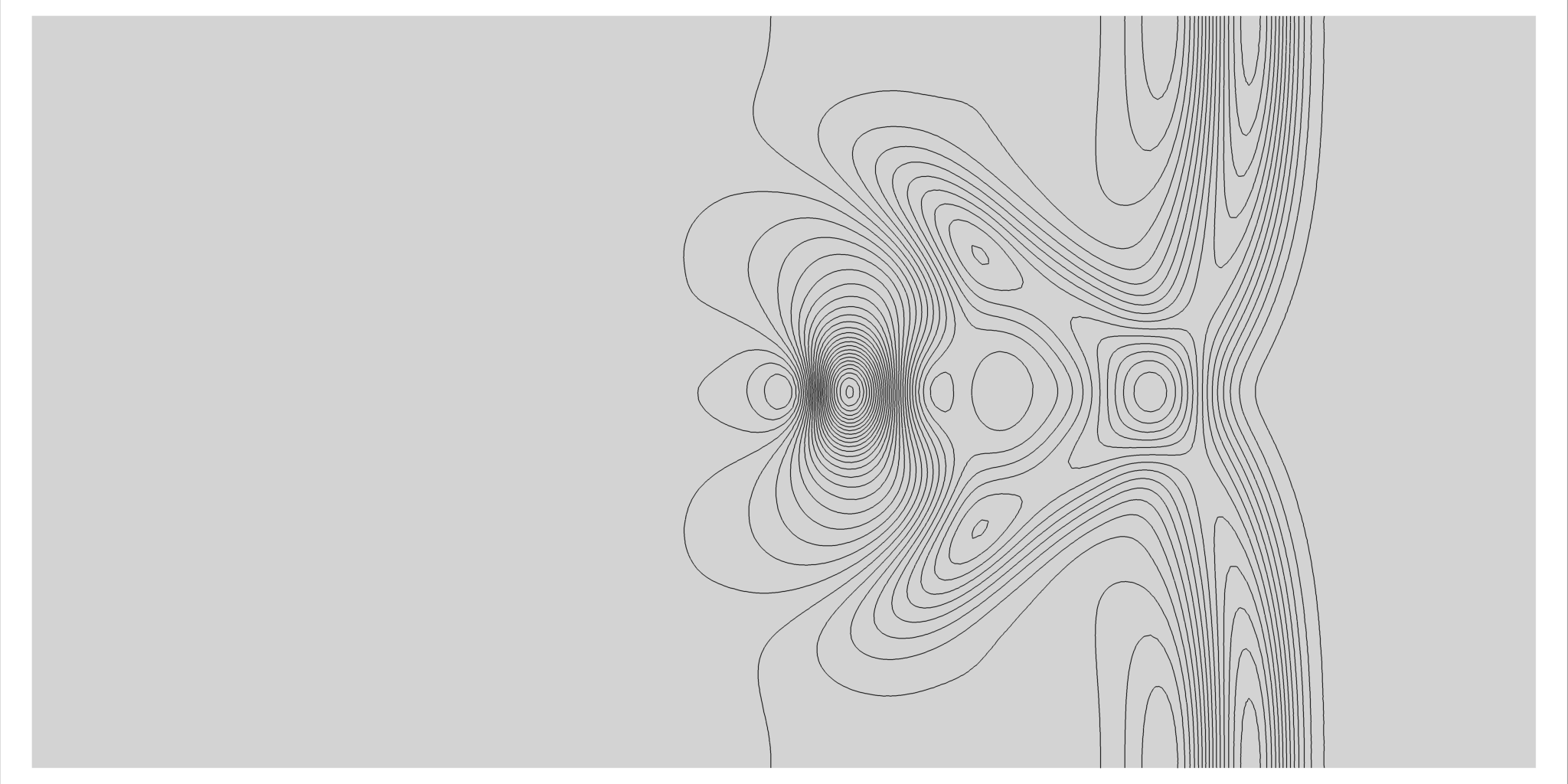} 
	\end{center}
	\caption{2D surface perturbation. Free surface obtained for the well balanced test with an initial surface perturbation of $\epsilon = 10^{-2}$ at times $t\in \left\lbrace 0.12, 0.24, 0.36, 0.48\right\rbrace$ (from top left to bottom right) using LADER scheme with $\theta = 0.51$.}
	\label{fig.WB2D_perturbation}
\end{figure}

As second test, we introduce a small perturbation on the water height by defining $\epsilon = 10^{-2}$.
The solution, depicted in Figure \ref{fig.WB2D_perturbation} for times $t\in\left\lbrace 0.12, 0.24, 0.36, 0.48 \right\rbrace$, shows the expected propagation of the perturbation without the generation of spurious oscillations due to the presence of the Gaussian bottom bump. Moreover, the pattern obtained presents a good agreement with the solutions available in the bibliography, see e.g. \cite{CSDT09,TD14sw}.

\begin{figure}
	\begin{center}
		\includegraphics[width=0.49\linewidth]{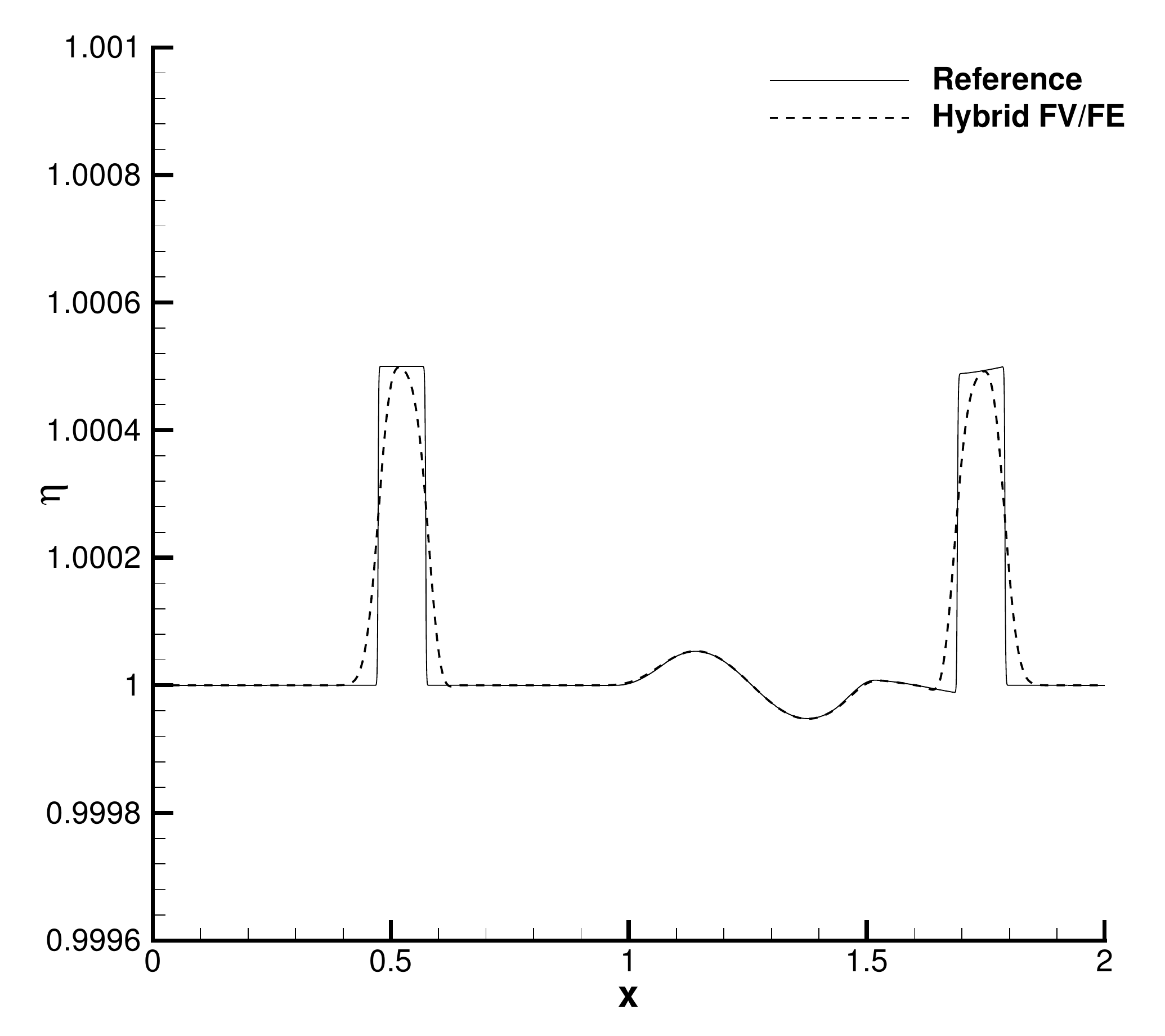}
		\includegraphics[width=0.49\linewidth]{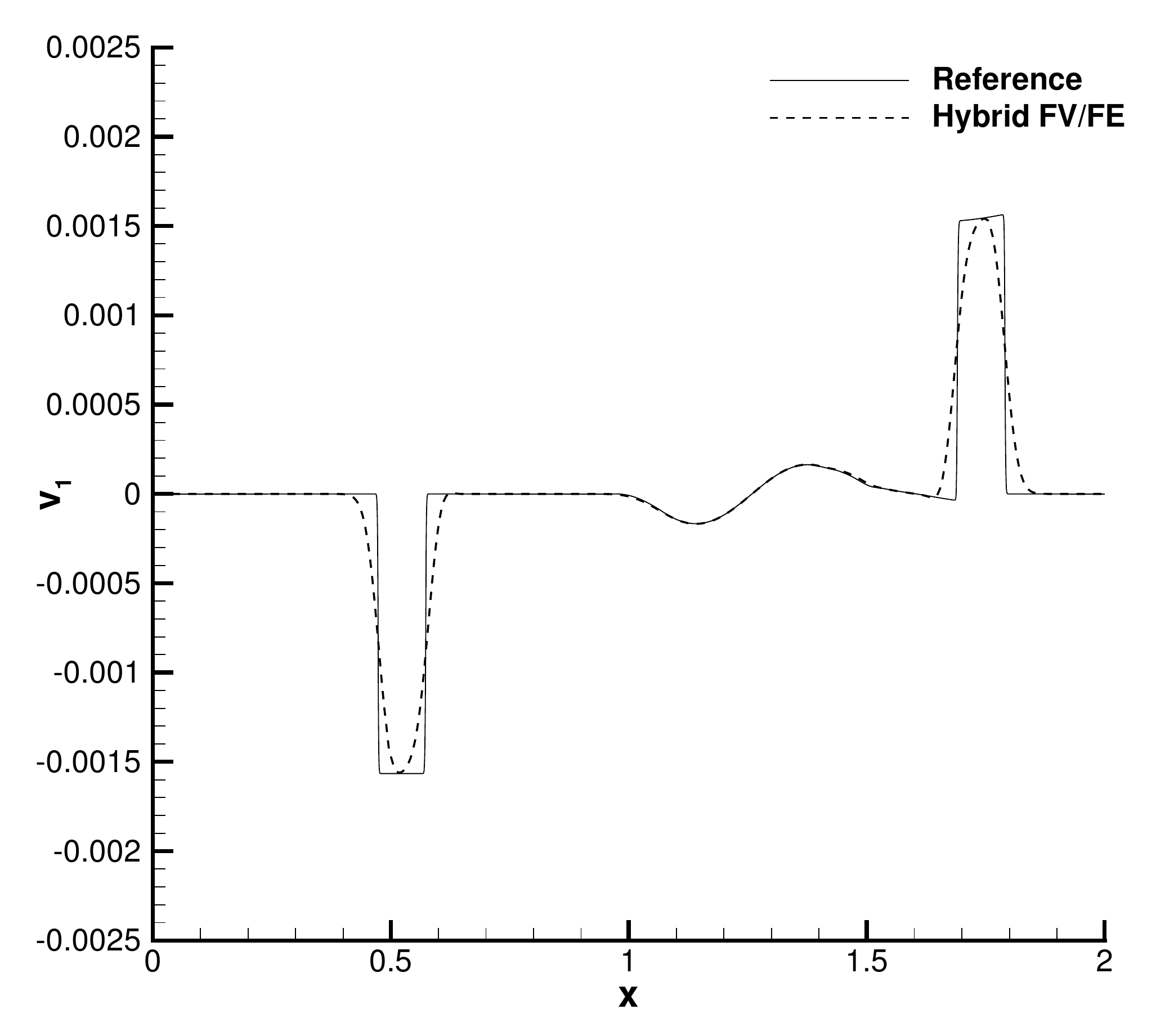}
	\end{center}
	\caption{Small perturbation of a free surface. 1D cut along the $x$ axis through the free surface profile $\eta$ (left) and the velocity in $x-$direction $\vel_{1}$ (right), obtained with the  second order hybrid FV-FE scheme with $\theta=0.51$ (discontinuous black line) and the numerical reference solution obtained with a well-balanced second order TVD-FV scheme (continuous line).} 
	\label{fig.WB1D}
\end{figure}

The third test considered assumes an initial fluid at rest, $\vel(\x,0) = \mathbf{0}$, over a fixed variable bottom 
\begin{equation*}
	b(\x,t) =  \left\lbrace \begin{array}{lc}
		\frac{1}{4}\left( \cos( 10 \pi (x-\frac{3}{2}) ) + 1 \right) & \mathrm{ if } \; -1.6 \leq x \leq  -1.4,\\
		0  & \mathrm{ elsewhere, } 
	\end{array}\right.
\end{equation*}
and an initial free surface given by
\begin{equation*}		
	\eta\left(\mathbf{x},0\right) = \left\lbrace \begin{array}{lc}
		1 + \epsilon & \mathrm{ if } \; 1.1 \leq x \leq  1.2,\\ 
		1  & \mathrm{ elsewhere, } 
	\end{array}\right. \qquad \epsilon = 10^{-3}.
\end{equation*}
The computational domain, $\Omega=[0,2]\times[0,0.2]$, is discreticed with $N_{x}=400$ divisions along the $x$-direction. Since the test has only one-dimensional effects the solution can be compared against a 1D reference solution computed on a very fine mesh of $10^{4}$ cells using the well-balanced second order TVD finite volume scheme proposed in \cite{DT11}. The solution obtained using the LADER hybrid FV-FE scheme  with $\theta=0.51$ is shown in Figure \ref{fig.WB1D}. We can note an excellent agreement between the numerical solution obtained with the new semi-implicit hybrid FV/FE scheme and the fine grid reference solution. 

\subsection{Riemann problems}
A set of Riemann problems both with flat and variable bottom is analysed in this section, see \cite{toro-book-swe,BTT08,DB16} for a detailed description of the resulting wave structures. The initial condition is given by
\begin{equation*}		
	\eta \left(\mathbf{x},0\right) = \left\lbrace \begin{array}{lc}
		\eta_{L} & \mathrm{ if } \; x \le  x_{c},\\
		\eta_{R} & \mathrm{ if } \; x >  x_{c};
	\end{array}\right.
	\quad
	{\vels}_{1} \left(\mathbf{x},0\right) =\left\lbrace \begin{array}{lc}
		\vels_{L}  & \mathrm{ if } \; x \le  x_{c},\\
		\vels_{R} & \mathrm{ if } \; x >  x_{c};
	\end{array}\right. \quad
	{{\vels}_{2}} \left(\mathbf{x},0\right) =0;\quad
 	b\left(\mathbf{x},0\right) = \left\lbrace \begin{array}{lc}
 		b_{L} & \mathrm{ if } \; x \le  x_{c},\\
 		b_{R} & \mathrm{ if } \; x>  x_{c}
 	\end{array}\right.
\end{equation*}
with the left and right states given in Table \ref{tab.RPdef}. The simulations are run in a 2D computational domain with $x\in\left[x_{L},x_{R}\right]$. The values of $x_{L}$ and $x_{R}$ are also defined in Table \ref{tab.RPdef} for each test case. The width of the computational domain in $y$ direction is variable, guaranteeing a small number of layers in $y$-direction while keeping a good aspect ratio according to the number of cells in $x$-direction, $N_{x}$, and leading to a mesh spacing of $\Delta x = 1/N_{x}$. Moreover, the table contains the final time considered for each simulation as well as the initial position of the discontinuity. For all tests we impose the exact solution on the left and right boundaries, while periodic boundary conditions are employed in $y$-direction. 

\begin{figure}[h]
	\centering
	\includegraphics[width=0.49\linewidth]{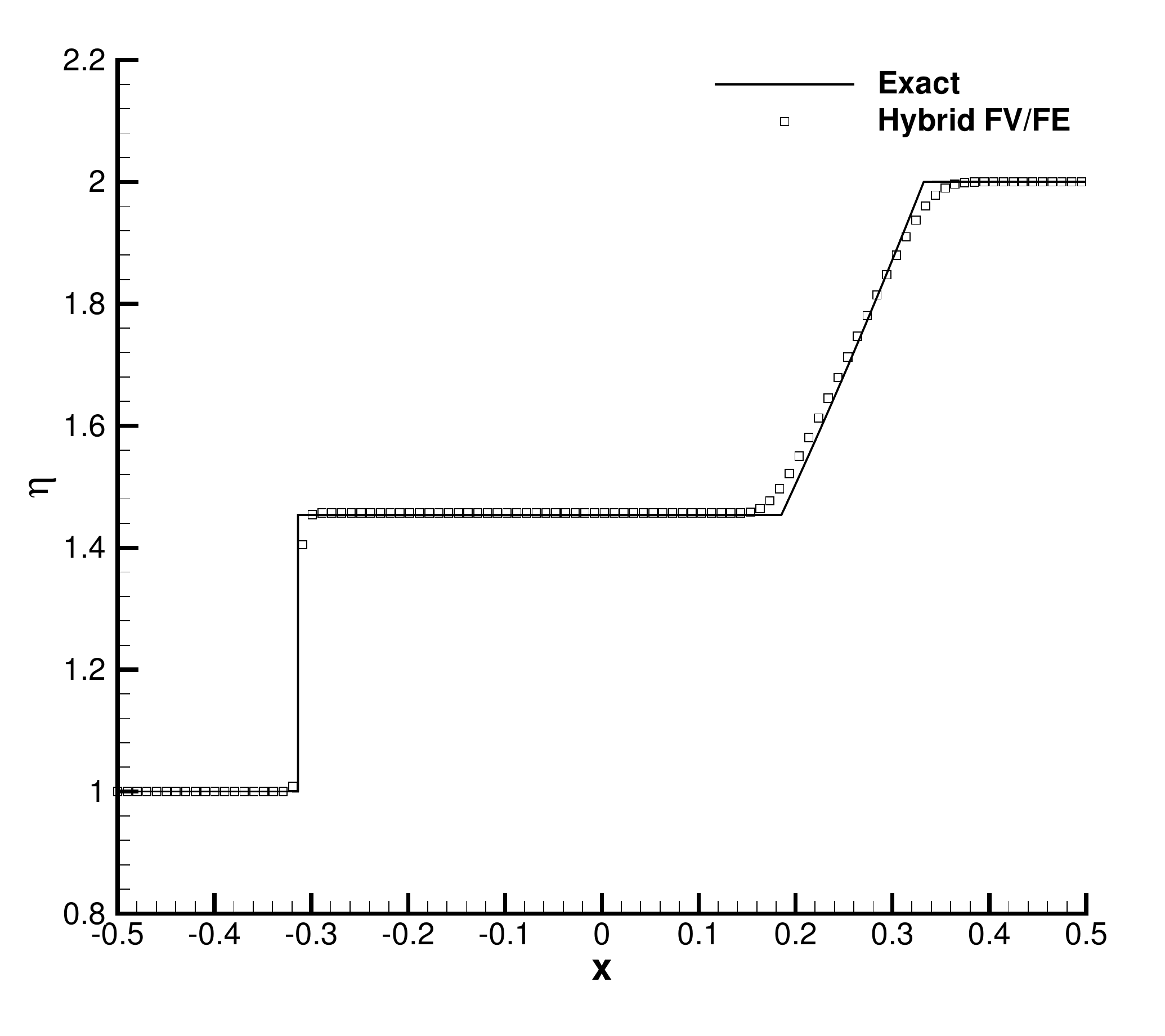}
	\includegraphics[width=0.49\linewidth]{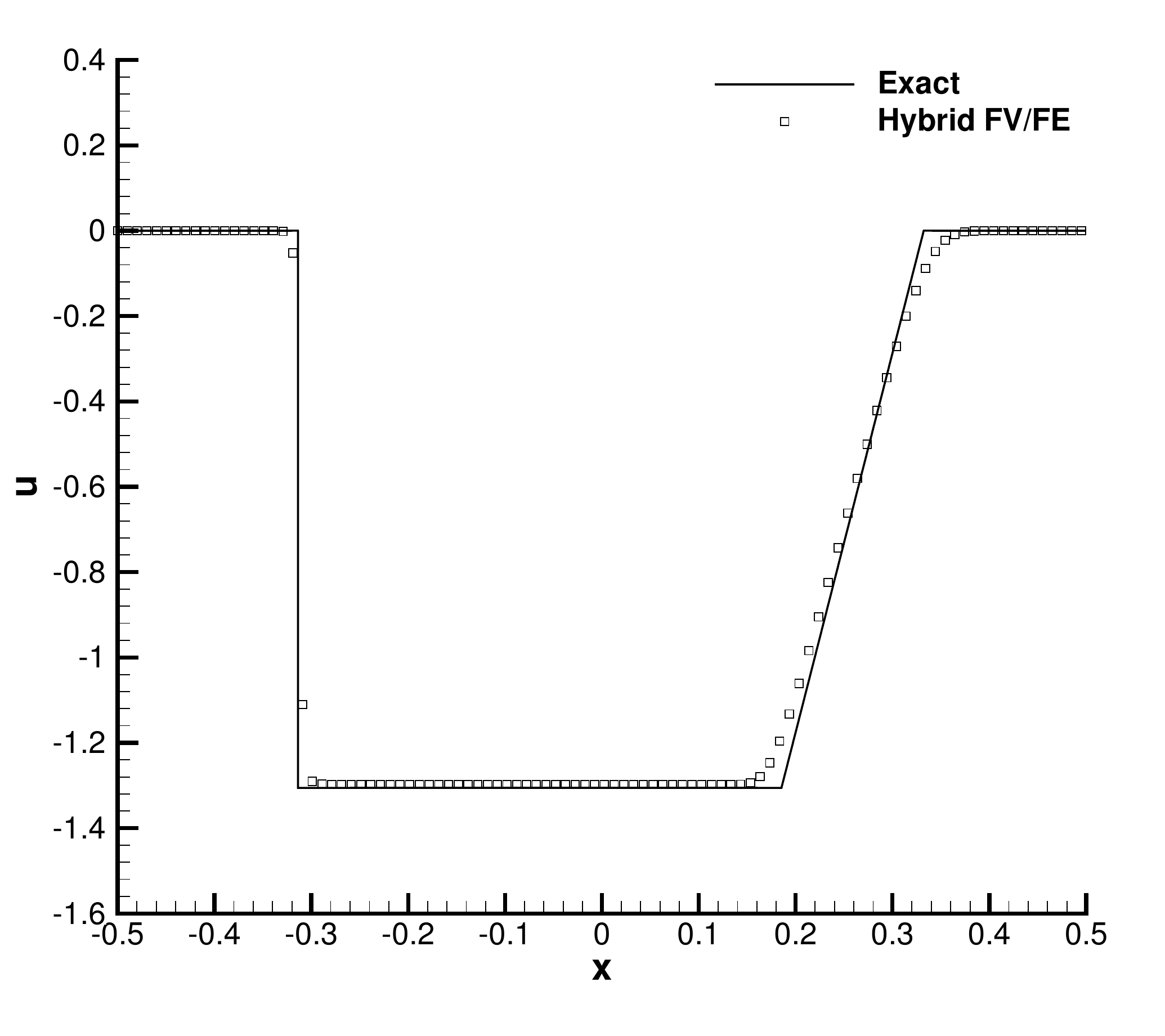}
	\caption{Free surface elevation (left) and horizontal velocity (right) of the first Riemann problem computed using LADER-ENO scheme (squares) and exact solution (continuous line).}
	\label{fig.RP1}
\end{figure}

\begin{table}[ht!]
	\renewcommand{\arraystretch}{1.1}
	\caption{Riemann problems. Left and right states of the initial condition, boundaries of the computational domain, $x_{L}$, $x_{R}$, position of the initial discontinuity, $x_{c}$, number of mesh cells along the $x$-axis, $N_{x}$, and final simulation time, $T$.}
	\label{tab.RPdef}
	\begin{center}
		\begin{tabular}{cccccccccccc}
			\hline 
			Test &  $\eta_{L}$ &  $\eta_{R}$  &  $u_{L}$ &  $u_{R}$ &  $b_{L}$ &  $b_{R}$ & $x_{L}$ & $x_{R}$ & $x_{c}$  & $N_{x}$ & $T$\\ \hline
			RP1 & $ 1 $ & $ 2 $ & $ 0 $ & $ 0 $  & $ 0 $ & $ 0 $& $ -0.5 $ &  $ 0.5 $  & $0$ & $200$ & $0.075$\\
			RP2 & $ 1.46184 $ & $ 0.30873 $ & $ 0 $ & $ 0 $  & $ 0 $ & $ 0.2 $& $ -0.5 $ &  $ 0.5 $  & $0$ & $400$ & $1$\\
			RP3 & $ 0.75    $ & $ 1.10594 $ & $ -9.49365 $   & $ -4.94074    $& $ 0 $ & $ 0.2 $& $ -15 $ &  $ 5 $  & $0$ & $800$ & $1$\\
			RP4 & $ 0.75    $ & $ 1.10594 $ & $ -1.35624 $   & $ -4.94074    $& $ 0 $ & $ 0.2 $& $ -15 $ &  $ 5 $  & $0$ & $1600$ & $1$\\
			RP5 & $ 1       $ & $ 1\cdot 10^{-14}  $ & $ 0 $ & $ 0 $  & $ 0 $ & $ 0 $& $ -0.5 $ &  $ 0.5 $  & $0$ & $300$ & $0.075$\\
			\hline 
		\end{tabular} 
	\end{center}
\end{table}
The first Riemann problem assumes a constant bottom bathymetry and an initial jump in the water height. The results obtained for the free surface elevation and the horizontal velocity are depicted in Figure~\ref{fig.RP1}, we observe that both the left shock and the right rarefaction are correctly captured.

\begin{figure}
	\centering
	\includegraphics[width=0.49\linewidth]{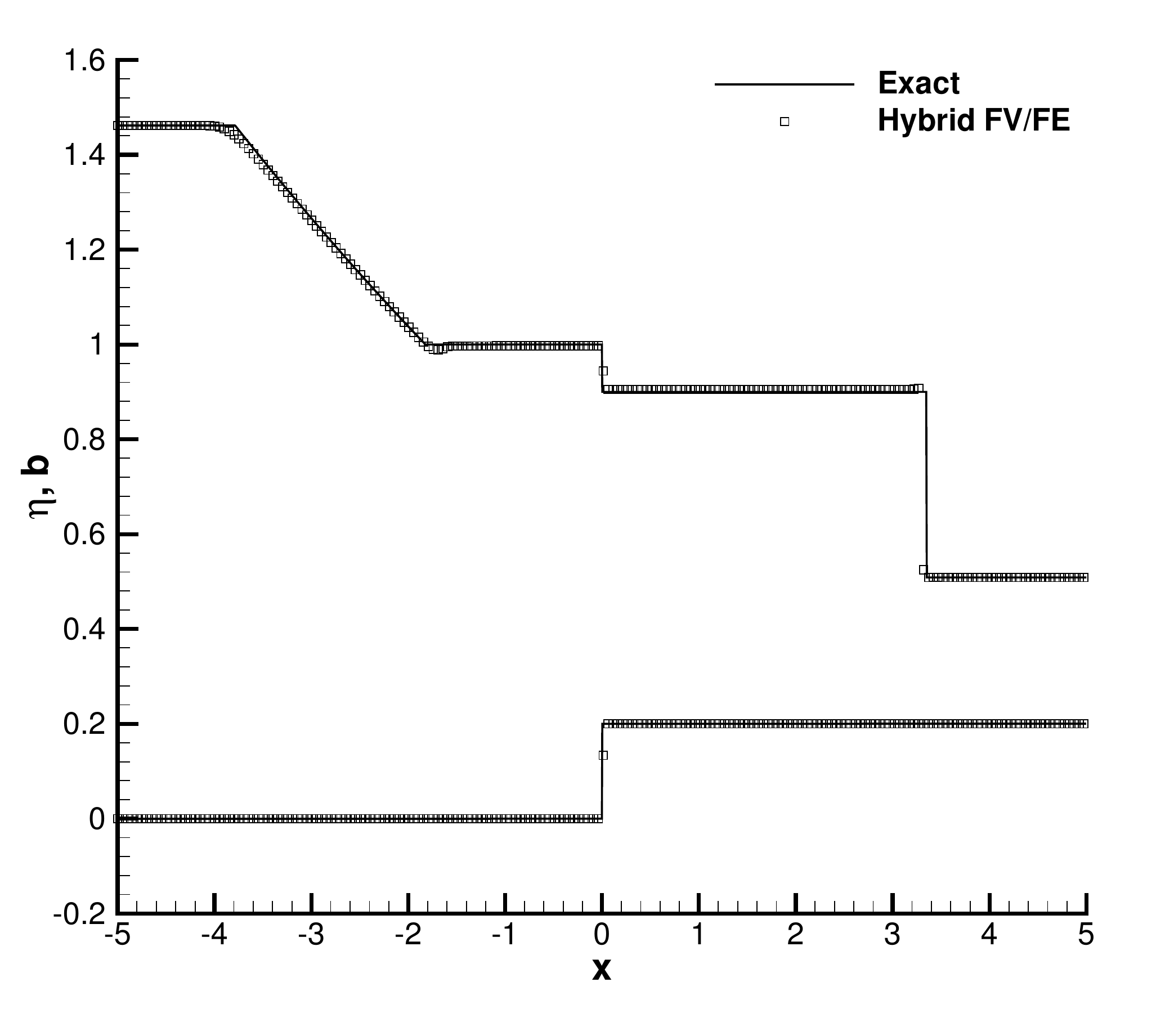} 
	\includegraphics[width=0.49\linewidth]{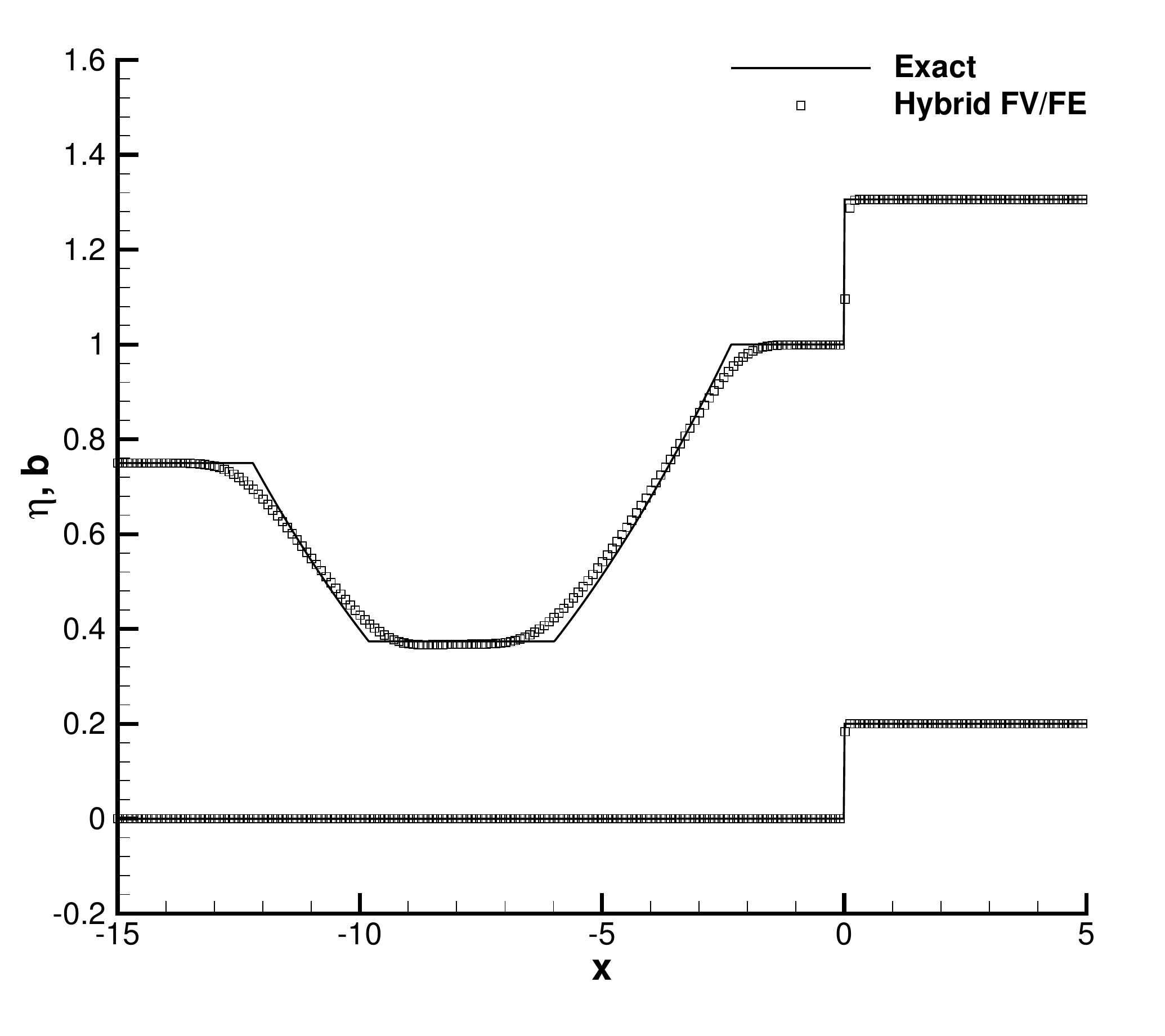}\\ 
	\includegraphics[width=0.49\linewidth]{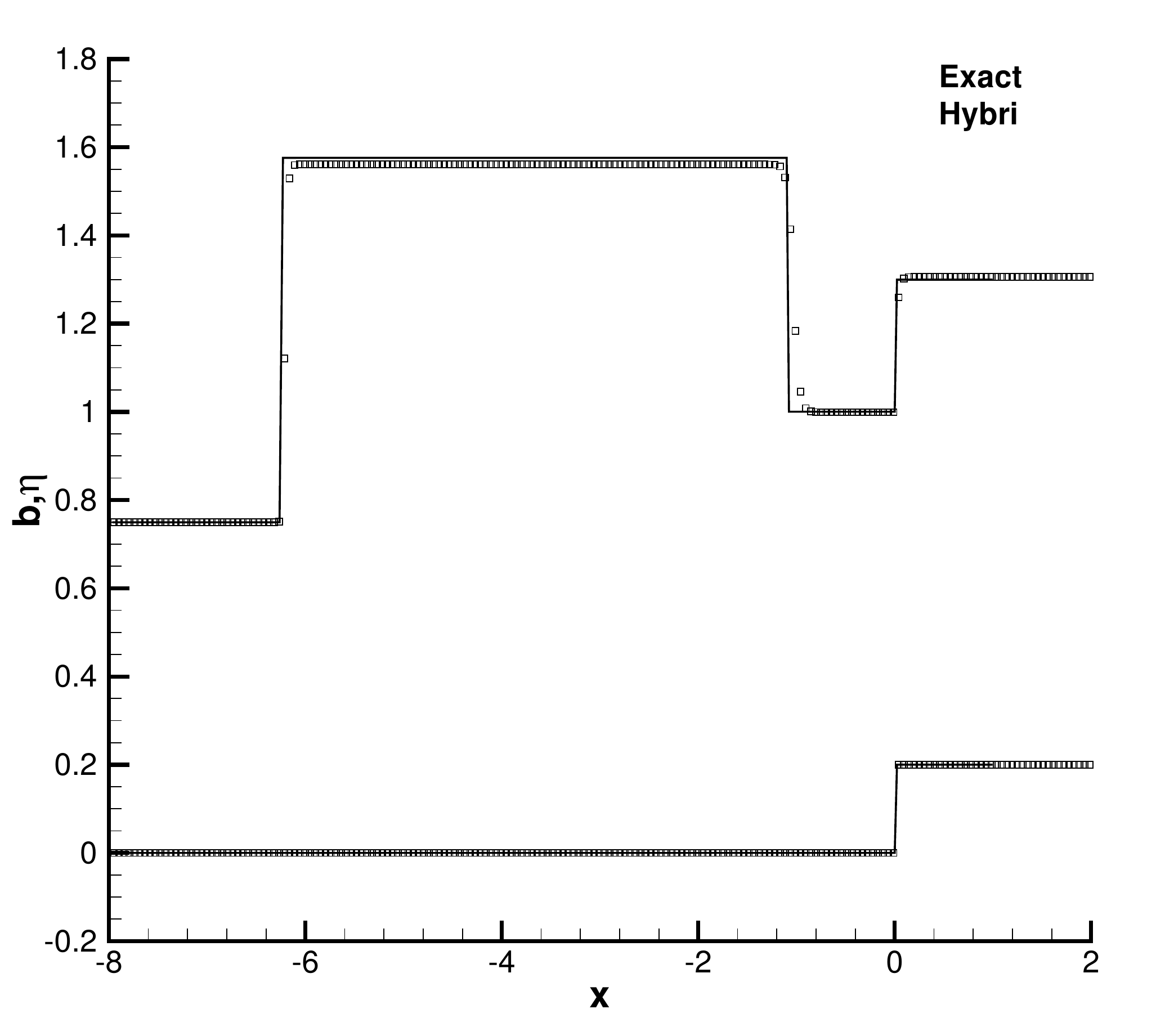}
	\includegraphics[width=0.49\linewidth]{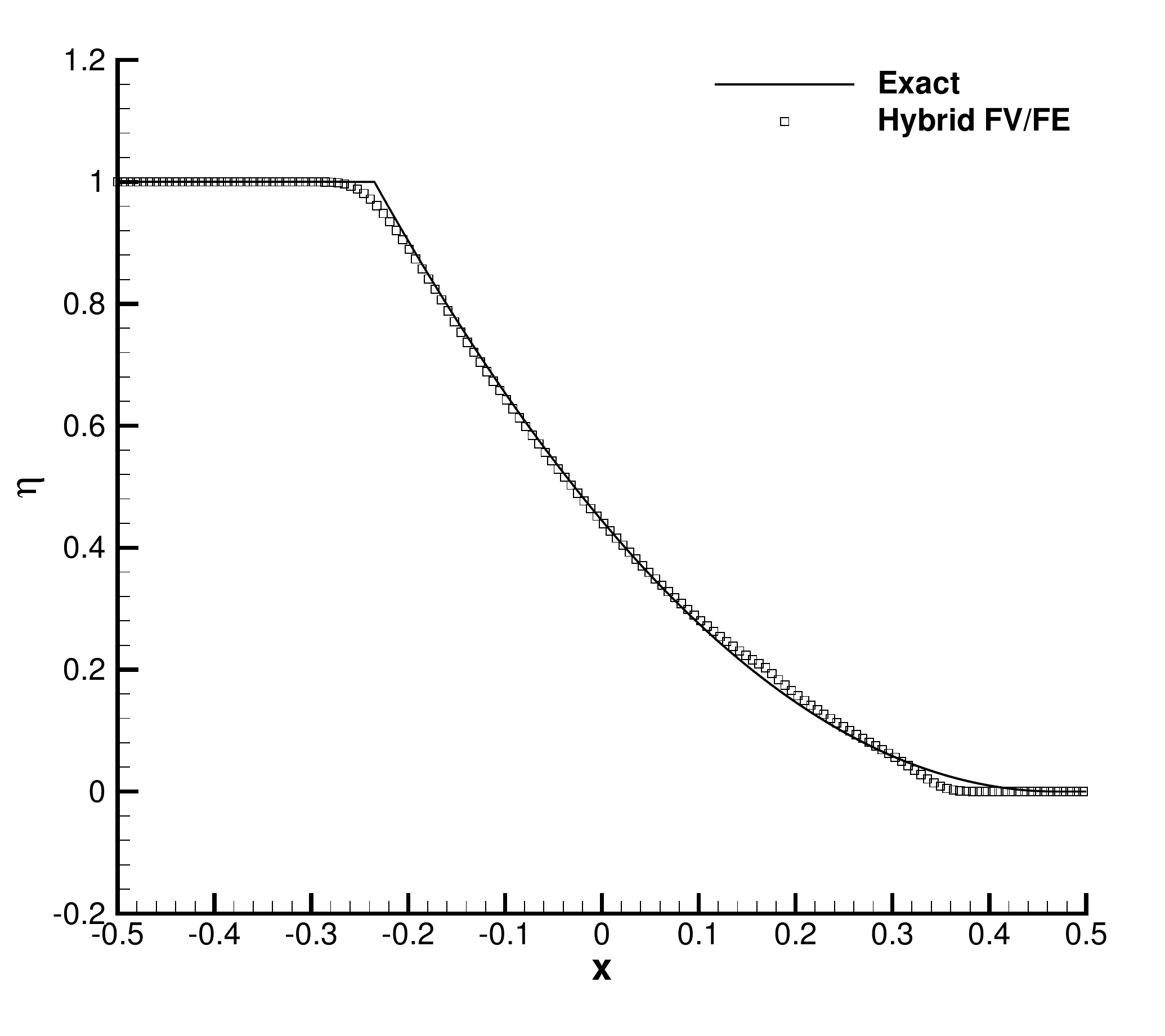}
	\caption{Free surface elevation and bottom bathymetry of the Riemann problems computed using the hybrid FV/FE method (squares) and reference solution (continuous line). From top left to bottom right: RP2 , LADER-ENO scheme; RP3, RP4, {RP5, LADER-ENO FV scheme with Rusanov dissipation in the FE algorithm.}}
	\label{fig:RP12}
\end{figure}

The free surface obtained for the remaining Riemann problems using the new hybrid FV/FE methodology is reported in Figure~\ref{fig:RP12}. RP2, RP3 and RP4 assume a constant in time and variable in space bathymetry, with a bottom step of height $\Delta b = 0.2$ located at $x_{c}=0$. In RP2 we perform a dambreak-type test, while the third and fourth tests are negative supercritical  motions with two rarefactions and two shock waves, respectively. Different setups have been considered to address the tests. RP2 has been run using the LADER-ENO scheme. 
RP3 requires for the activation of the Rusanov-type dissipation in the FE algorithm and the use of an artificial viscosity in the FV scheme so that we avoid spurious oscillations in the numerical solution. 
RP4 has also been simulated considering $c_{\alpha}=1$ and switching on the Rusanov-type dissipation in the FE algorithm. 
All three test show an excellent agreement with the exact solution; the moving wave fronts are properly captured, as well as the stationary shocks at the bottom step and the post shock states. 

The last Riemann problem, RP5, corresponds to a one-dimensional dambreak over a dry bottom, the exact solution of which has first been found by Ritter in \cite{Ritter}.  
{Figure~\ref{fig:RP12} report the results obtained using the second order ADER-ENO scheme for the nonlinear convective terms jointly with the Rusanov dissipation in the FE algorithm}.  
From the numerical point of view, to deal with wet and dry fronts, we compute the velocity field $\vel$ from the conservative momentum $\q$ as $\vel = \frac{h \q}{h^2+\epsilon}$ with $\epsilon = 10^{-7}$ to avoid division by zero when passing from conservative to primitive variables in dry areas of $\Omega$. {Moreover, an extra limiting strategy is applied to avoid artificial shocks in wet and dry fronts. Consequently, for small values of the water depth the first order scheme in space is retrieved.}

\subsection{Circular dambreak}
\begin{figure}
	\centering
	\includegraphics[width=0.49\linewidth]{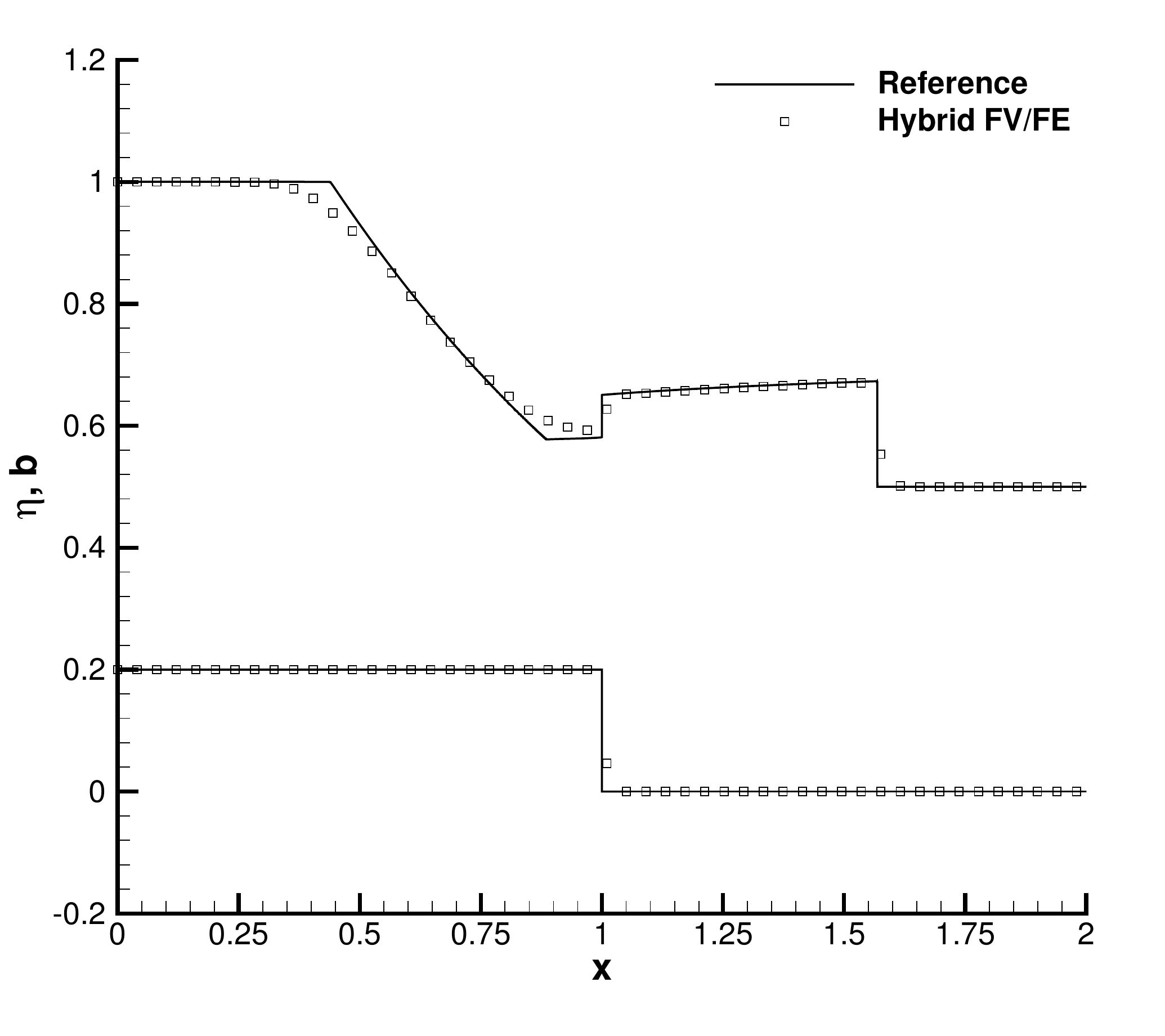}
	\includegraphics[width=0.49\linewidth]{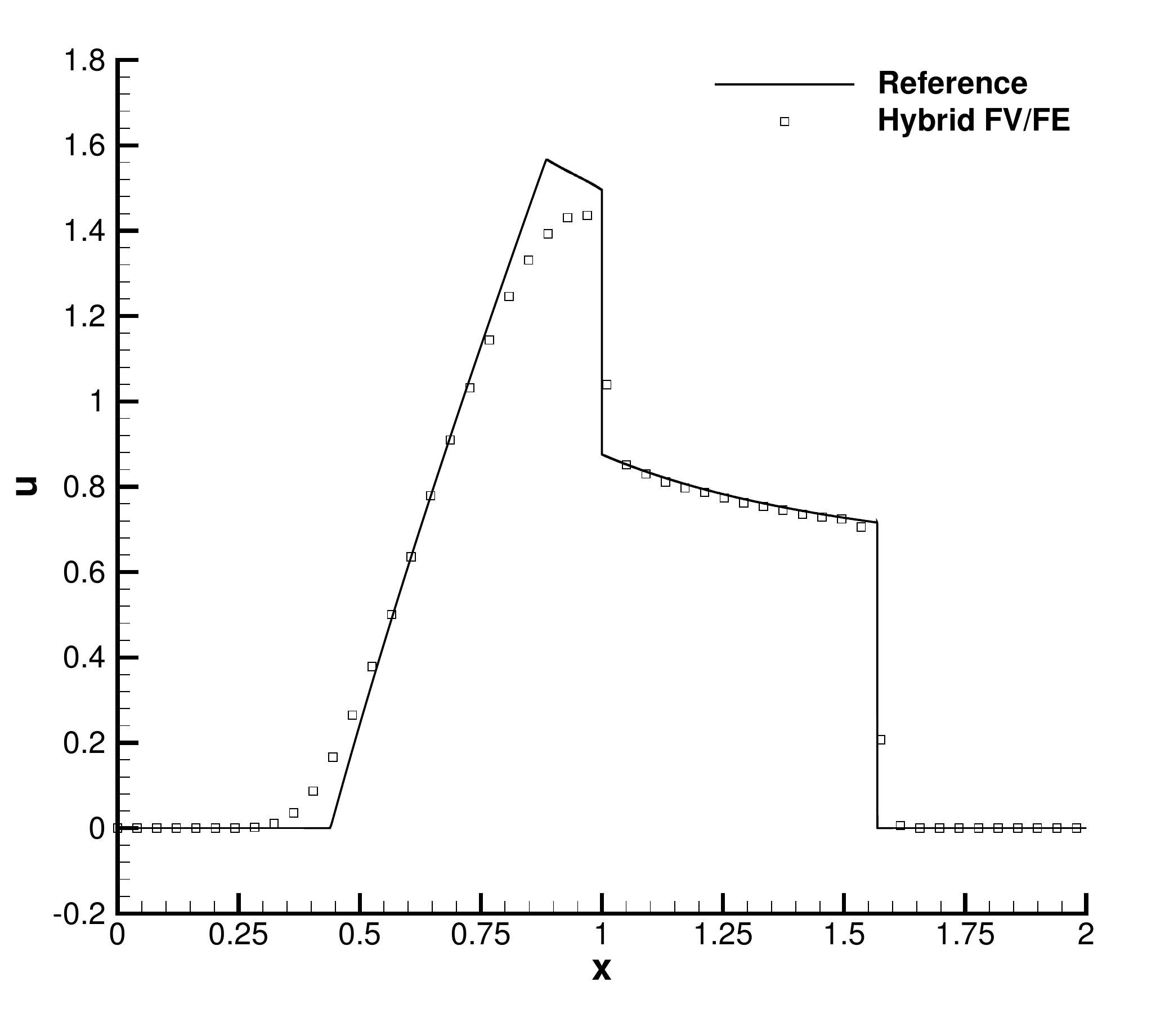}
	\caption{Circular dambreak. Comparison of the numerical (squares) and reference (continuous line) solutions obtained for the free surface and bottom profile (left) and velocity field (right) obtained using the hybrid FV/FE scheme with Rusanov dissipation in the FE algorithm.}
	\label{fig:CD}
\end{figure}

\begin{figure}
	\centering
	\includegraphics[trim=10 10 20 10,clip,width=0.49\linewidth]{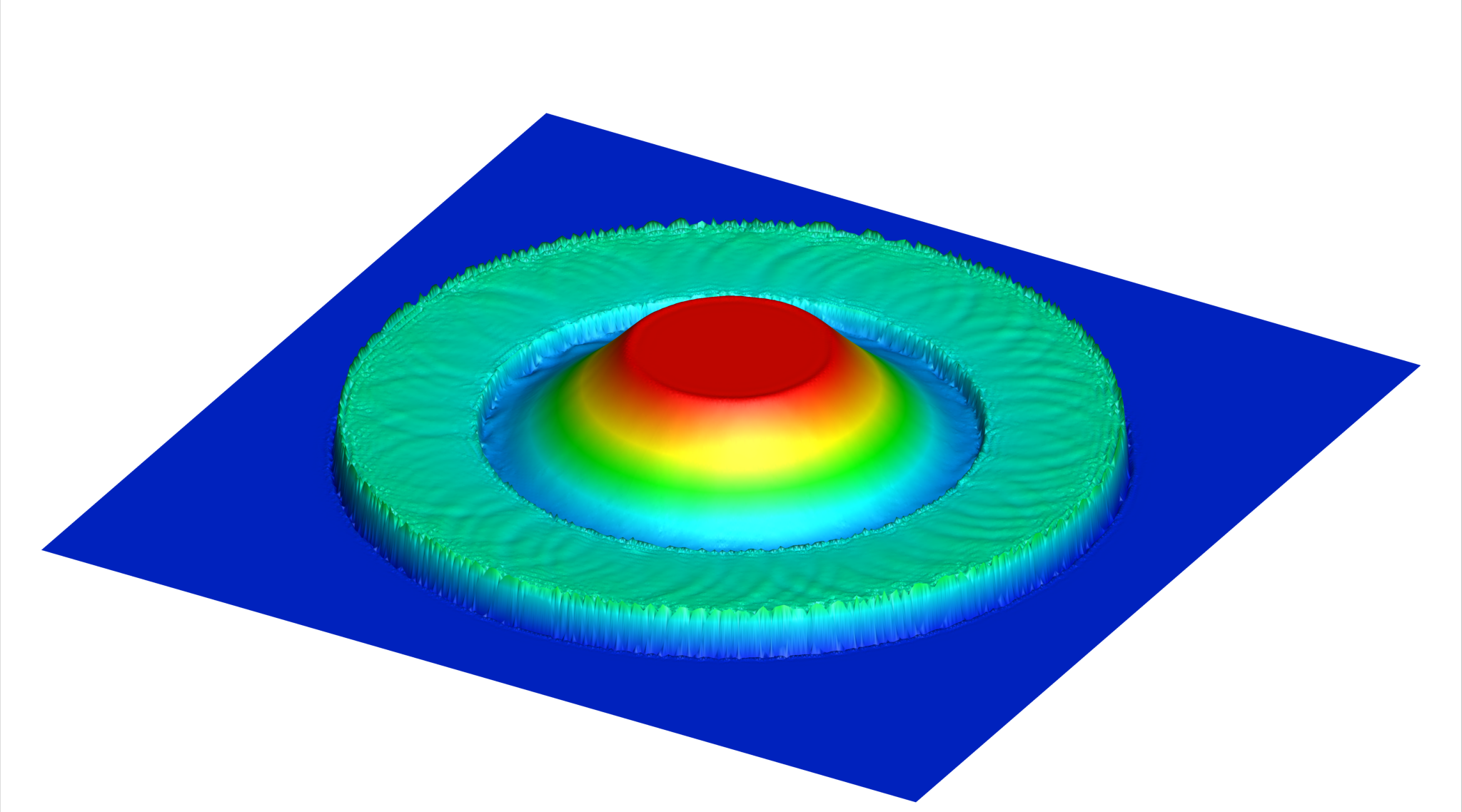}
	\includegraphics[trim=10 10 20 10,clip,width=0.49\linewidth]{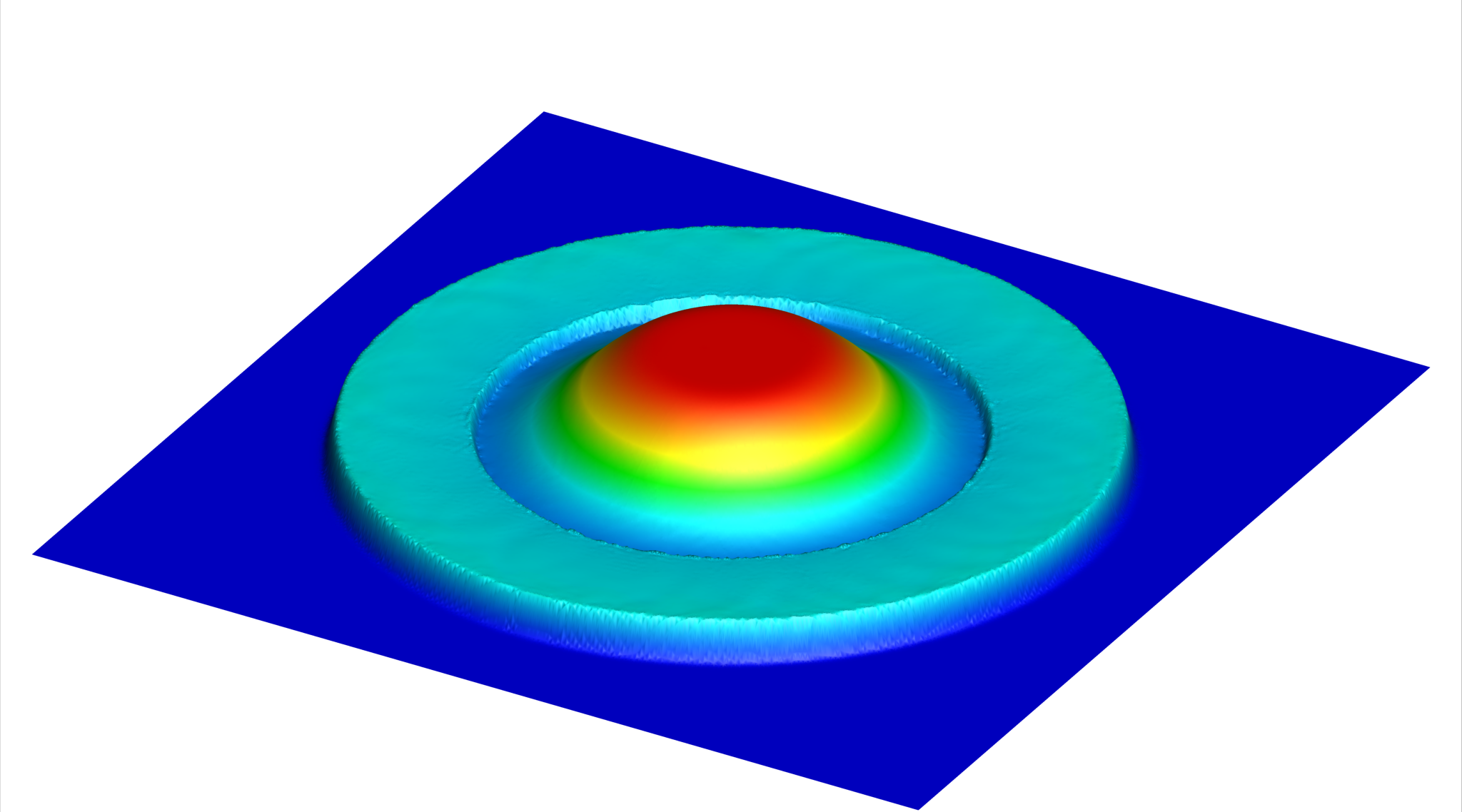}
	\caption{Circular dambreak. Free surface obtained at $T=0.2$ using the hybrid FV/FE scheme without Rusanov dissipation in the projection stage (left) and activating it (right).}
	\label{fig:CD3DRcomp}
\end{figure}

The circular dambreak over a bottom step defined in \cite{DC13,TD14sw} is employed to further study the robustness of the developed methodology in addressing 2D problems with initial discontinuities. We consider the computational domain $\Omega=[-2,2]^{2}$ and the following circular Riemann problem as initial condition
\begin{equation*}
	h(\x,0) = \left\lbrace \begin{array}{lr}
		0.8 & \mathrm{if }\, \left\| \x\right\| \leq 1, \\
		0.5 & \mathrm{if }\, \left\| \x\right\| \ge 1,
	\end{array} \right. \qquad
	b(\x,0) = \left\lbrace\begin{array}{lr}
		0.2 & \mathrm{if }\, \left\| \x\right\| \leq 1, \\
		0 & \mathrm{if }\, \left\| \x\right\| \ge 1,
	\end{array}\right.\qquad
	\vel(\x,0) = \mathbf{0}.
\end{equation*}
The simulation is run on a mesh composed of $106606$ primal elements using the hybrid FV-FE methodology with the LADER-ENO scheme for the nonlinear convective terms and setting $\theta=0.5$. In Figure \ref{fig:CD3DRcomp} we show the 3D plots of the free surface with and without Rusanov dissipation active in the FE algorithm. As expected, the extra dissipation avoids the spurious oscillations arising in the vicinity of discontinuities at the cost of an increase in the numerical diffusion.
A 1D cut through the numerical solution of the free surface and the velocity field obtained with Rusanov-type dissipation at time $T=0.2$ along $y=0$ is depicted in Figure \ref{fig:CD} together with the numerical reference solution computed again using the second order TVD-FV scheme employed in the well balance test, Section \ref{sec:wellbalance}. For details on how to obtain this numerical reference solution in cylindrical coordinates, see \cite{toro-book-swe}. The obtained numerical results show a good agreement between the new hybrid FV/FE methodology and the reference solution. 

\subsection{Low Froude number flow around a circular cylinder}
\begin{figure}
	\centering
		\includegraphics[width=0.49\linewidth]{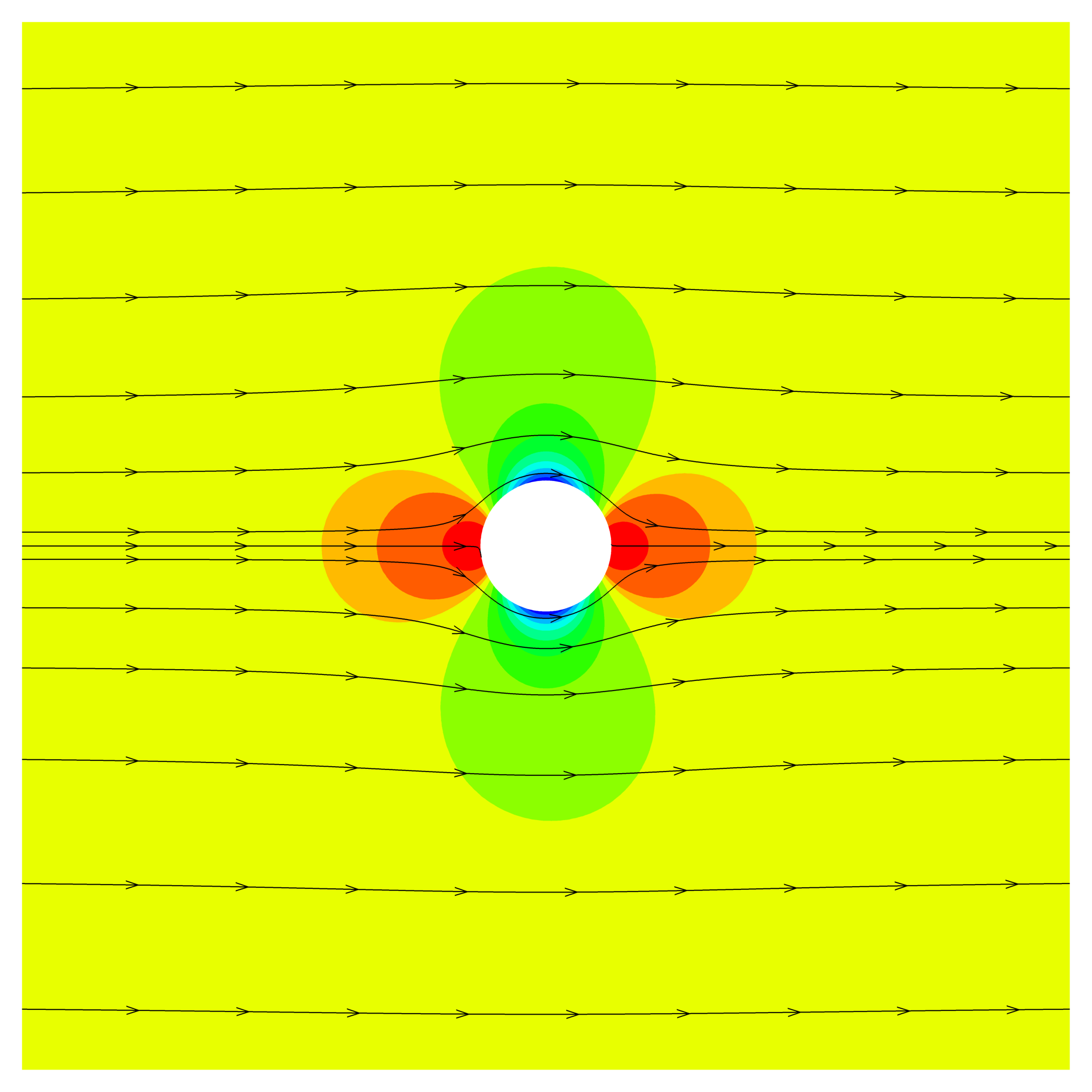}
	\caption{Low Froude number flow around a circular cylinder. { Zoom of the} contour plot of the free surface elevation and streamlines around the cylinder obtained using the hybrid FV/FE LADER scheme ($t= 10$, CFL$=0.5$).}
	\label{fig:cylinder}
\end{figure}

\begin{figure}[h]
	\centering
	\includegraphics[width=0.49\linewidth]{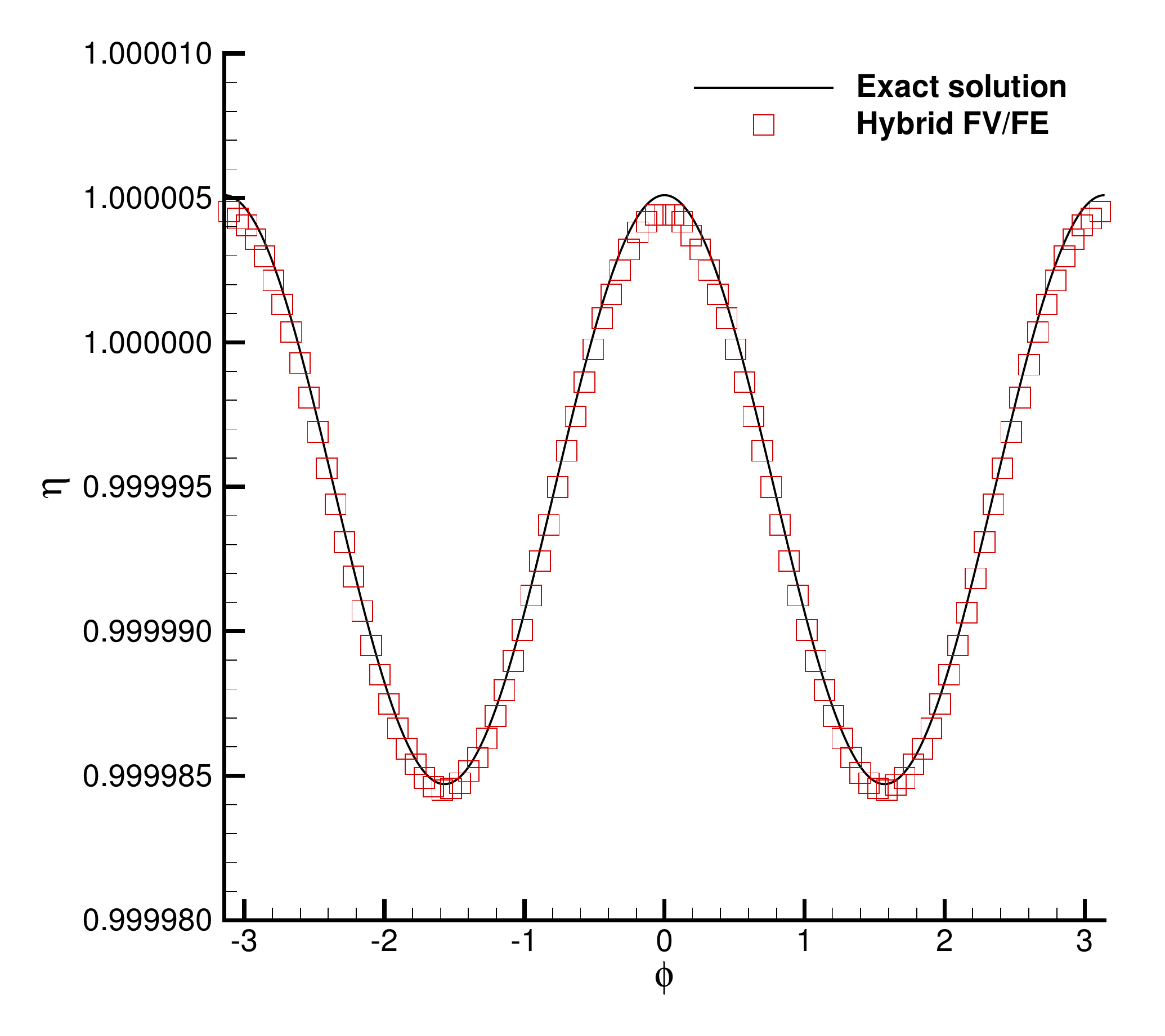}
	\includegraphics[width=0.49\linewidth]{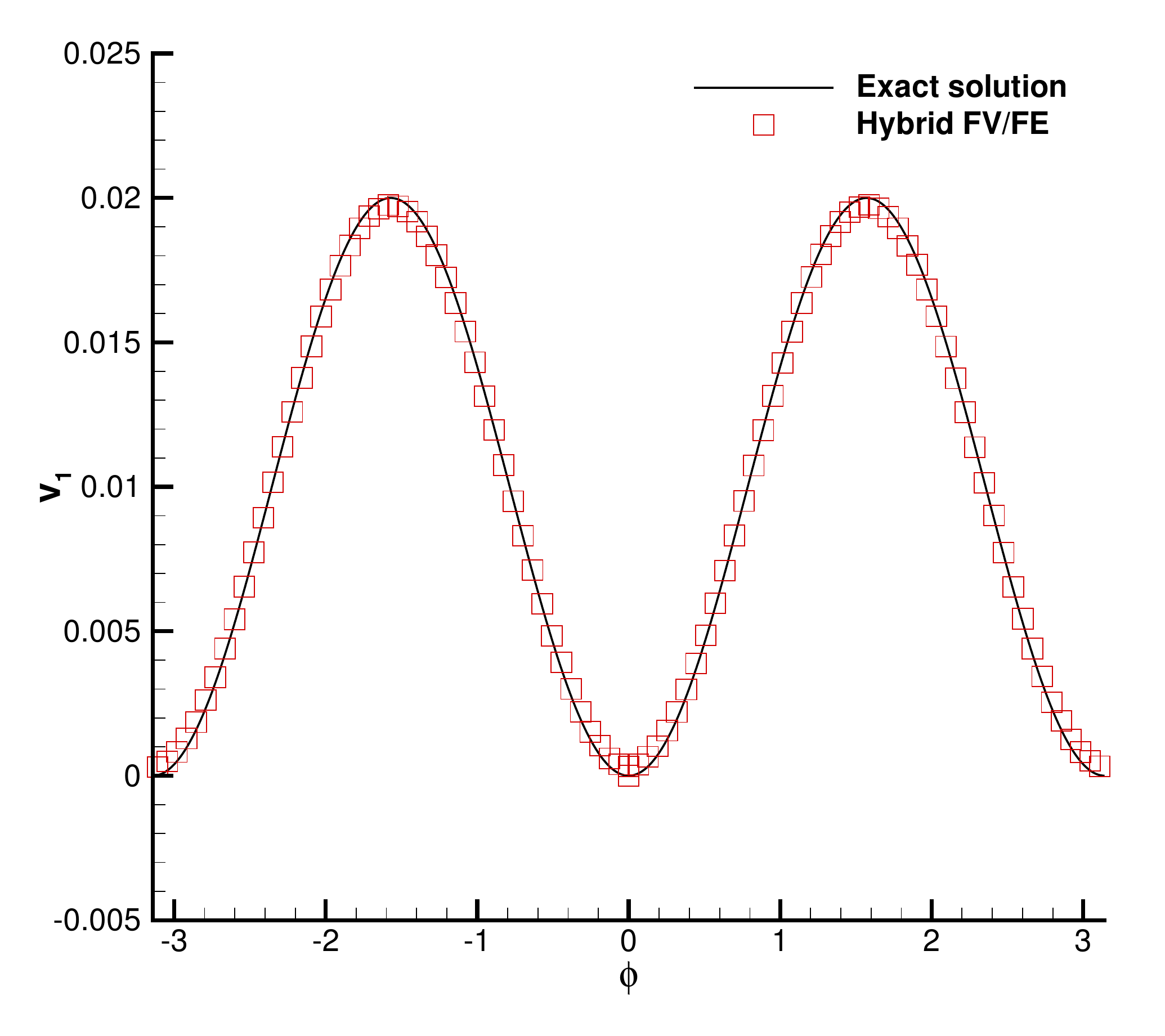}
	\caption{Flow over a cylinder. Comparison of the numerical (red squares) and exact (black line) solutions obtained for the free surface (left) and velocity field (right) at radius $r=1.01$.}
	\label{fig:cylinder_comparison}
\end{figure}

In this section, we analyse the results obtained with the proposed hybrid FV/FE method for the solution of a steady low Froude number flow around a circular cylinder, see \cite{TD14sw} for details on the setup of this test. The computational domain considered here is 
\mbox{$\Omega = { \left[ -8,8\right]^{2}} \setminus \left\lbrace \x=(x,y) \in\mathbb{R}^{2}, \, \left\| \x \right\|^{2}\leq r_c = 1 \right\rbrace$}, the equilibrium free surface is set to $\eta_{0} = 1$ { and the bottom elevation is $b=0$}. 
A uniform horizontal velocity, $\vel_{m} = (v_m,0)^T$ with $v_m=10^{-2}$, is imposed at the left, top and bottom boundaries of the domain, while outflow boundary conditions are set on the right boundary. Consequently, the corresponding Froude number is Fr$ = \frac{\vels_{m}}{\sqrt{g { h_{0}}}} = 3.19\cdot 10^{-3}$. The nature of this test problem allows the computation of an analytical solution of the corresponding potential flow as 
\begin{equation*} 
	r   = \left\| \x \right\|, \qquad
	\tan(\phi)  = \frac{y}{x}, \qquad
	\vels_{r}     =  \vels_{m} \left(1-\frac{r_c^2}{r^2}\right)\cos(\phi), \qquad
	\vels_{\theta} = -\vels_{m} \left(1+\frac{r_c^2}{r^2}\right)\sin(\phi)
\end{equation*}
with $\vels_{r}$ and $\vels_{\theta}$ the radial and the angular velocities, respectively. Moreover, taking into account Bernouilli's equation we also get the exact free surface elevation as 
\begin{equation*}
	\eta = \eta_{0} + \frac{1}{2} \vels_{m}^{2} g \left(2 \frac{r_c^2}{r^2} \cos(2\phi)-\frac{r_c^4}{r^2}\right).
\end{equation*} 
The numerical solution is computed on a primal grid made of $30384$ triangular elements, using the second order LADER scheme for the nonlinear convective terms and setting $\theta=1$ in the time discretization. No nonlinear limiting strategy is needed in this test thanks to the smoothness of the solution. The initial condition for the velocity field is set to the exact solution, thus avoiding the generation of strong initial transient waves. In Figure \ref{fig:cylinder}, we depict the contour plot of the free surface together with the streamlines of the velocity field obtained at $T=10$, when the stationary state has already been reached.
For comparison against the exact solution, Figure  \ref{fig:cylinder_comparison} shows the values of the free surface and the velocity obtained along the circumference of radius $r=1.01$ centred at the origin. A good agreement is observed for both fields.

\subsection{High Froude number flow around a blunt body}
\begin{figure}
	\begin{center}
		\begin{tabular}{cc} 
			\includegraphics[trim=10 10 20 10,clip,width=0.41\textwidth]{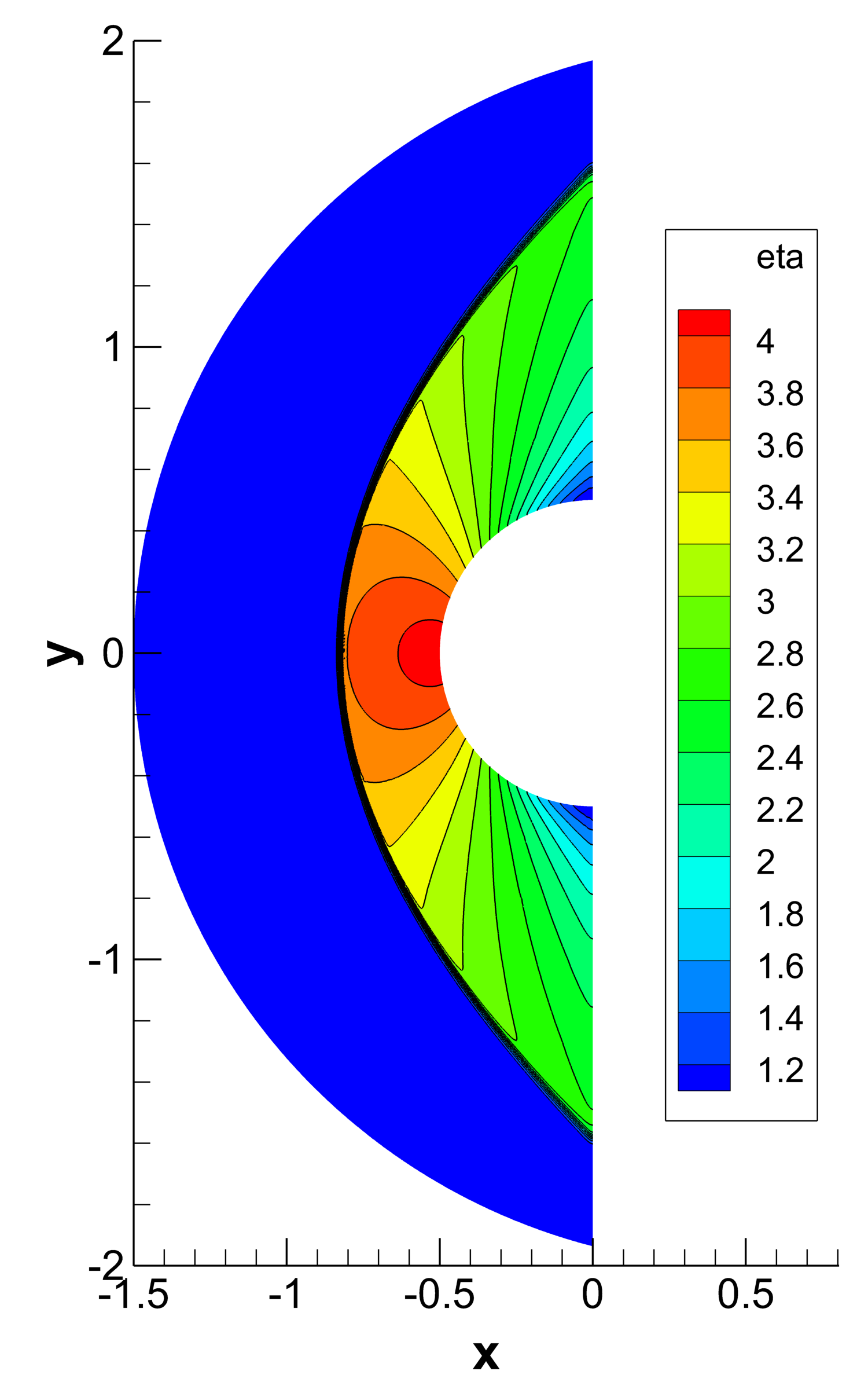}   & 
			\includegraphics[trim=10 10 20 10,clip,width=0.41\textwidth]{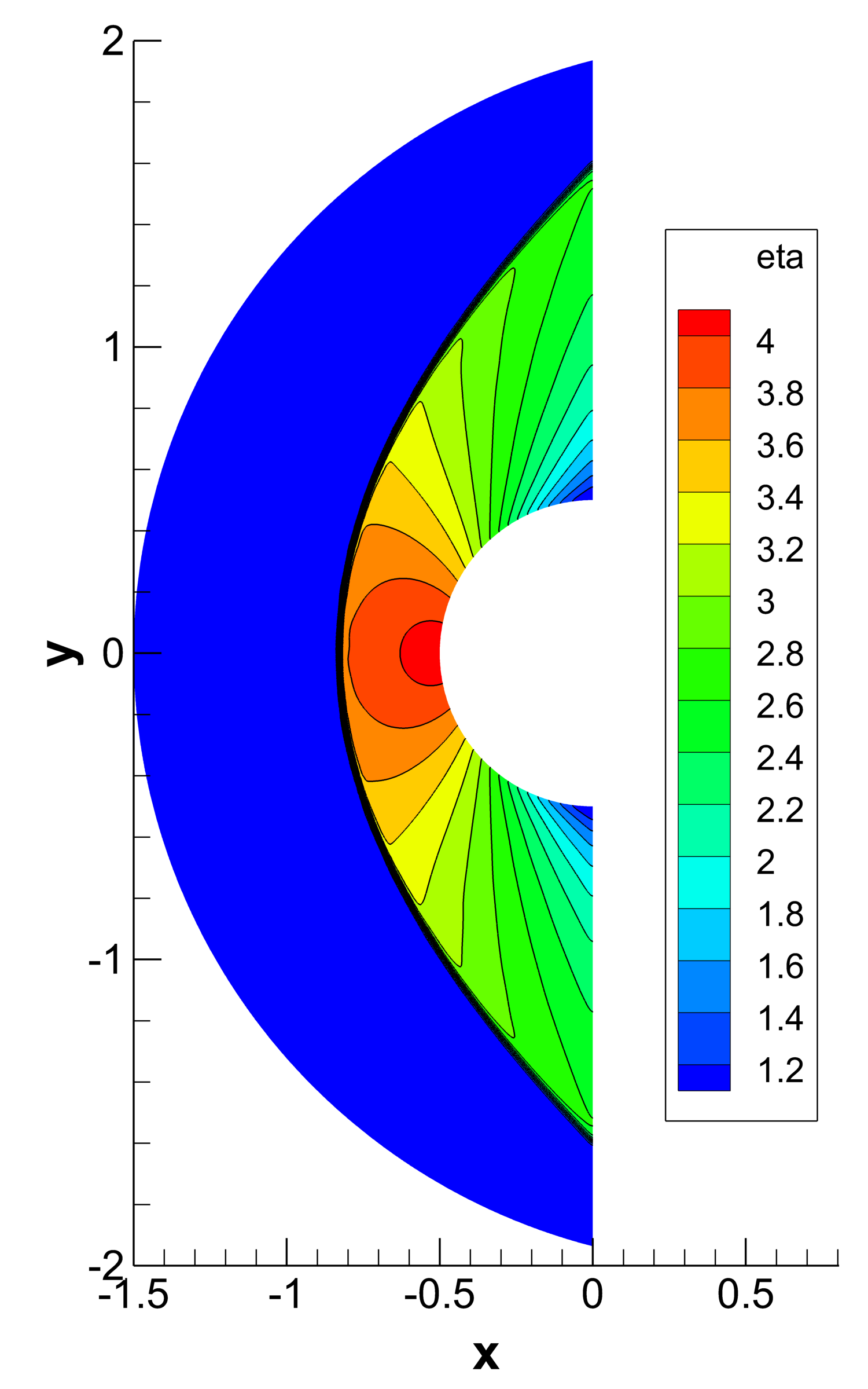}     
		\end{tabular}
		\caption{ Free surface elevation contour lines for the supercritical flow (Fr $=3$) over a blunt body at time $T=2$. Left: the new hybrid FV/FE scheme presented in this paper. Right: reference solution obtained with an explicit Godunov-type finite volume scheme. } 
		\label{fig.bb2d} 
	\end{center}
\end{figure}

\begin{figure}[!htbp]
	\begin{center}
		\begin{tabular}{cc} 
			\includegraphics[trim=10 10 20 10,clip,width=0.45\textwidth]{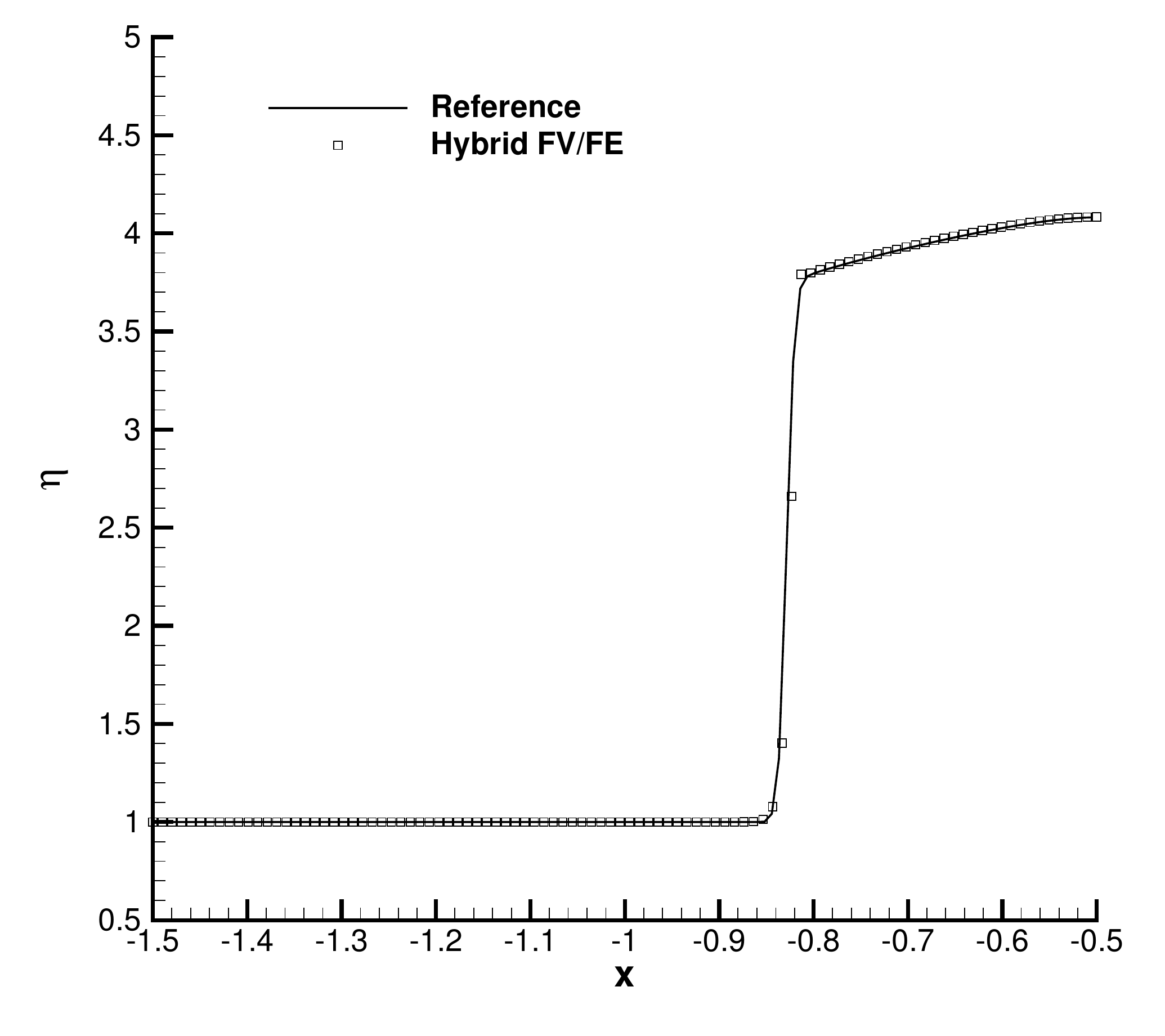}   & 
			\includegraphics[trim=10 10 20 10,clip,width=0.45\textwidth]{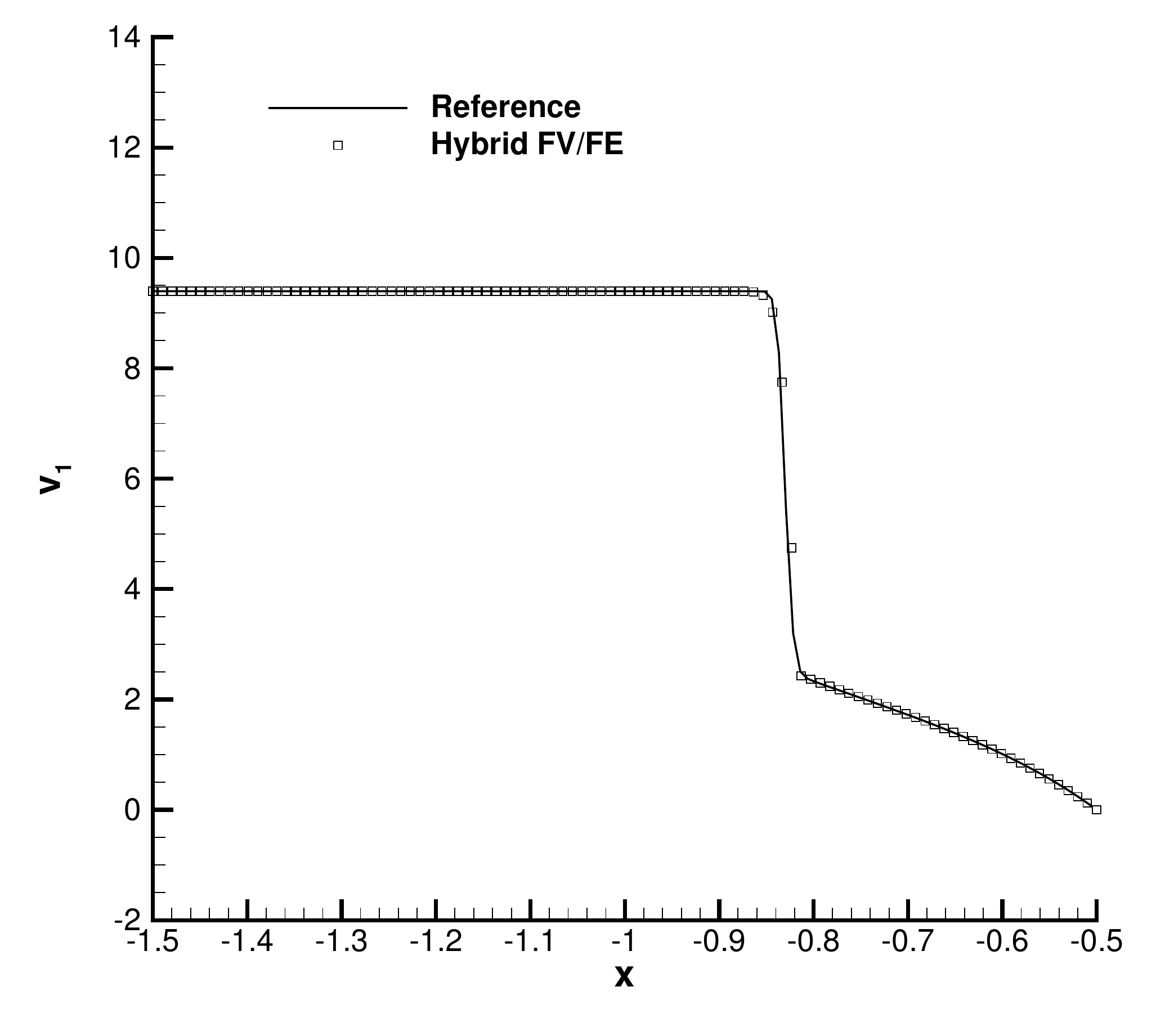}     
		\end{tabular}
		\caption{ Supercritical flow over a blunt body at time $t=2$. 1D cut through the numerical solution along the $x$ axis. Free surface elevation (left) and velocity component $v_1$ (right). } 
		\label{fig.bb2d.cut} 
	\end{center}
\end{figure}

To show the capability of the proposed method to deal also with high Froude number flows, we now focus on the simulation of a supercritical flow around a blunt body. We consider the computational domain defined by the intersection of the half-plane $x<0$ and the annulus centred at $\x_c=(0.5,0)$ with internal and external radius $r=0.5$, $R=2$ respectively, that is discretized using $348964$ triangles. The initial flow is characterized by a water depth $h=1$ and a velocity field $\vel=\left(3 \sqrt{gh},0\right)^{T}$, {corresponding to a Froude number of Fr$=3$.} To simulate the presence of the wall in the inner circumference we use symmetry boundary conditions, while inlet boundary conditions are set in the outer circumference and outflow boundary conditions are considered in the remaining boundaries. The simulation is run until $T=2$ when the stationary solution has already been reached. In the left plot of Figure \ref{fig.bb2d} we can observe the solution obtained using the LADER-ENO scheme with artificial viscosity $c_{\alpha}=2$ and activated Rusanov dissipation in the finite element part. The right subplot reports the reference solution obtained using a explicit Godunov-type explicit finite volume scheme \cite{Dambreak3Dexp}. To ease comparison also a 1D cut, along the $x$ axis, of the water depth and the velocity is depicted in Figure \ref{fig.bb2d.cut}. An excellent agreement between the new semi-implicit hybrid FV/FE scheme and the reference solution can be observed at the bow shock. 


\subsection{3D dambreak on a dry plane}
\begin{figure}
	\begin{center}
		\begin{tabular}{cc} 
			\includegraphics[width=0.5\textwidth]{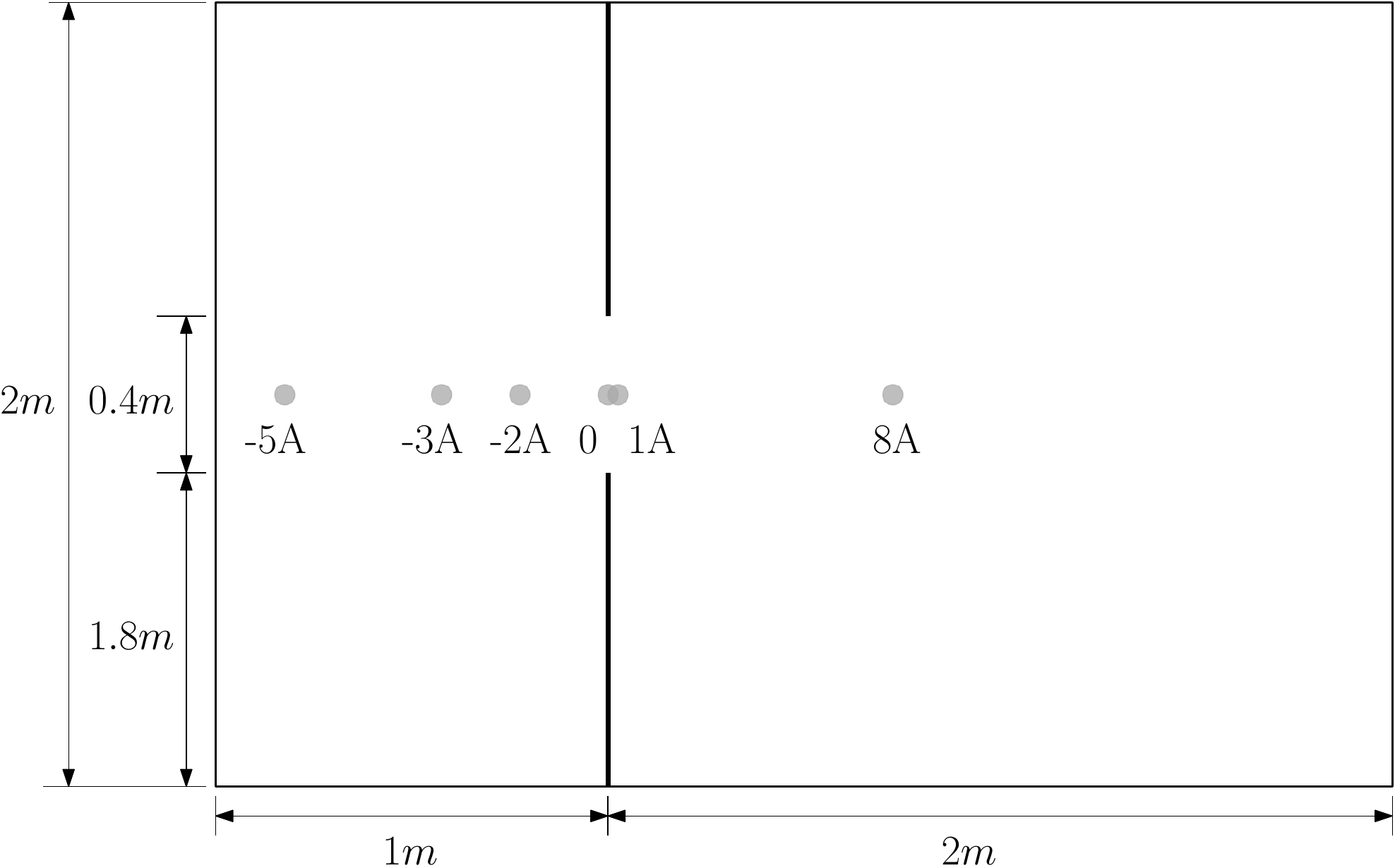}
		\end{tabular}
		\caption{Computational domain and wave gauges locations for the 3D dambreak test case. } 
		\label{fig.dbf.3d.domain}
	\end{center}
\end{figure}

To analyse the behaviour of the developed methodology in practical applications, we now consider the three-dimensional dambreak on a dry plane proposed in \cite{Dambreak3Dexp,Dambreak3D} and also considered in \cite{DIM3D}. This test is based on a laboratory experiment which considers a tank of dimensions $1m \times 2m \times 0.7m$ connected to a flat plate through a removable gate, see Figure \ref{fig.dbf.3d.domain}. The tank is initially filled with water at rest up to a height of $h=0.6m$. Then, the gate is opened in less than $0.1s$, hence simulating an instantaneous rupture of a dam.   Consequently, the water front flows downstream, while a rarefaction wave travels in the opposite direction inside the tank. 

\begin{figure}
	\begin{center}
		\begin{tabular}{cc} 
			\includegraphics[trim=10 10 20 10,clip,width=0.45\textwidth]{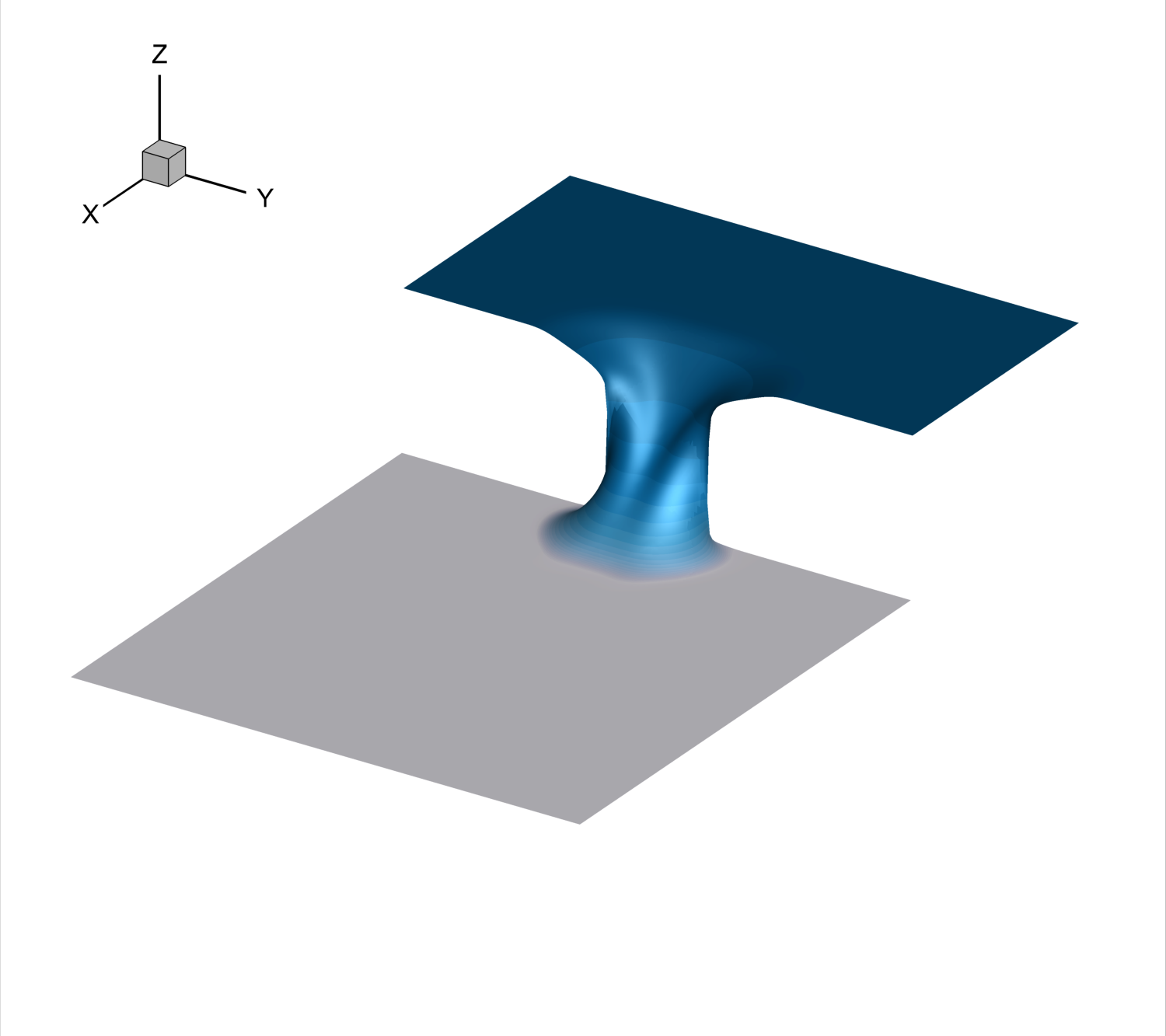}   & 
			\includegraphics[trim=10 10 20 10,clip,width=0.45\textwidth]{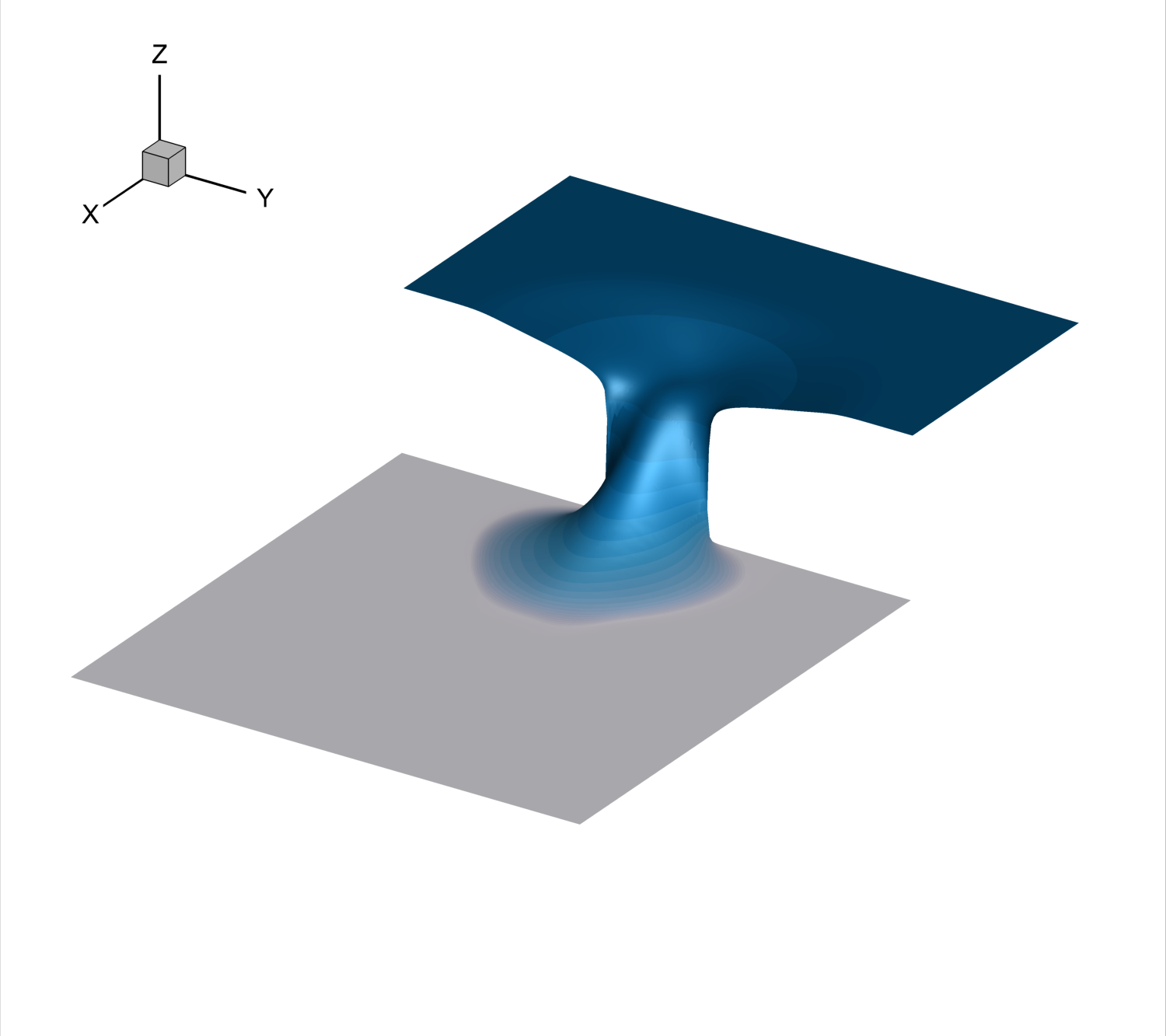} \\  
			\includegraphics[trim=10 10 20 10,clip,width=0.45\textwidth]{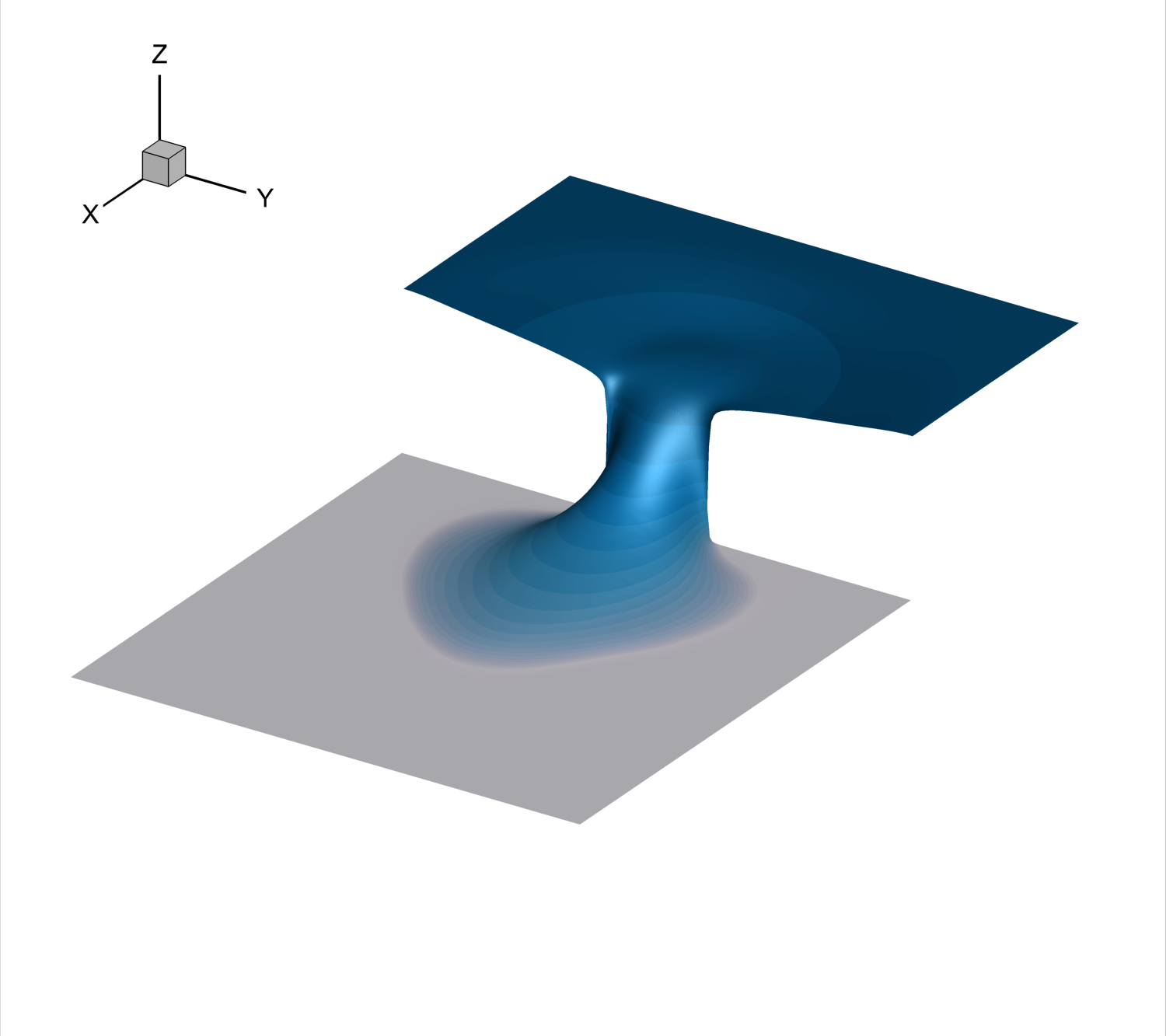} & 
			\includegraphics[trim=10 10 20 10,clip,width=0.45\textwidth]{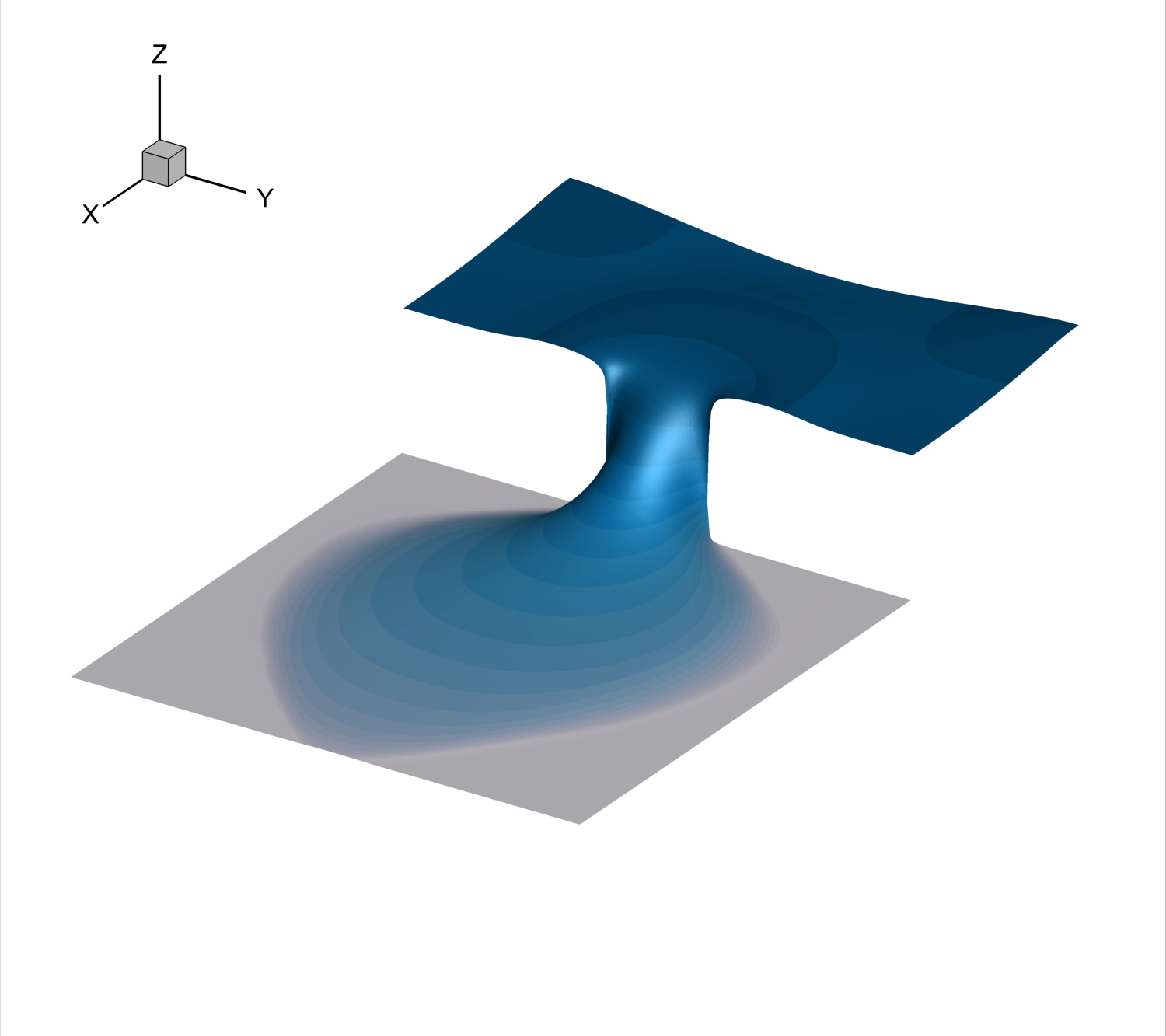} \\
			\includegraphics[trim=10 10 20 10,clip,width=0.45\textwidth]{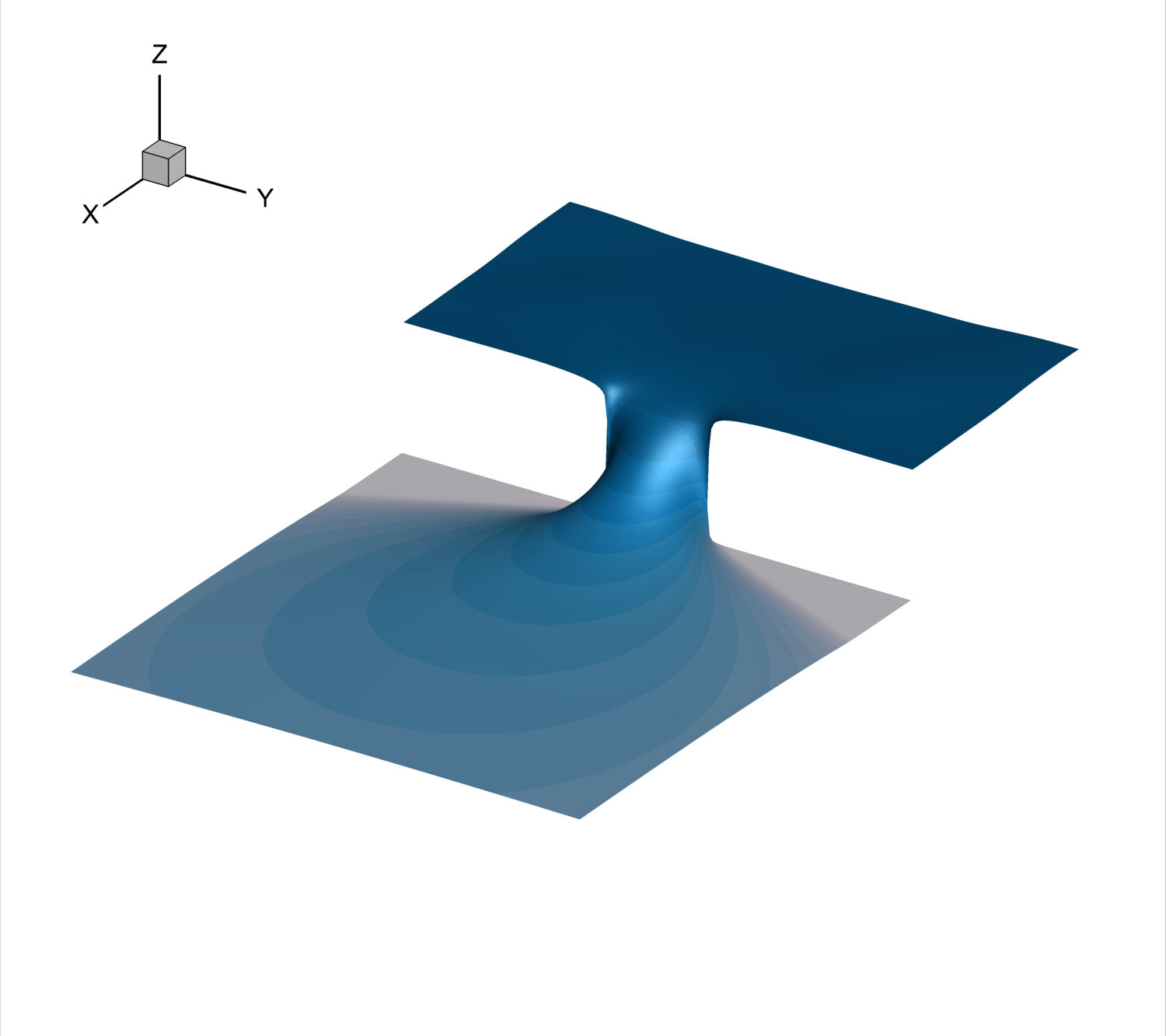}  & 
			\includegraphics[trim=10 10 20 10,clip,width=0.45\textwidth]{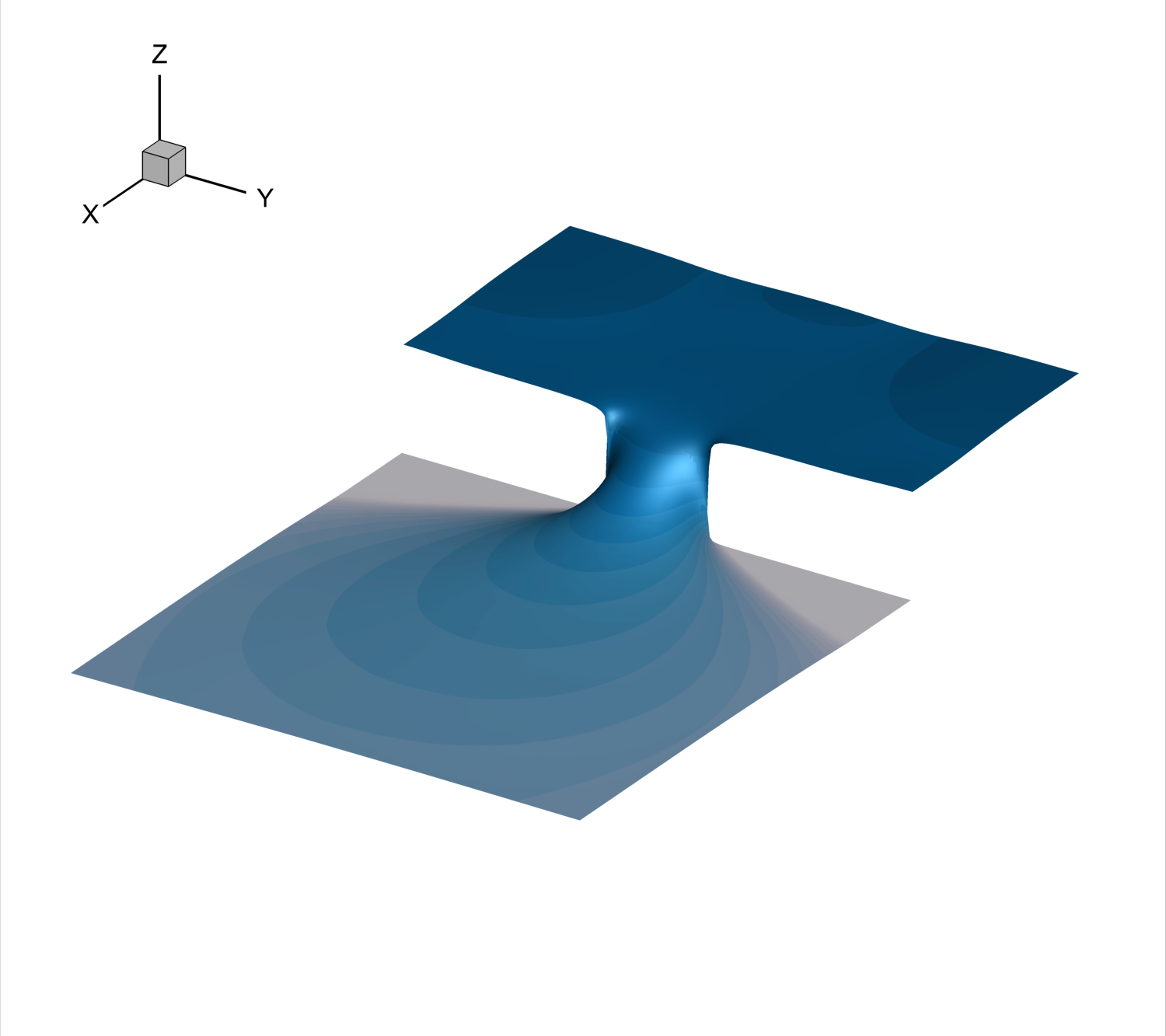}  
		\end{tabular}
		\caption{Temporal evolution of the 3D free surface profile  $\eta(\mathbf{x},t)$ obtained with the new hybrid FV/FE scheme at times $t=0.1$, $t=0.2$, $t=0.3$, $t=0.5$, $t=1.0$ and 
			$t=2.0$, from top left to bottom right. } 
		\label{fig.dbf.3d}
	\end{center}
\end{figure}

For the numerical simulation of this benchmark problem using the hybrid FV/FE methodology for the shallow water equations proposed in this paper, we define symmetry boundary conditions for the tank walls and free outflow in the three remaining boundaries delimiting the plane. Note that, in this test case, the friction coefficient plays and important role to correctly reproduce experimental data and is set to $n= 10^{-3}$. The 2D computational domain is discretized using a primal triangular mesh made of $135166$ cells. Regarding the time step computation a CFL number of $0.5$ has been defined 
and the Rusanov-type dissipation was used in the FE discretization to avoid spurious oscillations.  
The temporal evolution of the obtained free surface elevation can be seen in Figure \ref{fig.dbf.3d} for time instants $t\in\left\lbrace 0.1,0.2,0.3,0.5,1.0,2.0 \right\rbrace$. 

\begin{figure}
	\begin{center}
		\begin{tabular}{cc} 
			\includegraphics[width=0.45\textwidth]{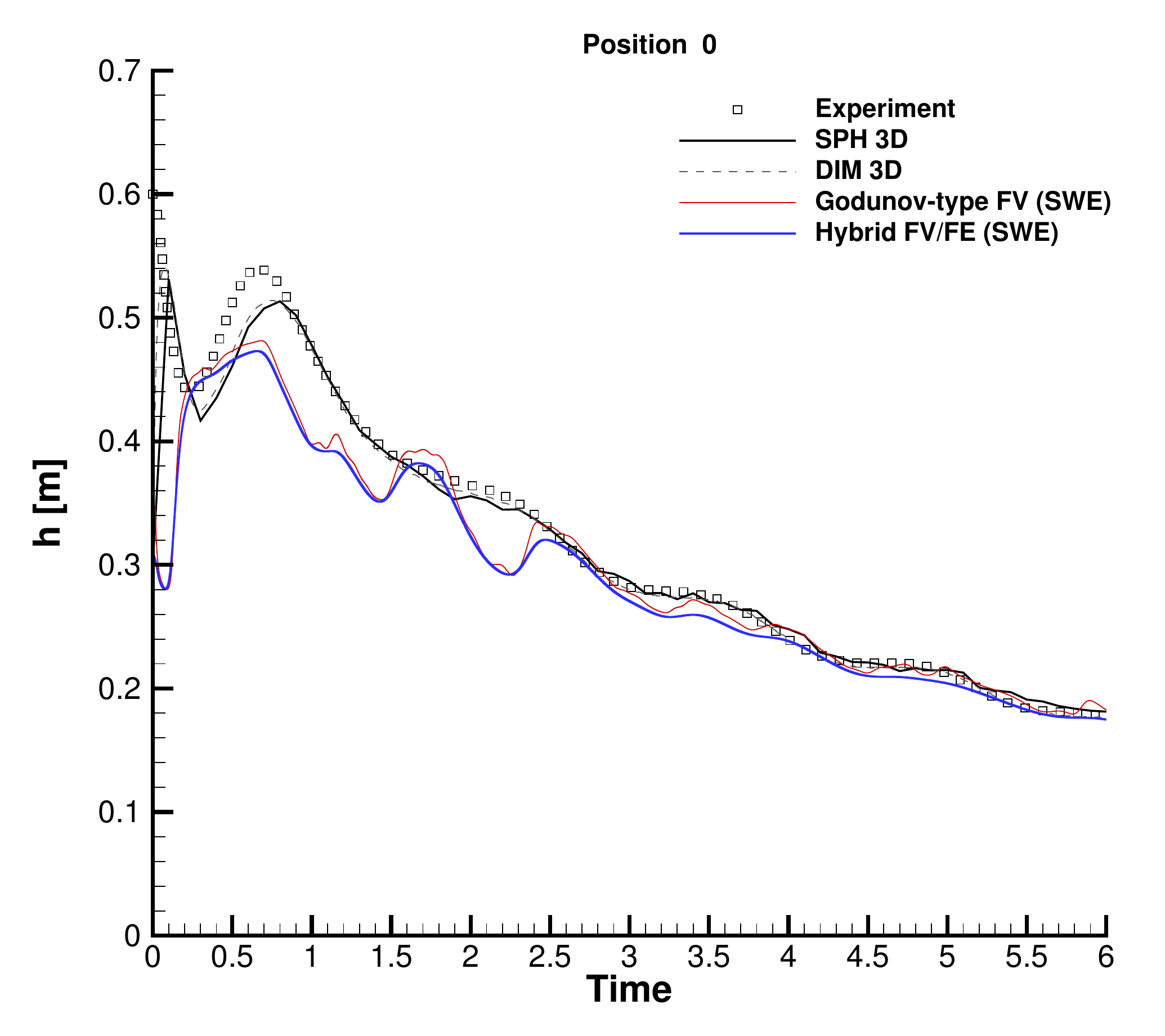}   & 
			\includegraphics[width=0.45\textwidth]{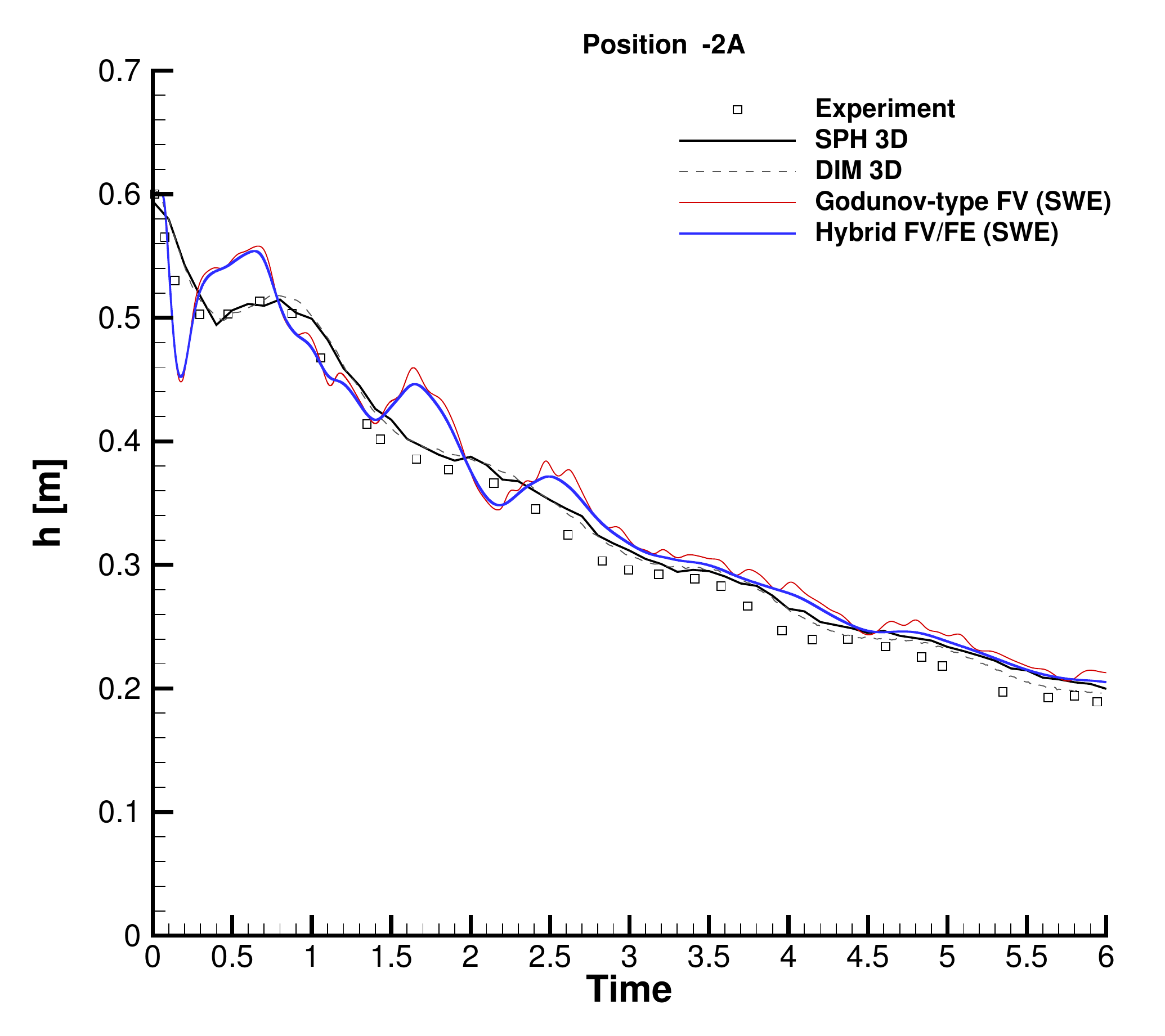} \\  
			\includegraphics[width=0.45\textwidth]{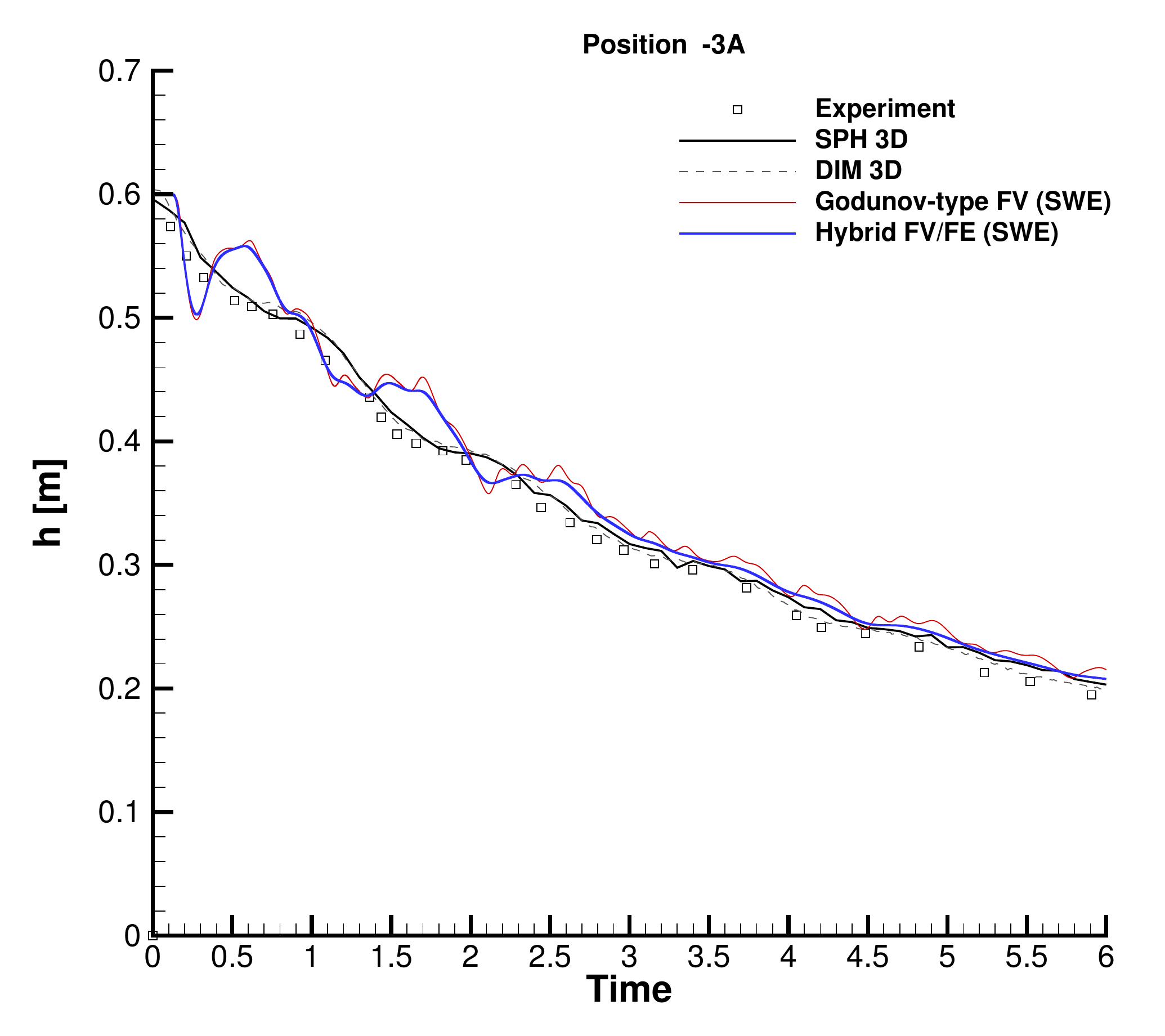} & 
			\includegraphics[width=0.45\textwidth]{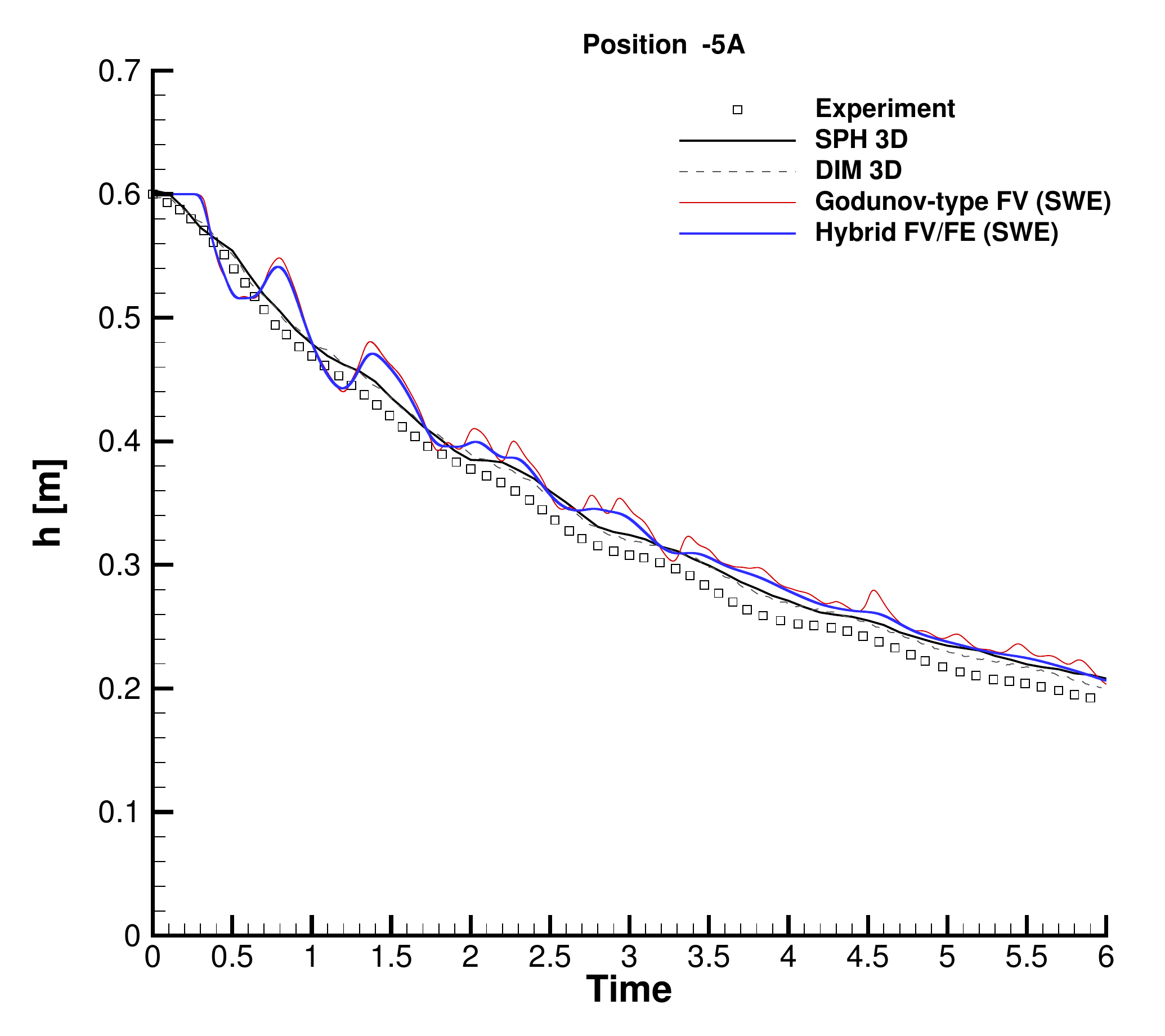} \\
			\includegraphics[width=0.45\textwidth]{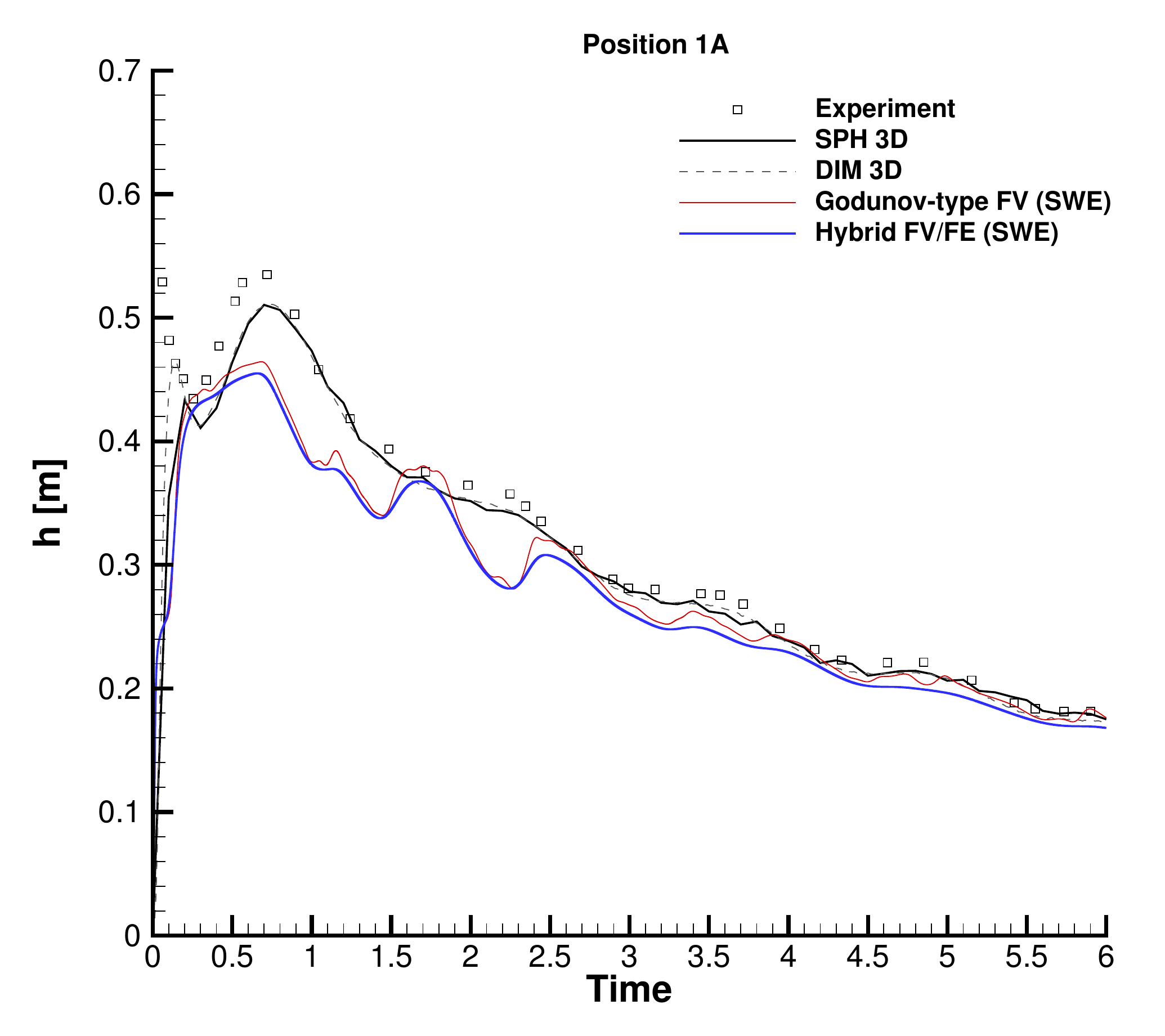}  & 
			\includegraphics[width=0.45\textwidth]{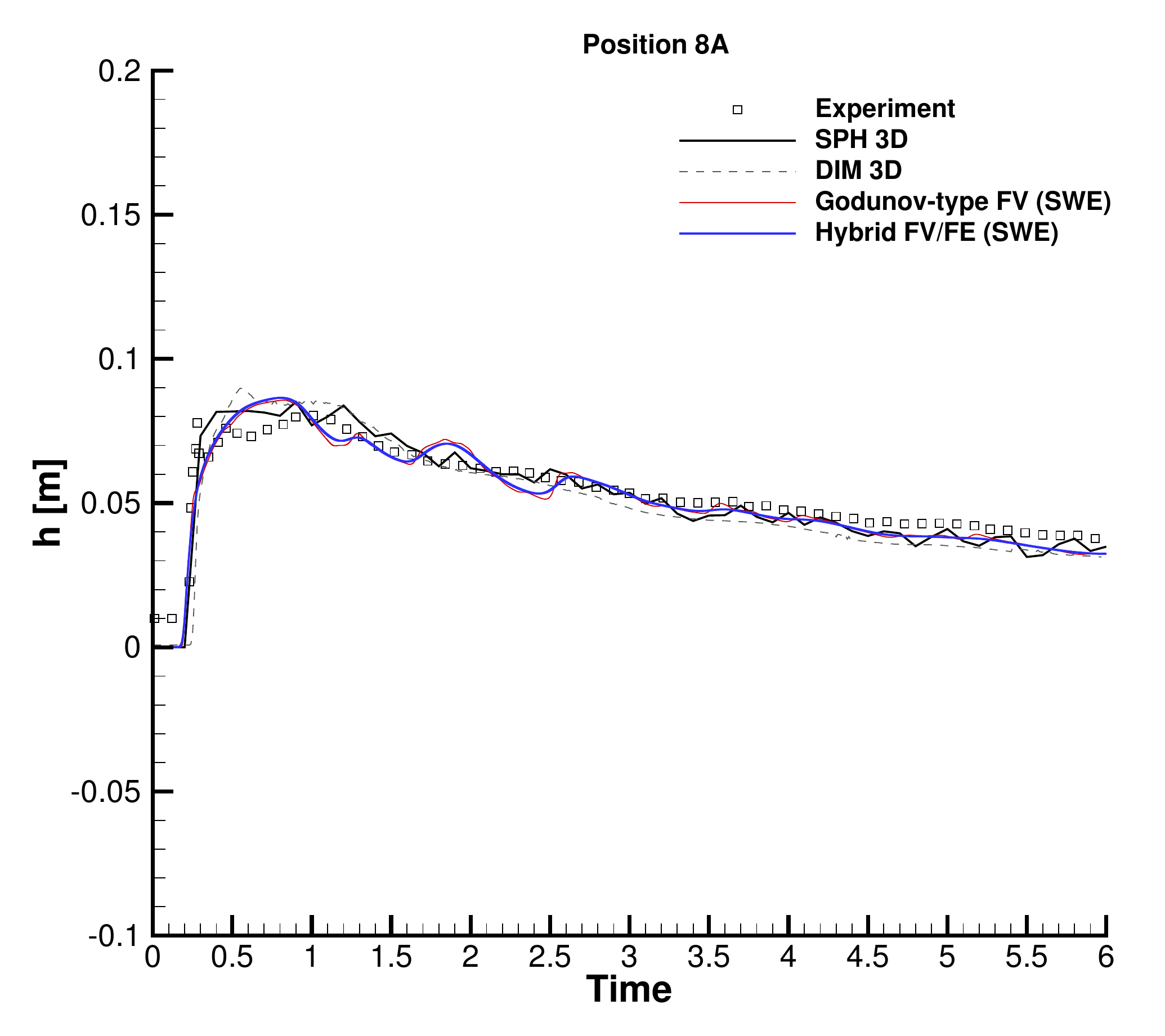}  
		\end{tabular}
		\caption{Time series of the free surface at wave gauges $0$, $-2A$, $-3A$, $-5A$, $1A$ and $8A$. The experimental results of Fraccarollo and Toro \cite{Dambreak3Dexp} are compared with those of an explicit Godunov-type finite volume scheme and the new hybrid FV/FE scheme applied to the shallow water equations. For further comparison, also the numerical results obtained with a fully nonhydrostatic 3D SPH scheme (SPH) and a 3D diffuse interface method (DIM) are shown.  } 
		\label{fig.dbf.h}
	\end{center}
\end{figure}

\begin{table}[h!]
	\renewcommand{\arraystretch}{1.1}
	\caption{Coordinates of the wave gauges of the three-dimensional dambreak for a bidimensional domain $\Omega=[-1,1]\times[-1,2]$.}
	\label{tab.wavegauge.dbf.3d}
	\begin{center}
		\begin{tabular}{ccccccc}
			\hline 
			Wave Gauge &   -5A & -3A & -2A  & 0 & 1A & 8A \\ \hline
			$x$ & $ -0.82 $ & $ -0.62 $ & $ -0.42 $ & $ 0 $  & $ 0.02 $ & $ 0.722 $\\
			$y$ & $ 0 $ & $ 0 $ & $ 0 $ & $ 0 $  & $ 0 $ & $ 0 $\\
			\hline 
		\end{tabular} 
	\end{center}
\end{table}
To better monitor the flow evolution, different wave gauges were located both inside the tank and on the flat plate. In particular, we focus on the positions indicated in Table \ref{tab.wavegauge.dbf.3d}. Figure \ref{fig.dbf.h} reports the time evolution obtained using the hybrid FV/FE methodology together with the experimental data reported in \cite{Dambreak3Dexp}. In the same figures we also report further numerical results obtained for this test problem: a three dimensional smooth particle hydrodynamics (SPH) scheme for the weakly compressible free surface Navier-Stokes equations \cite{Dambreak3D}, as well as a diffuse interface method (DIM) for the weakly compressible Euler equations, \cite{DIM3D}, and also a pure Godunov-type finite volume scheme for the shallow water equations, \cite{Dambreak3Dexp}. As expected, the solution provided by the shallow water solvers is less accurate at the beginning of the simulation, since the hydrostatic pressure assumption inherent in the classical shallow water equations is not verified. However, the matching with the different schemes and the experimental data improves at larger times, providing a good agreement with the reference data at a smaller computational cost than the fully nonhydrostatic 3D solvers. Further details on the description of the physical phenomena,  the comparison between 3D Navier-Stokes results and shallow water solvers and a description of the SPH and DIM methodologies can be found in \cite{Dambreak3D,DIM3D}. 


\section{Conclusions}\label{sec:conclusions}
In this paper, we have presented an extension of the  modern family of semi-implicit hybrid FV/FE methods \cite{BFTVC17,BBDFSVC20,BRVD21} to the solution of two-dimensional shallow water flows at all Froude numbers. A time semi-discretization of the PDE system written in conservative variables allows the decoupling of the computation of the conserved variables and the free surface elevation. The resulting method is provably well balanced and perfectly adapts to the hybrid finite volume/finite element framework developed previously for the compressible and incompressible Navier-Stokes equations. The computational domain is then discretized at the aid of staggered unstructured grids, thus allowing complex domain configurations that are relevant for environmental and coastal engineering. The proposed backward and forward interpolation between meshes is designed to avoid the generation of spurious numerical diffusion. The hyperbolic equations are solved using a second order in space and time finite volume methodology providing an intermediate approximation for the momentum vector in the dual grid.  Next, the second order wave equation for the free surface is addressed using classical continuous finite elements and the velocity field is finally corrected with the updated free surface elevation. To ameliorate the robustness of the developed methodology, avoiding spurious oscillations of the solution in the presence of large discontinuities, a Rusanov-type dissipation term has been added to the FE system. Despite the considered LADER methodology is second order accurate in space and time in a pure FV framework, the splitting employed in this paper decreases the temporal accuracy of the overall scheme. To improve it, we propose an extension of the algorithm making use of the theta method. 
The final algorithm is carefully validated including a numerical convergence analysis, the numerical verification of the C-property and the simulation of a set of Riemann problems. The flow over a cylinder and flow over a blunt body test cases show the capability of the method to deal with a wide range of Froude numbers. Finally, a 3D dambreak over a dry plane test is presented comparing the numerical solution with experimental data and other numerical reference solutions obtained with 3D codes for the free-surface Navier-Stokes equations and a pure Godunov-type finite volume discretization of the shallow water system.

Future work will concern the practical application of the developed methodology to more complex and more realistic problems, which will require an efficient parallel implementation of the proposed algorithm. {We will also consider the use of general polygonal and in particular of general quadrilateral grids, which are well suited to capture the behaviour of viscous flow in the vicinity of a wall.} Moreover, we also plan an the extension of the developed methodology to other hyperbolic PDE systems including involution constraints, such as the MHD equations \cite{DBTF19_semiimplicitMHD}, or the unified GPR model of continuum mechanics \cite{GPRmodel,SIGPR}.

%
\section*{Acknowledgements}

S.B. and M.D. are members of the INdAM GNCS group and acknowledge the financial support received from the Italian Ministry of Education, University and Research (MIUR) in the frame of the Departments of Excellence  Initiative 2018--2022 attributed to DICAM of the University of Trento (grant L. 232/2016) and in the frame of the 
PRIN 2017 project \textit{Innovative numerical methods for evolutionary partial differential equations and  applications}. S.B. was also funded by INdAM via a GNCS grant for young researchers and by an \textit{UniTN starting grant} of the University of Trento. 

The authors would like to thank Prof. L. Fraccarollo for kindly providing the experimental reference data of the 3D dambreak problem.

\bibliographystyle{plain}
\bibliography{./mibiblio}

\end{document}